\documentclass[12pt,a4paper]{article}
\usepackage{a4wide}
\usepackage{amsfonts}
\usepackage{amsmath}
\usepackage{amssymb, amsthm}
\usepackage{multirow}
\usepackage{latexsym}
\usepackage{graphicx}
\usepackage{color, epsfig} 

\newtheorem{defi}{Definition}[section]{\bf}{\rm}
\newtheorem{rema}{Remark}[section]{\bf}{\rm}
\newtheorem{theorem}{Theorem}[section]{\bf}{\rm}
\newtheorem{proposition}{Proposition}[section]{\bf}{\rm}
\newtheorem{corollary}{Corollary}[section]{\bf}{\rm}
\newtheorem{lemma}{Lemma}[section]{\bf}{\rm}
\newtheorem{acknowledgement}{Acknowledgement}[section]{\bf}{\rm}
\numberwithin{equation}{section} 
\newcommand{\auth}{\textsc} 

\newcommand \be   {\begin{equation}}
\newcommand \ee   {\end{equation}}
\newcommand \RR    {\mathbb{R}}  
\newcommand \Scal    {\mathcal{S}}
\newcommand \Ccal    {\mathcal{C}} 
\newcommand \Tcal    {\mathcal{T}} 
\newcommand \Ocal    {\mathcal{O}} 
\newcommand \Ncal    {\mathcal{N}} 
\newcommand \Ical    {\mathcal{I}} 
\newcommand \hpi {\widehat \pi} 

\newcommand \sbold   {\mathbf{s}}  
\newcommand \RN    {{\RR^N}} 
\newcommand \eps     {\epsilon} 
\newcommand \del     {\partial} 
\newcommand \lamb  {{\overline \lam}} 
\newcommand \Sigmab  {{\overline \Sigma}} 
\newcommand \Vb  {{\overline V}} 
\newcommand \lam   {\lambda}  
\newcommand \gam   {\gamma}  
\newcommand \hlam   {\lambda_{N+1}}

\newcommand \Lam   {\Lambda} 
  
\newcommand \Bzero {\RR^N}
\newcommand \conv {\text{conv}} 
\newcommand \Lip   {{\mathcal Lip}} 
\newcommand \dist  {\text{dist}} 
\def     \conv    {\text{conv}}
\newcommand \conc  {\text{conc}}

\newcommand \strength   {{\eps}} 
\newcommand \tildeu {{\widetilde u}}

\newcommand \Bcal    {\mathcal{B}}
\newcommand \Bone    {{\mathcal{B}_{\delta_1}}}
\newcommand \Btwo    {{\mathcal{B}_{\delta_2}}}
\newcommand \GO {{\mathcal O}} 
\newcommand \OO {{\mathcal O}}
\newcommand{\HH}{{\mathcal H}} 
\newcommand  \rt  {\widehat{r}} 
\newcommand{\WW}{{\mathcal W}} 
  
\newcommand \Hcal    {\mathcal{H}}
\newcommand \Gcal    {\mathcal{G}}
\newcommand \Fcal    {\mathcal{F}}

\newcommand \Jcal    {\mathcal{J}}
\newcommand{\Id}{\text{Id}}

\def\s{\mathbf s}

\newcommand \ONB {\overline{\mathcal N}}
\def\osc{\mathop{osc}}  

\newcommand \tr {\widetilde r} 
\newcommand \ts {\widetilde s} 
\newcommand \tf {\widetilde f} 
\newcommand \fhat {\widehat f}
\newcommand \tlam {\widetilde \lambda} 
\newcommand \hgam {\widehat \gam}  
\newcommand \supp {\text{supp}} 
\def\espaceverticalhaut{\parbox[b][0.5cm][c]{0cm}{\ }}
\def\espaceverticalhautb{\parbox[b][0.9cm][c]{0cm}{\ }}
\def\espacevertical{\parbox[][0.7cm][c]{0cm}{\ }}

\newcommand \tu {{\tilde{u}}}
\newcommand \lh {{\widehat l}} 
\newcommand \kappah {{\widehat \kappa}} 
\newcommand \fh {{\widehat f}} 
\newcommand \gh {{\widehat g}} 
\newcommand \tm {{\widetilde m}} 
\newcommand \QBB {Q} 
\newcommand \IcalBB {\Ical}
\begin{document}
\title{Nonlinear hyperbolic systems: 
\\  
Non-degenerate flux, inner speed variation, 
\\
and graph solutions}
\author{Olivier Glass and Philippe G. LeFloch \\
Laboratoire Jacques-Louis Lions\\
Centre National de la Recherche Scientifique, U.M.R. 7598\\
University of Paris 6\\
75252 Paris Cedex 05, France.\\
glass@ann.jussieu.fr, \, lefloch@ann.jussieu.fr}
\date{} 
\maketitle 
\begin{abstract} 
We study the Cauchy problem for
general, nonlinear, strictly hyperbolic systems of partial differential equations in one space variable. First,  
we re-visit the construction of the solution to the Riemann problem and 
introduce the notion of a {\sl nondegenerate} (ND) system. This is the optimal condition  
guaranteeing, as we show it, that the Riemann problem can be solved with {\sl finitely many} waves, only; 
we establish that the ND condition is generic in the sense of Baire (for the Whitney topology), so that 
any system can be approached by a ND system. 
Second, we introduce the concept of {\sl inner speed variation}
and we derive new interaction estimates on wave speeds. 
Third, we design a wave front tracking scheme and establish its strong convergence to the entropy
solution of the Cauchy problem; this provides a new existence proof 
as well as an approximation algorithm. As an application, we investigate the time-regularity of the 
graph solutions $(X,U)$ introduced by the second author, and propose a {\sl geometric version} of our scheme; 
in turn, the spatial component $X$ of a graph solution can be chosen to be continuous in both time 
and space, while its component $U$ is continuous in space and has bounded variation in time.  
\end{abstract} 
 
 \vfill
 
 To appear in : Archive for Rational Mechanics and Analysis
 
 \eject 
     

\section{Introduction}   

In this paper, we are interested in entropy solutions 
to general, nonlinear, strictly hyperbolic systems of partial differential equations 
in one space variable, and we investigate 
several important issues of the theory, especially the Riemann problem, the wave interaction estimates, and  
the graph solutions. 
In particular, a version of the wave front tracking scheme is introduced below, which extends schemes 
originally introduced by \auth{Dafermos} \cite{Dafermos1}, \auth{DiPerna}~\cite{DiPerna}, \auth{Bressan} \cite{Bressan1}, 
and \auth{Risebro} \cite{Risebro} for {\sl genuinely nonlinear} systems of conservation laws 
(see also \cite{BaitiJenssen}). The new version applies to {\sl general} systems that need not be genuinely nonlinear
nor in conservative form. 

Recall that the front tracking scheme is a variant of 
the random choice scheme introduced by \auth{Glimm} \cite{Glimm} for genuinely nonlinear, conservative systems. 
Recently, Glimm's scheme was extended to conservative systems with {\sl general} flux 
by \auth{Bianchini} \cite{Bianchini}, \auth{Iguchi and LeFloch} \cite{IguchiLeFloch}, and 
\auth{Liu and Yang} \cite{LiuYang}. 
The key contributions in the above works were the derivation of sharp estimates for the solution of the Riemann problem,
together with the introduction of suitable interaction functionals controling the total variation of solutions to the Cauchy problem. 
The analysis in \cite{Bianchini,IguchiLeFloch,LiuYang} took its roots in earlier works, 
pioneered by \auth{Liu} in 1981, on particular systems
(see~\auth{Liu}~\cite{Liu1a,Liu1b,Liu1c,Liu2}, \auth{LeFloch et al.} \cite{LeFloch2,HayesLeFloch2,LeFloch4},
\auth{Chern} \cite{Chern}, 
and \auth{Ancona and Marson} \cite{AnconaMarson1,AnconaMarson2,AnconaMarson3}), 
and
in the recent developments on the vanishing viscosity method by \auth{Bianchini and Bressan}~\cite{BianchiniBressan}.

 In comparaison to the Glimm scheme, implementing a front tracking scheme is more demanding 
 since, in addition to the interaction estimates and functionals controling the growth of {\sl wave strenths} 
(provided in \cite{Bianchini,IguchiLeFloch,LiuYang}),  
 one also needs a precise control of {\sl wave speeds.} 
In a front tracking method,  
when dealing with a general hyperbolic system and in order to guarantee that the limiting function
is the entropy solution of interest, it is often necessary (at interactions) to split a given front into several fronts 
(of smaller strength, say). Clearly, it is preferable to avoid to repeat this splitting step too often during the evolution 
--- for otherwise this could lead to a finite-time blow-up of the algorithm with the number of wave fronts or of interactions points
becoming infinite.  

We found it convenient to introduce wave fronts that are not just single shocks or 
single rarefactions but, rather, {\sl wave packets} 
consisting of (finitely or possibly infinitely) many shock and rarefaction waves propagating at the {\sl same} speed. 
To control the speed of a $j$-wave packet $(u_-,u_+)$ we introduce 
the notion of {\sl inner speed variation} $\vartheta_j(u_-,u_+)$, defined
as the largest minus the smallest wave speeds within the wave packet under consideration. 
One contribution in the present work is to provide {\sl sharp interaction estimates} on $\vartheta_j(u_-,u_+)$, 
via suitable convexity and wave interaction arguments.

The present paper also provides a significant 
generalization of the earlier analysis by \auth{Iguchi and LeFloch}~\cite{IguchiLeFloch}, 
who dealt with general systems approached by piecewise genuinely nonlinear (PGNL) ones.  
We introduce below the notion of a {\sl nondegenerate (ND)} system 
and we show that the ND condition is the optimal condition guaranteeing  
that the Riemann problem is solvable with {\sl finitely many} waves, only. 
Moreover, the nondegeneracy condition is shown to be fully generic in the sense of Baire (in the Whitney topology), 
so that {\sl any} system can be approached by a nondegenerate one.  
  
We also include below a discussion of hyperbolic systems in nonconservative form. Recall that, for such systems, 
the standard notion of a distributional solution does not apply, and it was recognized in the 90's 
by \auth{LeFloch}~\cite{LeFloch0,LeFloch1} and \auth{Dal~Maso, LeFloch, and Murat}~\cite{DLM}, 
that one must supplement the hyperbolic system with 
a {\sl prescribed} 
family of continuous paths, related with the ``interior structure'' of discontinuities in solutions and 
typically constructed from the set of all traveling wave solutions associated with a {\sl given} 
parabolic regularization. The Riemann problem was solved in~\cite{DLM} for nonconservative, genuinely 
nonlinear systems, and a suitable generalization to the Glimm scheme was proposed
 by \auth{LeFloch and Liu}~\cite{LeFlochLiu}. We observe in this paper that both techniques developed 
 in \cite{BianchiniBressan} and \cite{IguchiLeFloch} allow for generalizations to nonconservative systems. 

In summary, our results provide a new proof of the {\sl global existence} of the entropy solution
(in the class of functions with tame variation) to the Cauchy problem
associated with a strictly hyperbolic system, and 
represent an alternative approach as well as a generalization 
to the proofs given in \cite{Glimm,LeFlochLiu} (for genuinely nonlinear systems) 
and, more recently, in~\cite{Bianchini,BianchiniBressan,IguchiLeFloch,LiuYang}. 

On the other hand, recall that the {\sl uniqueness} of the entropy solution
was established earlier by \auth{Bressan and LeFloch} \cite{BressanLeFloch1} (for genuinely nonlinear, conservative systems), 
\auth{Baiti, LeFloch, and Piccoli}  \cite{BaitiLeFlochPiccoli} (for general conservative or nonconservative systems
and arbitrary jump relations), and 
\auth{Bianchini and Bressan} \cite{BianchiniBressan} (for vanishing viscosity limits based the identity viscosity matrix). 

Finally, we apply our front tracking scheme and 
investigate the time-regularity of the {\sl graph solutions}  $(X,U)$ 
introduced in \auth{LeFloch}~\cite{LeFloch5,LeFloch6}. 
We propose a {\sl geometric version} of our scheme and, in turn,
provide here a new, {\sl more regular} parametrization of the graph solutions, which 
is not only continuous in space but also continuous in time, 
except at countably many times where wave cancellations/interactions take place. 

A brief outline of this paper follows. Section~\ref{Riem-0} deals principally with the Riemann problem
and generalizes the approach developed by Iguchi and LeFloch's; in particular, 
it includes the discussion of the nondegeneracy condition, wave interaction estimates, and 
inner speed variation estimates. Section~\ref{Sec:BB} provides 
a second proof of the inner speed variation estimates which follows 
the approach developed by Bianchini and Bressan. 
Section~\ref{NC-0} contains a brief discussion of generalizations of our results to nonconservative systems. 
Next, Section~\ref{WF-0} provides our new version of the front tracking scheme, together with the convergence proof. 
Finally, Section~\ref{SY-0} shows how to apply the previous framework in order to investigate the time-regularity of graph solutions.

 
\section{Non-degenerate hyperbolic systems of conservation laws}
\label{Riem-0}
 
\subsection{Notation}  

In the present section and in the following one we are primarily interested in {\sl conservative} systems.  
However, as we will explain in Section~\ref{NC-0} many of the forthcoming arguments carry over to 
{\sl nonconservative} systems, with only minor modifications. 
We thus consider a strictly hyperbolic system of conservation laws, 
\be
\del_t u + \del_x f(u) = 0, \qquad u=u(t,x) \in \RN, \quad  t \geq 0, \, x \in \RR,
\label{Class.1}
\ee
where, as is usual, we assume that all solutions take values in 
a neighborhood of some constant state in $\RN$ --~which is normalized to be the origin. 
We denote by $\Bcal_\delta$ the open ball in $\RN$ centered at the origin and having radius $\delta>0$, 
and we assume that all values $u(t,x) \in {\overline\Bcal}_{\delta_1}$ for some $\delta_1>0$.  
In \eqref{Class.1}, the map ${f : \Bone \to \RN}$ is smooth  
and, for each $u \in \Bone$, $A(u):=Df(u)$ admits $N$ real and distinct eigenvalues $\lam_1(u) < \cdots < \lam_N(u)$.
We denote by $r_j(u)$ and $l_j(u)$ the left- and right-eigenvectors associated with $\lam_j(u)$, i.e., 
$$ 
A(u) \, r_j(u) = \lam_j(u) \, r_j(u), 
\quad l_j(u) \cdot A(u) = \lam_j(u) \, l_j(u),  
$$
and normalized so that 
$$
\aligned
&  |r_j(u)| = 1, \quad l_j(u) \cdot r_j(u) =1, 
\\
&  l_i(u) \cdot r_j(u) = 0 \quad \mbox{ for } i \neq j.  
\endaligned
$$

In \cite{IguchiLeFloch}, the authors advocated to use, for each $j$-wave family, a foliation of $\Bone$ based on  
a {\bf global parameter} $u \mapsto \mu_j(u)$, satisfying 
$$
\nabla\mu_j(u) \cdot r_j(u) \neq 0.  
$$
Each wave curve ${s \mapsto \psi_j(m)=\psi_j(m;u_0)}$ is parametrized to ensure that the state $\psi_j(m)$ 
lies on the submanifold $\big\{ \mu_j(u)= m \big\}$, i.e., 
$$ 
\mu_j(\psi_j(m)) = m. 
$$  
In the present paper we adopt a particular global parametrization.  
We fix a constant vector $\lh$ such that $\lh \cdot r_j(u)$ never vanishes, 
$$
\lh \cdot r_j(u) \neq 0, 
$$ 
and we set 
$$
\mu_j(u) = \mu(u) := \lh \cdot u. 
$$
Hence, our global parameter is also {\sl independent} of the specific wave family under consideration. 
The following formulas will involve the renormalized eigenvectors 
$$
\rt_j := \frac{r_j }{ \nabla\mu_j \cdot r_j} = \frac{r_j}{ \lh  \cdot r_j}. 
$$

Recall that, for the {\bf $j$-integral curve}  
$$
\OO_j(u_-) = \bigl\{ w_j(m;u_-) \, /\, m_\star \leq m \leq m^\star \big\}
$$ 
defined by 
$$
\del_m w_j = \rt_j(w_j), \quad w_j(m_-;u_-) = u_-,
$$
where $m_-:=\mu(u_-)$, 
one can check that there exists $\delta_2 < \delta_1$ such that the following property holds. 
For all $u_- \in \Btwo$, the map 
$[m_\star, m^\star] \ni m \mapsto w_j(m_-;u_-)$ is smooth and takes  
its values in $\Btwo$, and the end points $m_*, m^*$ (depending on $u_-$) 
satisfy  $w_j(m_*;u_-), w_j(m^*;u_-)  \in \del\Btwo$. 
Moreover, one has 
$$
w_j(m;u_0) = u_0 + (m - m_0) \, \rt_j(u_0)  + \GO(m - m_0)^2
$$
with $m_0:=\mu(u_0)$. 

Similarly, the {\bf $j$-Hugoniot curve}  
$$ 
\HH_j(u_-) := \big\{v_j(m;u_-) \, / \, m_\star \leq m \leq m^\star \big\},
$$  
defined by 
$$
- \lamb_j(u_-,v_j(m;u_-)) \, \big(v_j(m;u_-) - u_- \big) + f(v_j(m;u_-)) - f(u_-) = 0, 
$$
satisfies  
$$
v_j(m;u_0) =  u_0 + (m - m_0) \, \rt_j(u_0) + \GO(m - m_0)^2,  
$$
whereas the {\bf $j$-shock speed} $\lamb_j(m;u_-) := \lamb_j(u_-, v_j(m;u_-))$ satisfies 
$$
\lamb_j(m;u_0) 
= \lam_j(u_0) +  \frac{1}{2} \bigl( \nabla\lam_j \cdot \rt_j \bigr)(u_0) \, (m -m_0) + \GO(m - m_0)^2.
$$

As usual, we restrict attention to propagating discontinuities $(u_-,u_+)$ with $u_+ \in \HH_j(u_-)$, 
satisfying the following {\bf entropy criterion} \cite{Oleinik,Wendroff,Liu1a,Liu1b,Liu1c,Liu2}. 
Setting $m_-:= \mu(u_-)$, $m_+:= \mu(u_+)$, and
$$ 
u_+ = v_j(m_+;u_-),  
$$ 
the shock $(u_-,u_+)$ is said to be {\bf admissible} if and only if   
$$
\lamb_j(u_-,u_+) \leq \lamb_j(m;u_-) \quad \mbox{ for all $m$ between $m_-$ and $m_+$.} 
$$
In other words, the shock speed achieves its minimum value at the point $m_+$. 
Following \cite{Lax1,Liu2,IguchiLeFloch}, we will be interested in 
combining rarefaction curves and shock curves and constructing
the {\bf $j$-wave curve} issuing from a given left-hand state $u_0 \in \Btwo$,
$$
\WW_j(u_0) := \bigl\{ \psi_j(m;u_0) \, / \, m_{\star} \leq m \leq m^\star \bigr\},
$$
where the end points of the curve satisfy  
$$
m_\star = m_{j\star}(u_0) < \mu(u_0)  < m_j^\star(u_0)=m^\star. 
$$


\subsection{Riemann problem for nondegenerate systems}

Recall that, following Liu's pioneering contribution to the subject (\auth{Liu} \cite{Liu1a,Liu1b,Liu1c,Liu2}),
\auth{Iguchi and LeFloch}~\cite{IguchiLeFloch} constructed and investigated the regularity of 
the entropy solution to the Riemann problem, 
for systems with piecewise genuinely nonlinear (PGNL) flux (see below for a definition). 
Here, we introduce a larger class and we establish that all of the results in \cite{IguchiLeFloch} remain valid. 
Recall that all of the key estimates depend on the $C^2$ norm of the flux only, 
and allowed the authors in \cite{IguchiLeFloch} to cover any flux realized as a limit of PGNL functions. 
Our more general condition will be shown to be fully {\sl generic},   
so that, by density, we can cover systems with fully arbitrary flux. 

It will be convenient to introduce the notation ($1 \leq j \leq N$, ${u \in \Bone}$) 
$$
\aligned 
& \pi_j^{(1)}(u) :=  r_j(u) \cdot \nabla \lam_j(u), 
\\
&  \pi_j^{(k+1)}(u) : = r_j(u) \cdot \nabla \pi_j^{(k)}(u), \qquad k = 1,2, \ldots, 
\endaligned 
$$
together with the following definition. Note that the definition clearly makes sense for a general 
matrix-valued mapping $A$ which need not be a Jacobian matrix. 

\begin{defi} 
\label{Riem-1} 
The matrix-valued mapping $A =A(u)$ ($u \in \Bone$) is called {\bf nondegenerate} (ND, in short) 
if for every ${j=1, \ldots, N}$ and every $u \in \Bone$, not all of the values $\pi_j^{(k)}(u)$
(${1 \leq k \leq N+1}$) are zero,
$$
\big( \pi_j^{(1)}(u), \pi_j^{(2)}(u), \ldots, \pi_j^{(N+1)}(u) \big) \neq \big( 0,0, \ldots, 0 \bigr). 
$$
\end{defi}

When the ND condition holds, we will also say that \eqref{Class.1} is a {\sl nondegenerate system}
or that, in the conservative case, $f$ is a {\sl nondegenerate flux.}

For instance, in the scalar case $N=1$ a function $f:\RR \to \RR$ is nondegenerate 
if and only if $f''$ and $f'''$ do not vanish simultaneously. 

Given any nondegenerate flux $f$  
we can associate to each point $u \in \Bone$ its {\bf critical exponent} $p=p_j(u)$ as 
the smallest index $k \geq 1$ such that $\pi_j^{(k)}(u) \neq 0$. 
Our definition includes the genuinely nonlinear flux for which $\pi_j^{(1)}(u)$ never vanishes
(so that $p_j(u)=1$), 
as
well as the {\bf piecewise genuinely nonlinear} (PGNL) flux introduced in \cite{IguchiLeFloch} 
where it was assumed that $\pi_j^{(1)}(u), \pi_j^{(2)}(u)$ never 
vanish simultaneously (so that $p_j(u)$ equals $1$ or $2$). 
The case of linearly degenerate fields for which $\pi_j^{(1)}(u)$ vanishes identically 
could be easily included in the present discussion, but for simplicity in the presentation we prefer to 
cover this case latter via a general density argument. 

Our first objective is to generalize the construction in \cite{IguchiLeFloch} and to 
prove that the Riemann problem for ND flux can be solved uniquely and that the solution
contains finitely many waves, only. 

\begin{theorem}{\rm (Riemann problem for nondegenerate flux.)}
\label{Riem-2}

\noindent
Let \eqref{Class.1} be a strictly hyperbolic system of conservation laws with nondegenerate flux $f:\Bone \to \RN$. 
Then there exists $\delta_2< \delta_1$ depending only upon the $C^2$ norm of $f$, such that 
the following properties holds.

1. For all $u_0 \in \Btwo$ and all ${j=1,\ldots,N}$ there exist $m_\star= m_{j\star}(u_0)$ 
and $m^\star= m_j^\star(u_0)$  
and a mapping $\psi_j=\psi_j(m;u_0)$ defined for $m_\star \leq m \leq m^\star$ 
and Lipschitz continuous in both arguments, such that  
$$
\aligned 
& \psi_j(m;u_0) \in \Btwo, \quad m_\star < m < m^\star, 
\\
& \psi_j(\mu(u_0);u_0) = u_0, 
\\
& \psi_j(m_\star;u_0), \, \psi_j(m^\star;u_0)  \in \del \Btwo. 
\endaligned 
$$

2. For each $m \in [m_\star,m^\star]$ there exist an integer $q$ 
and a finite sequence of states $u_1,u_2,\ldots,$ $u_{2q+1} \in \Btwo$ such that 
$$
u_{2q+1} = \psi_j(m;u_0) 
$$
and $($for all relevant values of $k$) $u_{2k}$ is connected to $u_{2k+1}$ by 
a $j$-rarefaction wave while $u_{2k+1}$ is connected to $u_{2k+2}$ by an admissible $j$-shock, 
with 
$$
\aligned  
\lam_j(u_0) \leq \lam_j(u_1) \leq & \lamb_j(u_1, u_2) 
\leq \lam_j(u_2) \leq \cdots  
\\
& \cdots \leq \lamb_j(u_{2q-1},u_{2q}) \leq \lam_j(u_{2q}) \leq \lam_j(u_{2q+1}).
\endaligned 
$$

3. For any Riemann data $u_l,u_r \in \Bcal_{\delta_2}$ the Riemann problem \eqref{Class.1},
\be
u(0,x) = \begin{cases}
u_l, & x<0, 
\\
u_r, & x>0, 
\end{cases} 
\nonumber 
\ee 
admits a unique, self-similar solution made of finitely many rarefaction waves and admissible shock waves. 
\end{theorem}

To establish these results we will explain how 
the arguments in \cite{IguchiLeFloch} can be extended to cover nondegenerate flux.  
Certain monotonicity properties of the shock speed along the Hugoniot curve and its critical values
play a central role in the construction of  the wave curves, and this will be the subject of 
the following section.  

Once the Riemann problem is solved and provided finitely many waves only arise in the construction, 
one can follow \cite{IguchiLeFloch} and derive additional regularity and interaction estimates.
Indeed, we emphasize that the proofs therein rely only on the facts 
that the Riemann solution contain finitely many waves and that all Hugoniot and integral curves
associated with \eqref{Class.1} are smooth maps. We summarize here the results that we will need later.

\begin{proposition} {\rm (Regularity of the wave curves.)}
\label{regu} 
Under the assumptions of Theorem~\ref{Riem-2}, for each $j= 1, \ldots, N$ the map 
$\psi_j= \psi_j(m;u)$ is Lipschitz continuous with respect to both $m,u$ (with Lipschitz constant depending only on the 
$C^2$ norm of $f$). Moreover, the first-order derivatives of $\psi_j$ with respect to $m,u$ 
are Lipschitz continuous at the point $m=\mu(u)$, that is 
\be
\begin{split} 
& \psi_j'(m;u) = \rt_j(u) + \GO(m-\mu(u)), 
\\
& D_u\psi_j(m;u) = \Id - \rt_j(u) \otimes \lh + \GO(m-\mu(u)),  
\end{split}
\nonumber 
\ee 
where $\otimes$ denotes the tensor product of two vectors.  
\end{proposition}

To any left-hand state $u_-$ and a right-hand state $u_+=\psi_j(m_+; u_-)$, we associate  
the {\bf wave strength} 
$$
\strength_j(u_-,u_+):=\mu(u_-) - \mu(u_-)
$$ 
and the associated {\bf wave speed function} $\lamb_j(s, m_+;u_-)$, which 
by definition yields the speed of waves making up the wave fan.  
Given two $j$-waves $(u_-,u_+)$ and 
$(u_-',u_+')$ we introduce the {\bf generalized angle} between them, as follows  
$$
\theta_j(u_-,u_+; u_-',u_+') 
:= 
\, \int_{m_-}^{m_+} \hskip-.2cm \int_{m_-'}^{m_+'} 
\frac{\big( \lamb_j(s', m_+'; u_-') - \lamb_j(s, m_+; u_-)  \big)_- }{(m_+ - m_-) \, (m_+' - m_-') } \, ds ds', 
$$
where $m_\pm := \mu(u_\pm)$ and $m_\pm' := \mu(u_\pm')$.

\begin{proposition} {\rm (Interaction estimates.)} 
\label{interactions} 
For all $u_l$, $u_m$, and $u_r \in \Bcal_{\delta_2}$ and $1 \leq i,j \leq N$
the following property holds. Suppose that $u_l$ is connected to $u_m$ by 
an $i$-wave fan and that $u_m$ is connected to $u_r$ by a $j$-wave fan. 
Then, the wave strengths $\strength_k(u_l, u_r)$ of the outgoing Riemann solution connecting $u_l$ to $u_r$ satisfy
($1 \leq k \leq N$)  
$$
\strength_k(u_l, u_r) 
= \strength_k(u_l, u_m) + \strength_k(u_m, u_r) + \GO(1) \, Q_{}(u_l, u_m, u_r), 
$$
where Iguchi-LeFloch's {\bf interaction potential} is defined by 
$$
Q_{}(u_l, u_m, u_r) : = \sum_{i \geq j} \Theta_{ij}(u_l, u_m, u_r) \, |\sigma_i(u_l, u_m) \, \sigma_j(u_m, u_r)| 
$$
and 
\be
\Theta_{ij}(u_l, u_m, u_r) 
:= \begin{cases}
  0,  & i < j, 
  \\
  1,  & i > j, 
  \\
  1,  & i=j, \, \sigma_j(u_l, u_m) \, \sigma_j(u_m, u_r) \leq 0, 
  \\
  \theta_j(u_l,u_m; u_m,u_r), & i=j, \,  \sigma_j(u_l, u_m) \, \sigma_j(u_m, u_r) > 0.
 \end{cases}
\nonumber 
\ee
\end{proposition}

In the rest of this section, we derive a key estimate on wave speeds. Let us introduce the following notion.

\begin{defi} The {\bf inner speed variation} of a $j$-wave connecting a left-hand state $u_-$ 
to a right-hand state $u_+=\psi_j(m_+; u_-)$, is defined as
$$
\aligned 
\vartheta_j(u_-, u_+) 
& = \vartheta_j(m_+; u_-)
\\
& : = \lamb_j^{\max}(u_-, u_+) - \lamb_j^{\min}(u_-,u_+)
\\ 
& : = 
\lamb_j(m_+,m_+;u_-) - \lamb_j(m_-, m_+;u_-).
\endaligned 
$$ 
\end{defi} 

For instance, $\vartheta_j(u_-, u_+) =0$ if and only if $u_-$ and $u_+$ are connected by a single admissible
discontinuity, while 
 $\vartheta_j(u_-, u_+) = \lam_j(u_+) - \lam_j(u_-)$ if (but not only if) $u_-$ and $u_+$ are connected by a single 
 rarefaction wave.

\begin{theorem} {\rm (Properties of the inner speed variation.)} 
\label{PropISV}

1. For any $u_-, u_-',m$ it holds  
$$
\vartheta_j(\mu(u_-')+m; u_-')  = \vartheta_j(\mu(u_-) + m; u_-) + \OO(m) \, |u_-' - u_-|. 
$$

2. For any $u_l, u_m,u_r$ with $u_m = \psi_j(\mu_m; u_l)$ and  $u_r = \psi_j(\mu_r; u_m)$, 
setting $u_r' := \psi_j(\mu_r; u_l)$,
$$
\vartheta_j(u_l, u_r') 
= \vartheta_j(\mu_r; u_l)  = \vartheta_j(\mu_m; u_l) + \GO( |\mu_r - \mu_m|). 
$$

3. For any $u_l, u_m,u_r$ with $u_m = \psi_j(\mu_m; u_l)$ and  $u_r = \psi_j(\mu_r; u_m)$ 
and monotonically ordered along the wave curve, 
\be
 \begin{split} 
\vartheta_j(\mu_r; u_l)  
\leq 
& \max\big( \vartheta_j(u_l,u_m), \vartheta_j(u_m, u_r) \bigr) 
\\
& + \Big( \lamb_j^{\min}(u_m,u_r) - \lamb_j^{\min}(u_l,u_m) \Big)_+ 
\\
& + \GO(1) \, Q_{}(u_l, u_m, u_r). 
\end{split} 
\label{56}
\ee
\end{theorem}

In the application of Theorem~\ref{PropISV} in the context of the front tracking scheme, 
every front $(u_-,u_+)$ will propagate with the speed $\lamb_j^{\min}(u_-,u_+)$, 
so that the 
term $\Big( \lamb_j^{\min}(u_m,u_r) - \lamb_j^{\min}(u_l,u_m) \Big)_+$ 
vanishes if the fronts $(u_l,u_m)$ and $(u_m,u_r)$ interact. 

Finally, after establishing a density result of the class of non-degenerate flux (Theorem~\ref{ApproxFlux} below)
we can conclude that:

\begin{corollary} 
\label{lecorollaire}
All of the results stated in the present section remain valid for general, strictly hyperbolic flux
that need not be nondegenerate, with the modification that a Riemann solution may 
contain (not finitely many but) countably many waves. 
\end{corollary}


\subsection{Fundamental properties at critical points} 
 
From now it will convenient to work with the following variant of the maps $\pi_j$, obtained by replacing 
the vectors $r_j$ by the normalized vectors $\rt_j$, i.e. 
$$
\aligned 
& \hpi_j^{(1)}(u) :=  \rt_j(u) \cdot \nabla \lam_j(u), 
\\
&  \hpi_j^{(k+1)}(u) : = \rt_j(u) \cdot \nabla \hpi_j^{(k)}(u), \qquad k = 1,2, \ldots, 
\endaligned 
$$
Clearly, for every $u$ and every $j$, the condition 
$$
\pi_j^{(1)}(u) = \ldots = \pi_j^{(k)}(u)=0, \qquad \pi_j^{(k+1)}(u) \neq 0
$$
is equivalent to 
$$
\hpi_j^{(1)}(u) = \ldots = \hpi_j^{(k)}(u)=0, \qquad \hpi_j^{(k+1)}(u) \neq 0, 
$$
so that a statement involving the coefficients $\hpi_j$ can be immediately restated with the coefficients $\pi_j$.

It is a well known fact that rarefaction curves and Hugoniot locus have a second-order 
tangency property at their base point. Here, we prove a new, higher-order tangency property which is satisfied 
at critical points. Recall that by differentiating the Rankine-Hugoniot it follows that 
\be
(m-m_-) \, \del_m \lamb_j(m;u_-) =  \kappa_j(m;u_-) \, \big(\lam_j(v_j(m;u_-)) - \lamb_j(m;u_-) \big),
\label{Pro.1}
\ee  
where the function $\kappa_j =\kappa_j(m;u_-)>0$ is smooth, bounded, and bounded away from zero, and satisfies 
$\kappa_j(m_-;u_-)=1$.

\begin{lemma} {\rm (Tangency property at critical points (I).)}  
\label{Riem-Tan} 
Suppose that, at some point $u_\star \in \Bone$ and for some $p \geq 1$,
\be
\label{Degen}
\aligned 
& \hpi_j^{(k)}(u_\star) = 0, \quad k=1, \ldots, p-1, 
\\
& \hpi_j^{(p)}(u_\star) \neq 0. 
\endaligned 
\ee
Then, the shock curve $v_j$ and the rarefaction curve $w_j$ issuing from $u_\star$ are tangent at $u_\star$ 
up to the order $p+1$, 
\be
\label{Tangence}
\del_m^{(k)}v_j(m_\star;u_\star) = \del_m^{(k)}w_j(m_\star;u_\star), \quad k=1, \ldots, p+1,
\ee
while the shock speed satisfies
\be
\label{TanVitChoc}
\aligned 
& \del_m^{(k)}\lamb_j(m_\star;u_\star) = 0, \quad k=1, \ldots, p-1, 
\\
& \del_m^{(p)}\lamb_j(m_\star;u_\star) = \frac{\hpi_j^{(p)}(u_\star) }{p+1}.
\endaligned 
\ee 
\end{lemma}

\begin{proof} Proceeding by induction on $k$ we first note that, for {\sl all} $m$,
\be
\aligned 
\del_m^{(k+1)} \lam_j(w_j(m;u_\star)) 
& = \del_m \hpi_j^{(k)}(w_j(m;u_\star)) 
\\
& = \nabla \hpi_j^{(k)}(w_j(m;u_\star)) \cdot \del_m w_j(m;u_\star),  
\\
& = \hpi_j^{(k+1)}(w_j(m;u_\star)). 
\endaligned 
\label{VitCar}
\ee

Now, \eqref{Tangence} and \eqref{TanVitChoc} are proven by induction: we consider the induction hypothesis at 
the rank $q$ (for $0 \leq q \leq p$):
\begin{quote}
$\del_m^{(k)} \lamb_j(m^\star,u^\star)=0$ holds for $k=1,\ldots, q-1$ and
$\del_m^{(k)}(v-w)(m^\star,u^\star)=0$ holds for $k=0,\ldots, q$.
\end{quote}
This is clearly satisfied at rank $q=0$. Suppose this is true at rank $q \leq p-1$;  let us establish it at rank $q+1$.

For the claim on $\del_m^{(q)} \lamb_j$ we proceed as follows.
In view of \eqref{Pro.1} and after differentiation with respect to $m$ we obtain
\be 
\del_m^{(k)} \lamb_j(m;u_\star) + \del_m^{(k)} \Big( \frac{m-m_\star}{\kappa_j(m;u_\star)} \, \del_m \lamb_j(m;u_\star) \Big)
= \del_m^{(k)} \lam_j(v_j(m;u_\star)). 
\label{Pro.8}
\ee  
Considering \eqref{Pro.8} at the point $m=m_\star$, 
using \eqref{VitCar} and the induction hypothesis,  
it follows easily that $\del_m^{(k)}\lamb_j(m_\star,u_\star) = 0$ for  $k=q$ with $q<p$.

We now consider the second claim about $\del_m^{(q+1)}(v(m^\star,u^\star)-w(m^\star,u^\star))=0$. 
By differentiating the Rankine-Hugoniot relation, we find 
$$
\del_m\lamb_j(m) \, (v_j(m)-u_-) = \bigl( A(v_j(m)) - \lamb_j(m)\bigr) \, \del_m v_j(m),
$$
which is to be compared with the equation characterizing the integral curve $w_j(m) = w_j(m;u_0)$, 
$$
0 = \bigl( A(w_j(m)) - \lam_j(w_j(m))\bigr) \, \del_m w_j(m). 
$$
The map $z(m) := v_j(m; u_-) - w_j(m; u_0)$ satisfies $z(m_0)=0$. From the above identities
we deduce that 
\be
\begin{split} 
& \bigl( A(w_j) - \lam_j(w_j)\bigr) \, \del_m z 
\\
& = \del_m\lamb_j \, (v_j - u_-) + \bigl(  A(w_j) -  A(v_j) - \lam_j(w_j) + \lamb_j \bigr) \, \del_m v_j.
\label{EqDiffVW}
\end{split}
\ee

Using \eqref{Degen}, \eqref{VitCar}, and the induction hypothesis on $\lamb_j$ 
(now proven up to the $q$-th order derivative) and by differentiating in $m$, the equation \eqref{EqDiffVW} yields 
$$
\bigl( A(u_0) - \lam_j(u_0)\bigr) \, \del_m^{(q+1)} z(m_\star,u_\star) = 0.
$$
(In the case $q=0$, no assumption on $\del_m^{(k)} \lamb_j(m_\star,u_\star)$ is needed, but 
only $\lamb_j(m_\star;u_\star) = \lam_j(u_\star)$.) We deduce that, for some $\beta^{(q+1)} \in \RR$, 
we have $\del_m^{(q+1)} z(m_\star,u_\star) = \beta^{(q+1)} r_j(u_\star)$.
Now, the parameter along the curves is such that 
\be
\label{RappelParam}
\mu(w_j(m)) = m = \mu(v_j(m)). 
\ee
Differentiating this relation and using, again, the induction hypothesis we obtain 
$$
\nabla \mu(u_\star) \, \del_m^{(q+1)} w_j(m_\star) = \nabla \mu(u_\star) \, \del_m^{(q+1)} v_j(m_\star), 
$$ 
which yields $\del_m^{(q+1)} z (m_\star,u_\star) = 0$, and completes the induction.

It remains to prove \eqref{TanVitChoc} at the rank $p$ and \eqref{Tangence} at the rank $p+1$.
Using Leibniz identity and retaining only the non-zero terms, we deduce from \eqref{Pro.8}:
$$
\del_m^{(p)} \lamb_j(m;u_\star) + p \, \frac{\del_m^{(p)}\lamb_j(m_\star;u_\star)}{\kappa_j(m;u_\star)} 
= \del_m^{(p)} \lam_j(v_j(m_\star;u_\star)) = \hpi_j^{(p)}(u_\star), 
$$ 
which establishes the first claim. For the second claim, 
differentiating \eqref{EqDiffVW} $p$ times, using Leibniz identity, 
and keeping the relevant terms only, we obtain  
$$
\aligned 
& \bigl( A(w_j) - \lam_j(w_j)\bigr) \, \del^{p+1}_m z (m_\star,u_\star) 
\\
& = \Big( p \, \del^{p}_m\lamb_j (m_\star,u_\star) \, - \del_m^p \lam_j(w_j(m_\star,u_\star)) 
 + \del_m^p \lamb_j(m_\star,u_\star) \Big) \rt_j (u_\star) =0,
\endaligned 
$$ 
in view of \eqref{TanVitChoc}. Using again \eqref{RappelParam}, this establishes \eqref{Tangence}. 
\quad \end{proof}


It follows from \eqref{Pro.1} that, at a critical point $m_0$ where 
$$
\del_m \lamb_j(m_0;u_-)=0,
$$
the shock speed must coincide with the characteristic speed, 
$$
\lamb_j(m_0;u_-) = \lam_j(v_j(m_0;u_-)),
$$
and it can also be checked from the Rankine-Hugoniot relation that 
\be
\del_m v_j(m_0;u_-) = \rt_j(v_j(m_0;u_-)). 
\label{Pro.3}
\ee  
Generalizing this observation we now prove: 

\begin{lemma} {\rm (Tangency property at critical points (II).)}   
\label{Riem-3} 
If, for some $u_-, m_0$, and $p$,  
\be
\aligned
& \del_m^{(k)} \lamb_j(m_0; u_-)=0, \quad k=1, \ldots, p,
\\
& \del_m^{(p+1)} \lamb_j(m_0; u_-) \neq 0,  
\endaligned 
\label{Pro.4}
\ee  
then at the critical point $u_0:=v_j(m_0; u_-)$ one has 
$$
\aligned 
& \hpi_j^{(k)}(u_0)=0, \quad k=1, \ldots, p-1, 
\\
& \hpi^{(p)}_j(u_0) \neq 0, 
\endaligned 
$$ 
\be 
(m_0 - m_-)  \, \del_m^{(p+1)} \lamb_j(m_0;u_-) = \kappa_j(m_0;u_-) \, \hpi^{(p)}_j(u_0), 
\label{Pro.6}
\ee  
and the Hugoniot curve issuing from $u_-$ is tangent up to order $p$ to the integral curve issuing from $u_0$, 
that is 
\be 
\del_m^{(k)} v_j(m_0;u_-) = \del_m^{(k)}w_j(m_0;u_0), \quad k=0, \ldots, p. 
\label{Pro.7}
\ee  
\end{lemma}

\begin{proof} We first  rely on \eqref{Pro.8}.
By the assumption \eqref{Pro.4} the left-hand side vanishes at $m=m_0$ for all $k<p$. Hence, 
we deduce that, for all $k<p$, one has $\del_m^{(k)} \lam_{j}(v_j)(m_0;u_-) =0$. With \eqref{Pro.3}
we then deduce
$$
\del_m \lam_j(v_j(m_0;u_-)) = \nabla\lam_j(v_j(m_0;u_-)) \, \del_m v_j(m_0;u_-) = \nabla\lam_j \cdot \rt_j(v_j(m_0;u_-)),
$$
Hence, $\hpi^{(1)}_j(v_j(m_0; u_-))=0$. Hence, 
using \eqref{EqDiffVW} and \eqref{RappelParam}, one sees, as in Lemma~\ref{Riem-Tan} that \eqref{Pro.7} is valid for $p=2$.

More generally, if \eqref{Pro.7} is already established at the rank $k$ then using \eqref{VitCar} it follows that 
$$
\aligned 
\del_m^{(k)} \lam_j(v_j(m_0;u_-)) 
& = \del_m^{(k)} \lam_j(w_j(m_0;u_0)) 
\\
& =\hpi_j^{(k)}(u_0). 
\endaligned 
$$
Hence, we conclude that $\hpi_j^{(k)}(u_0)=0$. Using the same procedure as above, one deduces \eqref{Pro.7} at the rank $k+1$.

For $k=p$ the left-hand side of \eqref{Pro.8} at $m_0$ is a multiple of the left-hand side of \eqref{Pro.6}, 
$$
\frac{m_0 - m_- }{\kappa_j(m_0;u_-)} \, \del_m^{(p+1)} \lamb_j(m_0;u_-), 
$$
while the right-hand side equals $\hpi_j^{(p)}(u_0)$. Thus, \eqref{Pro.6} holds. 
\quad  \end{proof} 


Some further observations are in order. 

\begin{lemma} {\rm (Propagating discontinuities with coinciding speeds.)} 
\label{Riem-4} 
Suppose that $(u_0,u_1)$ and $(u_1,u_2)$ are shock waves satisfying the entropy criterion and 
propagating with the same speed $\Lam$, 
$$
\aligned 
& u_1 = v_j(m_1;u_0), \quad u_2 = v_j(m_2;u_1), 
\\ 
& \lamb_j(u_0,u_1) = \lamb_j(u_1,u_2) = \Lam, 
\endaligned 
$$
where $m_1$ and $m_2$ satisfy $\mu(u_0) < m_1 < m_2$. Then, the discontinuity $(u_0,u_2)$ is a shock satisfying the 
entropy criterion and propagating with the speed $\Lam$,  
$$
u_2 = v_j(m_2;u_0), \quad \lamb_j(u_0,u_2) = \Lam.
$$ 
\end{lemma}

\begin{lemma} {\rm (Equivalent formulation of the entropy criterion.)}
\label{Riem-5} 
A discontinuity $(u_-, u_+)$ with $u_+ \in \HH_j(u_-)$ satisfies the entropy criterion if and only if 
$$
\lamb_j(u_-,u_+) \leq \lamb_j(m; u_-)
\quad \mbox{ for all $m$ between $m_-$ and $m_+$,} 
$$
where $m_-= \mu(u_-)$ and $m_+= \mu(u_+)$.
\end{lemma}

We only give the proof of Lemma~\ref{Riem-5}, the proof of Lemma~\ref{Riem-4} being similar.

\medskip
\noindent
{\bf Proof.} 
Set $m_-:=\mu(u_-)$ and $m_+:=\mu(u_+)$ and, for definiteness,
assume that $m_-<m_+$. We must show that 
$$
\Lam := \lamb_j(u_-,u_+) \geq \lamb_j(m;u_+), 
\quad 
m_- \leq m \leq m_+.
$$
By contradiction let us assume that there exists a ``first'' point $m_0 \in [m_-,m_+)$ 
at which the above condition fails, that is 
\be
 \begin{array}{ll}
  \Lam \geq \lamb_j(m;u_+), \quad & m_- \leq m \leq m_0, \hfill \\[1ex]
  \Lam =    \lamb_j(m;u_+), \quad & m = m_0, \hfill \\[1ex]
  \Lam <    \lamb_j(m;u_+), \quad & 0 < m - m_0 \ll 1, \hfill
 \end{array}
 \label{Pro.10}
\ee
On the other hand, since the shock connecting $u_-$ to $u_+$ being entropy admissible, 
\be
\Lam \leq \lamb_j(m;u_-), \quad  m_- \leq m \leq m_+. 
\label{Pro.11}
\ee 
We treat the case where  $m_-<m_0<m_+$, the case where the point $m_0$ coincides with the left endpoint
$m_-$ of the interval being analogous.  
 
Let us expand $\lamb_j(m;u_+)$ in a neighborhood of $m=m_0$, 
$$
\lamb_j(m;u_+) =  \Lam + \frac{(m-m_0)^{p+1} }{ (p+1)!} \, \del_m^{(p+1)}\lamb_j(m_0;u_+) + \GO\bigl( (m-m_0)^{p+2}\bigr), 
$$
where $\del_m^{(p+1)}\lamb_j(m_0;u_+) \neq 0$. By Lemma~\ref{Riem-3} we know that up to a positive multiplicative constant
$$
\del_m^{(p+1)}\lamb_j(m_0;u_+)  = \hpi_j^{(p)}(u_0),
$$
while all $\hpi_j^{(k)}(u_0)=0$ for $k=1, \ldots, p-1$. 
On the other hand, since the shocks $(u_-,u_+)$ and $(u_0, u_+)$ propagate at the same speed, $\Lam$,
the shock speed of $(u_-,u_0)$ is also $\Lam$. (This is an elementary fact, already stated in Lemma~\ref{Riem-4}.)
Therefore $\lamb_j(m_0; u_-)=\Lam$ and, similarly as above, we can write 
$$
\lamb_j(m;u_-) =  \Lam + \frac{(m-m_0)^{q+1} }{(q+1)!} \, \del_m^{(q+1)} \lamb_j(m_0;u_-) 
+ \GO\bigl( (m-m_0)^{q+2}\bigr),
$$
where 
$$
\del_m^{(q+1)} \lamb_j(m_0;u_-)  = \hpi_j^{(q)}(u_0),
$$
while all $\hpi_j^{(k)}(u_0)=0$ for $k=1, \ldots, q-1$. Obviously, $q=p$. However, in view of conditions \eqref{Pro.10} 
on $\lamb_j(m; u_+)$ we see that $p$ must be even (with $\hpi_j^{(p)}(u_0)>0$), while 
the condition \eqref{Pro.11} on $\lamb_j(m; u_-)$ implies that $p$ must be odd (with $\hpi_j^{(p)}(u_0)>0$). This is a contradiction.  
\quad \qed 
\medskip 


\subsection{Construction of the wave curves}  
In the rest of this section we sketch the proof of Theorem~\ref{Riem-2}. Most of the arguments in 
\cite{IguchiLeFloch} carry over under our weaker assumption on $f$, and we will not repeat
them. We will only give the new ingredients of the proof.  
For definiteness, we assume that ${(\nabla\lam_j\cdot \rt_j)(u_0) > 0}$, and we construct the part $m > \mu(u_0)$ 
of the wave curve. Locally near $u_0$ we can use the integral curve $\OO_j(u_0)$.

Suppose that along the integral curve there exists a ``first point'' with coordinate $\mu^1(u_0)$ 
where $\nabla\lam_j\cdot r_j$ vanishes {\sl and changes sign,} that is by setting $u^1:=w_j(\mu^1(u_0); u_0)$ 
\be
 \begin{split}  
&  (\nabla\lam_j\cdot \rt_j)(w_j(m;u_0))>0, \quad \mu(u_0) \leq m < \mu^1(u_0), 
\\
&  (\nabla\lam_j\cdot \rt_j)(u^1) = 0, 
\\
&  (\nabla\lam_j\cdot \rt_j)(w_j(m;u_0)) < 0, \quad 0 < m - \mu^1(u_0) <<1, 
\\
& \hpi_j^{(k)}(u^1) =0 \, \, (k=1, \ldots, p), \quad \hpi_j^{(p+1)} (u^1) < 0. 
\end{split}
\label{Cons.1} 
\ee
Note that, clearly, $p$ must be an odd integer. 

It should be noted here that there can not be accumulation points in the critical set along an integral curve. 
Indeed, if $m^l \to m^\infty $ ($l \to \infty$) were a sequence such that 
$$
(\nabla\lam_j\cdot \rt_j)(w_j(m^l;u_0)) = 0. 
$$
Then, by induction on $k \geq 1$ and using the intermediate value theorem we could find sequences $m^{k,l}$ such that 
$$ 
\hpi_j^{(k)}(w_j(m^{k,l};u_0)) = 0. 
$$
Letting $l \to \infty$ we would conclude that, at the point $u^\infty:=w_j(m^\infty;u_0)$, 
$$
\hpi_j^{(k)}(u^\infty) = 0,  \qquad      k=1,2,\ldots, 
$$
which is impossible since $f$ is nondegenerate.

Note that the wave curve coincides with the integral curve
up to the value $\mu^1(u_0)$: 
$$
\psi_j(m;u_0) := w_j(m;u_0), \quad \mu(u_0) \leq m \leq \mu^1(u_0).
$$
We claim that $\mu^1$ is a continuous function of its argument. Namely, consider the function 
$$
F_j(m;u) := \hpi_j^{(1)}(w_j(m;u)), 
$$ 
which is smooth with respect to $m$ and $u$. Since by assumption the state $u_0$ is such that $F_j(\mu^1(u_0);u_0) = 0$,  
$$
(\del_m^{(p)} F_j)(\mu^1(u_0);u_0) = \hpi_j^{(p+1)}(w_j(m;u)) < 0.
$$ 
Therefore, one can locally solve the equation $F_j(\mu,u)=0$ in a neighborhood of $u_0$, in a continuous manner.


To extend the wave curve we need a pattern made of a rarefaction wave followed by a shock wave (actually a left-contact
wave).  For all meaningful values $(n,m;u)$ we set 
$$
\aligned 
& G_j(m,n;u) 
& := 
\begin{cases}
\displaystyle \frac{1 }{ m - n} \bigl(\lamb_j(m;w_j(n;u))
  -\lam_j(w_j(n;u))\bigr) ,   & m \neq n, 
\\  
\displaystyle \frac{1}{2} \, \Bigl( \nabla\lam_j \cdot \rt_j \Bigr) (w_j(m;u)) , 
& m=n, 
\end{cases}
\endaligned 
$$
where $u^1:=w_j(\mu^1(u); u)$. We will now apply the implicit function theorem to the equation 
\be
G_j(m,n;u) = 0, 
\label{Cons.01}
\ee 
and show that it defines a function $n=\nu^1(m;u)$ in the neighborhood
of the point $(m,u) = (\mu^1(u_0),u_0)$. As presented above, our construction is associated 
with the point $u^1$ at which $r_j \cdot \nabla \lam_j(u^1) = 0$.  
The special case of a PGNL flux is straightforward since it can be checked that 
$r_j \cdot \nabla \cdot (r_j \cdot \nabla \lam_j)(u^1) \neq 0$ 
implies that $\del_n G_j(\mu^1(u_0),\mu^1(u_0),u_0) \neq 0$. 
Handling a general nondegenerate flux is much more delicate and this is discussed now. 
We regard $u$ as a parameter and, for convenience, in the presentation 
we simply fix the point $u_0$ and skip the dependence in $u$. We are interested in 
the corresponding singular point $u^1$ along the rarefaction curve from $u_0$,
where $\hpi^{(p+1)}(u^1) \neq 0$ and we must analyze \eqref{Cons.01} near 
$$
m^1 := \mu^1(u_0). 
$$

Consider first the case of a scalar conservation law. The function $G_j$ is independent of the variable $u$, 
and the equation under consideration reads 
$$
\frac{f(n) - f(m)}{ n - m} - f'(n) =0. 
$$
Using that $f$ has a nondegenerate point at $m^1$, 
precisely, $f^{(p+2)}(m^1) \neq 0$ so that 
$$
f(m) \sim \frac{f^{(p+2)}(m^1) }{ (p+2)!} \, (m - m^1)^{p+2}, 
$$  
one can check 
geometrically on the graph of the function $f$ that the above condition determines a unique solution $n= \nu^1(m)$ 
which, furthermore, can be easily expanded near $m^1$: 
$$
\nu^1(m) \sim c_p \, ( m -m^1), 
$$
where the constant $c_p<0$ is the unique negative root of the equation
$$
\frac{1 - c_p^{p+2} }{ 1- c_p} = (p+2) \, c_p^{p+1}. 
$$ 

We now extend the above result to systems.

\begin{lemma} 
\label{regul}
For $\alpha, \beta>0$ sufficiently small, introduce the small cone 
$$
\Ccal := \big\{ |(n - m^1) - c_p \, (m - m^1)| 
\leq \alpha \, |m - m^1| \big\}. 
$$
Near the point $(m,n) = (m^1,m^1)$ within the cone $\Ccal$
the function $G_j$ admits the expansion 
\be
G_j(m,n) 
= Q_j(m,n) + \GO( 1) \, \big( (m-m^1)^{p+1} +  (n - m^1)^{p+1} \big), 
\label{Cons.expansion}
\ee
where 
\be
\aligned 
& Q_j(m,n) = \big( (n - m^1) - c_p \, (m - m^1) ) \, R_j(n,m), 
\\
& R_j(m,n) = \hpi^{(p+1)}(u^1) \, (m- m^1)^{p-1} \,  S_j(n,m), 
\endaligned  
\label{Cons.expansion2}
\ee 
and the function $S_j$ is smooth and close to $1$, 
\be
|S_j(m,n) -1| <\beta. 
\label{Cons.expansion3}
\ee 
\end{lemma} 

Once this lemma is established, the implicit function theorem straightforwardly applies to the 
mapping $\widetilde G_j$ defined by 
$$
\widetilde G_j(m,n) := \frac{G_j(m,n) }{ (m- m^1)^{p-1} \,  S_j(m,n)} \qquad \text{ in the cone } \Ccal, 
$$
and extended arbitrarily outside $\Ccal$ as a smooth map. By the above claim, 
we have 
$$
\widetilde G_j(m,n) = \hpi^{(p+1)}(u^1) \, \big( (n - m^1) - c_p \, (m - m^1) )  + \GO( 1)  \big( |m - m^1|^2 +  |n - m^1|^2 \big), 
$$
so that there exists a unique function $\nu^1=\nu^1(m)$ defined 
in a neighborhood of $m=m^1$ and such that 
$$
\widetilde G_j(m,\nu^1(m)) = 0,
$$
which, moreover, satisfies 
$$
\nu^1(m) = m^1 + c_p \, (m - m^1) + {o}(m - m^1).
$$
In turn, the solution remains within the cone $\Ccal$ (in a sufficiently small neighborhood of $m^1$, at least), 
and therefore $\nu^1(m)$ is a solution of the original equation \eqref{Cons.01}. 

\

It remains to prove Lemma~\ref{regul} above. 
Relying on the tangency property between the Hugoniot and integral curves (Lemma \ref{Riem-Tan}), one can replace
the Hugoniot curve by the rarefaction curve in the definition of $G_j$ while making an error of order 
$(m- m^1)^{p+1} +  (n - m^1)^{p+1}$. Let us introduce the reduced flux  $$
f_j(u) := \lh \cdot f(u) 
$$
and set
$$
\fhat_j (m) := f_j(w_j(m; u_0)). 
$$
We can rewrite the equation $G_j(m,n) = 0$ in the following form 
$$
\aligned 
& \frac{\fhat_j(m) - \fhat_j(n) }{ m-n} 
= \fhat'_j(n) +  
\GO( 1) \, \big( (m- m^1)^{p+1} +  (n - m^1)^{p+1} \big). 
\endaligned 
$$
The behavior of the flux in the neighborhood of $u^1$ leads to 
$$
\fhat_j(m) \sim \frac{\hpi^{(p+1)}(u^1) }{ (p+2)!} \, (m - m^1)^{p+2}, 
$$  
and the expression in \eqref{Cons.1} can be expanded exactly as in the scalar case, 
leading to the result stated in the claim.

The mixed curve is now defined locally, and we discuss its extension. Consider the shock speed of an arbitrary discontinuity connecting 
$w_j(n;u_0)$ to $v_j(m; w_j(n;u_0))$, that is, 
$$
\Lam_j(m,n) := \lamb_j(m;w_j(n;u_0)). 
$$
By construction, at ${n=\nu^1(m)}$ the shock speed coincides with the characteristic speed: 
$$
\Lam_j(m,\nu^1(m)) = \lam_j(w_j(\nu^1(m);u_0)).
$$ 
 
One can show that 
$$
\del_m \nu^1(m) = \frac{1}{ (\nabla\lam_j \cdot \rt_j)
 \bigl( w_j(\nu^1(m);u_0) \bigr)} \, 
 \bigl( \del_m\Lam_j \bigr) (m,\nu^1(m)).
$$
Since ${\nu^1}'(\mu^1) <0$, it follows that ${\nu^1}'(m) < 0$ for 
$m$ sufficiently close to $\mu^1$ and therefore 
$$
\bigl( \del_m \Lam_j \bigr)(m,\nu^1(m)) < 0, 
\quad 0 < m-\mu^1(u_0)  <<1. 
$$
This means that the speed of the left-contact decreases as $m$ increases. 
The function $\nu^1(m; u_0)$ is well-defined until $\del_m \Lam_j$ eventually vanishes (see in particular $\partial_m G_i$)
at some point with coordinate denoted by $\mu^2= \mu^2(u_0)$,
$$
(\del_m \Lam_j)(\mu^2,\nu^1(\mu^2; u_0)) =0. 
$$
The wave curve is determined from the map $m \mapsto v_j(m;w_j(\nu^1(m;u_0);u_0))$, 
by
\be
\psi_j(m;u_0) = \begin{cases} 
  w_j(m;u_0),         & \mu(u_0) \leq m  < \mu^1(u_0), 
  \\
  v_j(m;w_j(\nu^1(m;u_0);u_0)),  & \mu^1(u_0) \leq m \leq \mu^2(u_0). 
 \end{cases}
\nonumber  
\ee 
As in \cite{IguchiLeFloch} one can check that for each $m \in (\mu^1,\mu^2)$ the state $w_j(\nu^1(m);u_0)$ is connected to
the right-hand state $v_j(m;w_j(\nu^1(m);u_0))$ by a shock wave satisfying 
the entropy criterion.  All of the remaining arguments of the construction of the wave 
curve in \cite{IguchiLeFloch} can be extended in a similar way as we have explained above.

We observe that the construction requires finitely many waves only. 
Otherwise, if the wave curve was made of infinitely many shock and rarefaction pieces, then 
there would be an accumulation point of critical points, 
$m^l \to m^\infty $ ($l \to \infty$) such that 
$$
\hpi_j^{(1)}(\psi_j(m^l;u_0)) = 0. 
$$
Then, by induction on $k \geq 1$, 
using the intermediate value theorem and the higher-order tangency property between shock and rarefaction curves
as stated in Lemma ~\ref{Riem-3}, we would obtain 
$$ 
\hpi_j^{(k)}(\psi_j(m^\infty;u_0)) = 0
$$ 
for all $k$, which is impossible since $f$ is nondegenerate.

After the wave curves are constructed, the Riemann problem can be solved as in \cite{IguchiLeFloch}. 
This complete the discussion of Theorem~\ref{Riem-2}. Propositions~\ref{regu} and \ref{interactions}
are then immediate from the results in \cite{IguchiLeFloch}. 
From our construction, the following property of the Riemann solution follows immediately.

\begin{lemma} {\rm (Splitting property for the Riemann problem.)}
\label{LemCut}
Consider a wave fan $u_+=\psi_j(m_+;u_-)$, with $m_- < m_+$ for definiteness, 
associated with the intermediate states $u_0$, \dots, $u_{2q+1}$ as defined in Theorem~\ref{Riem-2}. Then the following 
two properties hold. 

1. For $k =0, \ldots, q$ and $m \in [\mu(u_{2k}),\mu(u_{2k+1})]$, 
the solution of the Riemann problem $(u_-, \psi_j(m;u_-))$ (resp.~$( \psi_j(m;u_-),u_+)$) 
coincides with the restriction of the solution of the Riemann problem $(u_-,u_+)$ to the interval 
$[\lamb_j(m_-,m_+;u_-), \lam_j(\psi_j(m;u_-))]$ 
(resp.  $[\lam_j(\psi_j(m;u_-)) , \lamb_j(m_+,m_+;u_-)]$).

2. For $k = 0, \ldots, q$ and $m \in [\mu(u_{2k+1}),\mu(u_{2k+2})]$, 
the solutions of the Riemann problems $(u_-, \psi_j(m;u_-))$ (respectively $(\psi_j(m;u_-),u_+)$) 
and $(u_-,u_+)$ coincide in the interval $[\lamb_j(m_-,m_+;u_-), \lam_j(u_{2k+1})]$
(respectively $[\lam_j(u_{2k+2}),\lamb_j(m_+,m_+;u_-)]$).
\end{lemma}

We end this section with another  property of Riemann solutions: roughly speaking,
the speeds within a wave fan can be determined by a convex or concave hull argument, 
which is similar to what is classically done for a scalar conservation law. 
This property will not be directly used in the rest of this paper.

\begin{lemma} {\rm (Construction of a wave packet using the convex hull.)}
\label{LemCConv} 
Let $u_+=\psi_j(m_+;u_-)$, with $m_-:=\mu(u_-) < m_+$, and 
introduce the intermediate states $u_1$,\ldots,$u_{2q+1}$ as in Theorem~\ref{Riem-2}. 
By setting 
$$
\aligned 
& f_j(m) := \lh \cdot f(\psi_j(m;u_-)), \qquad m_- \leq m \leq m_+, 
\endaligned 
$$
the convex hull $g_j := \conv_{[m_-,m_+]} f_j$ is precisely 
\be
\label{Convfi}
g_j (m) 
= 
\begin{cases} 
{\aligned 
f_j(m) = f_j(m_{2k}) & + \int_{m_{2k}}^m \lam_j(\psi_j(m';u_-)) \, dm', 
\\
& \hskip1.cm m_{2k} \leq m \leq m_{2k+1},  
\endaligned}
\\
{\aligned 
f_j(m_{2k+1}) + \lamb_j(u_{2k+1}, &u_{2k+2}) \, (m - m_{2k+1}), 
\\
& \hskip.6cm m_{2k+1} \leq m \leq m_{2k+2}.
\endaligned} 
\end{cases} 
\ee
Moreover, $f_j$ cannot be affine in a subinterval of $[m_-,m_+]$ where it coincides with $\conv_{[m_-,m_+]}f_j$
and, in consequence, the set of the intermediate values $m_k$ is nothing but 
the boundary of the set 
$$
\big\{ m \in [m_-,m_+]  / \  g_j (m) = f_j(m) \big\}. 
$$ 
When $m_+ < m_-$, the same statements hold by replacing the convex hull by the concave one.
\end{lemma}

Note that, even though the wave curves $\psi_j$ are only Lipschitz continuous, 
the function $g_j$ above is of class $C^1$, as this is clear from its definition.

\medskip
\noindent
{\bf Proof.} Let us here denote by $g_j$ the right-hand side of \eqref{Convfi}
(defined successively  on each interval $[m_k,mu_{k+1}]$. We will show that $g_j =\conv_{[m_-,m_+]} f_j$. 
Observe that $g_j$ is a convex function, as a consequence of the fact that
the wave speed function $\lamb_j(m,m_+;u_-)$ introduced earlier is non-increasing in $m$. 

First of all, within a rarefaction interval $[m_{2k}, m_{2k+1}]$ the identity 
\be 
\label{fiRar}
f_j(\psi_j(m,u_-)) = f_j(m_{2k}) + \int_{m_{2k}}^m \lam_j(\psi_j(m';u_-)) \, dm'
\ee
follows directly by differentiating this relation and using the definition of the rarefaction together with 
the normalization $\mu(v_j(m)) = m$. 

Second, considering a shock wave $(u_{2k}, u_{2k+1})$ and using the Rankine-Hugoniot relation,
we obtain (for the end point values) 
$$
\aligned 
f_j(u_{2k+2}) - f_j(u_{2k+1}) 
& =  \lamb_j(u_{2k+1},u_{2k+2}) \, \lh \cdot (u_{2k+2} - u_{2k+1}) 
\\
& = \lamb_j(u_{2k+1},u_{2k+2}) \, (m_{2k+2} - m_{2k+1}).
\endaligned 
$$
Hence, combining with the identity within rarefaction, we conclude that $f_j$ and $g_j$ 
coincide within rarefaction intervals, 
\be
\label{figi}
f_j(m) = g_j(m), \qquad  m \in [m_{2k},m_{2k+1}]. 
\ee

In consequence, to establish \eqref{Convfi} it is sufficient to prove that ${g_j \leq f_j}$. 
Moreover, since these functions coincide on the intervals $[m_{2k},m_{2k+1}]$, 
we need only show $g_j \leq f_j$ on each interval $[m_{2k+1},m_{2k+2}]$.  
So, we now establish that, for any $\tu= \psi_j(\tm;u_-)$ with 
$m_{2k+1} \leq \tm \leq m_{2k+2}$, 
\be
\label{fiaudessus}
\aligned 
f_j(\tm) 
& \geq f_j(u_{2k+1}) + \lamb_j(u_{2k+1},u_{2k+2}) (\tm - m_{2k+1}) 
\\
& = g_j(\tm).
\endaligned 
\ee

To this aim, we consider the Riemann problem $(u_-,\tilde{u})$ and the restriction of $f_j$ to the interval $[m_-,\tm]$, 
$$
\fh_j:=f_{j |[m_-,\tm]}
$$
and we define $\gh_j$ by \eqref{Convfi} with the intermediate states $u'_1,\dots,u'_{2q'+1}$ corresponding to this problem. 
It follows from previous observations that, at the end point $\tm$, 
$$
\fh_j(\tm) = \gh_j(\tm).
$$
In view of the splitting property in Lemma~\ref{LemCut}, it is clear that $\gh_j(m_{2k+1}) = f_j(m_{2k+1})$,
and 
\be
\label{figi2}
\gh_j' (m) \geq \gh_j'(m_{2k+1}), \qquad m_{2k+1} \leq m \leq \tm.
\ee
Now, if $u_{2k+1} \not = u_-$, we have 
$$
\aligned 
\gh_j'(m_{2k+1}) 
& = \fh_j'(u_{2k+1}) 
\\
& = \lam_j(u_{2k+1}) =\lamb_j(u_{2k+1},u_{2k+2}),
\endaligned 
$$
since the two Riemann solutions coincide on a whole interval. If $u_{2k+1}  = u_-$, we have 
$$
\gh'_j(m_{2k+1}) = \lamb_j(u_{2k+1},\tilde{u}) \geq \lamb_j(u_{2k+1},u_{2k+2}),
$$
as follows from the entropy criterion. This establishes \eqref{fiaudessus}.

Consider now the last statement in the lemma and 
suppose that $f_j$ is affine and coincides with $\conv f_j$ in some subinterval $[m_1,m_2]$ ($m_1<m_2$). 
Clearly, $(m_1,m_2)$ does not intersect a rarefaction interval 
$(m_{2k},m_{2k+1})$ since this would contradict the nondegeneracy assumption. 
Introducing $\tu$ as above, we see that 
the inequality \eqref{figi2} is strict unless $(u_{2k+1}, \tu)$ 
is a single shock of the same speed as $(u_{2k+1},u_{2k+2})$. 
But, if $v_j(\cdot;u_{2k+1})$ had this property on $[m_1,m_2]$, 
it would follow that $v_j(\cdot;u_{2k+1}) = w_j(\cdot;u_{2k+1})$ in that interval which, again, would 
contradict the nondegeneracy of the flux. \quad\qed


\subsection{Proof of the inner speed variation estimates}

\medskip
\noindent
{\bf Proof of Theorem~\ref{PropISV}.} \,  
1. We only treat the case $m \geq 0$ since the case $m<0$ is similar. From the estimates on the strengths we have: 
\be
\label{ReEstForces}
 \sup_{0 \leq n \leq m} | \psi_j(\mu(u_-) + n ; u_-) - \psi_j(\mu(u'_-) + n ; u'_-)  |
\lesssim |m| \, |u'_- - u_-|.
\ee
Observe also that the variation of the shock speeds is essentially equiavelent to the variation of 
the characteristics speeds, as follows 
\be
\label{LemmeDiffVitChoc}
\aligned 
& \lamb_j(\mu(u')+n; u') - \lamb_j(\mu(u)+n; u) 
\\
& = \lam_j(u') - \lamb_j(u) + \GO(n) \, |u' - u|. 
\endaligned 
\ee
This is a consequence of the fact that all functions under consideration are smooth, 
and that the equality holds (without remainder) when either $n=0$ or $u'=u$.

Denote $u'_+:= \psi_j(\mu(u'_-)+m;u'_-)$ and $u_+:= \psi_j(\mu(u_-)+m;u_-)$. To prove the first statement 
in Theorem~\ref{PropISV}, we will establish 
\be 
\label{ISVBis}
\aligned 
& \lamb_j^{\max}(\mu(u'_-)+m;u_-') - \lamb^{\max}_j(\mu(u_-)+m;u_-) 
\\
& = \lam_j(u'_+) - \lam_j(u_+) + \GO(m) \, |u'_+ - u_+|
\endaligned 
\ee 
\be
\label{ISVBis2}
\aligned 
& \lamb^{\min}_j(\mu(u'_-)+m;u_-') - \lamb^{\min}_j(\mu(u_-)+m;u_-) 
\\
& = \lam_j(u'_-) - \lam_j(u_-) + \GO(m) \, |u'_- - u_-|. 
\endaligned 
\ee 
This implies the desired estimate, since by \eqref{ReEstForces}, 
$$
\lam_j(u'_+) - \lam_j(u_+)  
= \lam_j(u'_-) - \lam_j(u_-) + \GO(m) \, |u'_- - u_-|
$$
and $|u'_+ - u_+| = |u'_- - u_-| +\GO(m) \, |u'_- - u_-|$.

We only consider \eqref{ISVBis}, since \eqref{ISVBis2} can be checked similarly. 
To prove \eqref{ISVBis}, it is sufficient to establish the {\sl inequality}
\be
\label{ISVBis3}
\aligned 
& \lamb^{\max}_j(\mu(u'_-)+m;u_-')  - \lamb^{\max}_j(\mu(u_-)+m;u_-)
\\
& \leq  \lam_j(u'_+) - \lam_j(u_+) + \GO(m) \, |u'_+ - u_+|. 
\endaligned 
\ee
This is so since $u_-$ and $u'_-$ play completely symmetric roles.

Each of the Riemann problems $(u_-,u_+)$ and $(u'_-,u'_+)$ are solved by a succession of 
rarefaction waves and shock waves. 
Call $\underline{m}$ (resp. $\underline{m}'$) the real in $[0,m]$ defined as follows: 
\begin{itemize}
\item If the problem $(u_-,u_+)$ (resp. $(u'_-,u'_+)$)  ends with a non-trivial rarefaction wave, 
we set $\underline{m}=m$ (resp. $\underline{m}'=m$),
\item If the problem $(u_-,u_+)$ (resp. $(u'_-,u'_+)$)  ends with a non-trivial shock wave 
$(\tilde{u},u^+)$ (resp. $(\tilde{u},u'^+)$),  we 
set $\underline{m}=\mu(\tilde{u})-\mu(u^-)$ (resp. $\underline{m}'=\mu(\tilde{u})-\mu(u'^-)$).
\end{itemize}

We distinguish between two main cases: 
\begin{itemize}
\item Case 1: if $\underline{m}' \geq \underline{m}$. By the entropy criterion, the shock speed enjoys a 
monotonicity property and we have
$$ 
\lamb_j(\mu(u_-) + \underline{m}' ; u_+) 
\leq \overline{\lambda}_j(\mu(u_-) + \underline{m} ; u_+) = \lam_j^{\max}(u_-,u_+), 
$$ 
and the conclusion follows easily from \eqref{LemmeDiffVitChoc}.

\item Case 2: if $\underline{m}' < \underline{m}$. In that case, there are two subcases:
\begin{itemize}
\item 2a. If the point $\hat{u}$ of parameter value $\mu(u_-)+\underline{m}' $ in the problem $(u_-,u_+)$ 
corresponds to a rarefaction wave. In that case, we clearly have 
$$
\overline{\lambda}_j\big(\mu(u_-)+m,\mu(u_-)+m;u_-\big) \geq \lambda_j(\hat{u}). 
$$
Hence, using Lemma \ref{LemCut}, we can prove \eqref{ISVBis3} as in Case~1, by considering 
the Riemann problems $(u_-,\psi_j(\mu(u_-)+\underline{m}';u_-))$ and $(u'_-,\psi_j(\mu(u'_-)+\underline{m}';u'_-))$. 
(In fact, in this case we only compare characteristic speeds.)

\item If the point $\hat{u}$ of parameter value $\mu(u_-)+\underline{m}' $ in the problem $(u_-,u_+)$ 
corresponds to a shock wave, say $(u_1,u_2)$, then 
$$
\overline{\lambda}_j(\mu(u_-)+m,\mu(u_-)+m;u_-) \geq \lambda_j({u_2}). 
$$
Again, using Lemma \ref{LemCut}, we can apply Case 1 to the problems $(u_-,u_2)$ and 
$(u'_-,\psi_j(\mu(u'_-)+\mu(u_2)-\mu(u_-);u'_-))$.
Note here that the latter problem coincides with $(u'_-,u'_+)$ up to the parameter 
$\mu(u'_-)+\underline{m}'$ and has speeds greater than the ones in $(u'_-,u'_+)$.
\end{itemize}
\end{itemize}

\noindent
2. Fix a left-hand state $u_l$, and let $u_m, u_r'$ be two states
on the $j$-wave curve from $u_l$, with $u_r' := \psi_j(\mu_r; u_l)$.
Using the standard Lipschitz continuity property of the 
wave speed along a given wave curve based at the point $u_l$ we obtain 
$$
\begin{array}{l}
\lamb_j^{\max} (u_l, u_r') = \lamb_j^{\max} (u_l, u_m) + \GO(\mu_m - \mu_r), \\
\lamb_j^{\min} (u_l, u_r') = \lamb_j^{\min} (u_l, u_m) + \GO(\mu_m - \mu_r), 
\end{array}
$$
which yields the desired statement on the inner speed variation
$\vartheta_j(u_l, u_r) = \lamb_j^{\max} (u_l, u_r') - \lam^{\min}_j (u_l,u_r')$. \\

\noindent
3. Introduce $\tu_r:= \psi_j(\mu_r; u_l)$. As we could construct ``left-hand'' wave curves as we did for (right-hand) 
wave curves, 
one can introduce $\tu_m$ such that
$$
\mu(\tu_m) = \mu_m, \qquad 
\tu_r = \psi_j(\mu_r;\tu_m).
$$
It follows from the regularity of the left-hand wave curves and the wave interaction estimates that  
$$
| \tu_m - u_m| =\GO(| \tu_r - u_r|) = \GO(1) \, Q_{}(u_l,u_m,u_r). 
$$
Now, using the first statement of the theorem, we deduce
$$
\lamb_j(\cdot,\mu_r; \tu_m) = \lamb_j(\cdot,\mu_r;u_m) + \GO(1) \, Q_{}(u_l,u_m,u_r).  
$$

We introduce two values of the parameter in the Riemann problem $(u_l,\tu_r)$:
\begin{eqnarray*}
&m_a :=  \max \big\{ m \in [\mu_l,\mu_m], \ /  \text{ the point in the problem } (u_l,\tu_r) \hfill \\
&\hspace{3cm}  \text{ with parameter } m \text{ corresponds to a rarefaction wave} \big \},  \\
&m_b := \min \big\{ m \in [\mu_m,\mu_r], \ / \  \text{ the point in the problem } (u_l,\tu_r)  \\
&\hspace{3cm} \text{ with parameter } m \text{ corresponds to a rarefaction wave} \big \},
\end{eqnarray*}
with the convention that $m_a=\mu_l$ (resp. $m_b=\mu_r$) when there is no rarefaction wave in the corresponding parameter range.
In the terminology of Lemma \ref{LemCConv}, the above conditions can be stated as $f_j(m)=\conv_{[\mu_l,\mu_r]} f_j(m)$. 

We now distinguish several cases that depend upon the values of $m_a$ and $m_b$.
\begin{itemize}
\item Suppose that $\mu_l < m_a \leq m_b < \mu_r$. In that case, using Lemma~\ref{LemCut},
we see that the Riemann problem $(u_l,\tu_r)$ can be obtained by gluing three Riemann problems together: 
$(u_l, \psi_j(m_a;u_l))$, $(\psi_j(m_a;u_l),\psi_j(m_b;u_l))$ (a single contact discontinuity), and $(\psi_j(m_b;u_l),\tu_r)$. 
We then deduce that
$$
\aligned 
& \vartheta_j(u_l,\tu_r) 
\\
& = \lamb_j(\mu_r,\mu_r;\tu_m) - \lamb_j(\mu_l,\mu_m;u_l) 
\\
& = \lamb_j(\mu_r,\mu_r;\tu_m) - \lamb_j(\mu_m,\mu_r;\tu_m) + \lamb_j(\mu_m,\mu_r;\tu_m) - \lamb_j(\mu_l,\mu_m;u_l) 
\\
& = \vartheta_j(u_m,\tu_r)  + \lamb_j(\mu_m,\mu_r;u_m) - \lamb_j(\mu_l,\mu_m;u_l) + \GO(1) \, Q_{}(u_l,u_m,u_r),
\endaligned
$$ 
and the desired inequality \eqref{56} follows.

\item Suppose $\mu_l = m_a < m_b < \mu_r$ (resp. $\mu_l < m_a < m_b = \mu_r$). 
In that case, by the same arguments, we have  $\vartheta_j(u_l,\tu_r) = \vartheta_j(\psi_j(m_b;u_l), \tu_r)$
(resp. $\vartheta_j(u_l,\tu_r) = \vartheta_j(u_l,\psi_j(m_a;u_l))$) and, again, \eqref{56} follows.

\item Suppose $\mu_l = m_a < m_b = \mu_r$. This case is obvious since the left-hand side of \eqref{56} vanishes. 
\end{itemize}
This completes the proof of Theorem~\ref{PropISV}. \quad\qed
\medskip


\subsection{Proof of the density property} 

We will now establish that the nondegeneracy condition is generic. This section will cover both 
{\sl conservative and nonconservative} systems, characterized by a flux $f$ or a matrix $A$, respectively. 

Let $\Omega$ be an open set in  $\RN$. We are interested in matrix-valued mappings $A=A(u)$ 
($u \in \Omega$)  
of class $C^\infty$ satisfying the strict hyperbolicity property, that is, 
for any $u \in \Omega$,  the matrix 
$A(u)$ admits $N$ distinct and real eigenvalues, 
\be
\label{Hyperb}
\lam_1(A,u) < \ldots <  \lam_N(A,u),
\ee 
and basis of left- and right-eigenvectors $l_j(A,u), r_j(A,u)$, normalized so that
\be
\label{NormalisationLR}
\aligned 
& |r_j(A,u)| =1, \quad   l_j(A,u) \cdot r_j(A,u) = 1, 
\\
&  l_i(A,u) \cdot r_j(A,u)= 0      \quad (i \neq j). 
\endaligned  
\ee
Let $\Hcal$ be the subset of $C^\infty(\Omega)$ consisting of all mappings $A$ 
satisfying the strict hyperbolicity condition \eqref{Hyperb} and such that the associated maps 
$\lam_j, l_j, r_j$ are smooth. 
Recall that, in a sufficiently small neighborhood of any constant matrix and 
under the assumption \eqref{Hyperb}, the maps $\lam_j, l_j, r_j$ indeed depend smoothly upon $u$. 
In this subsection however, we need not assume that the set $\Omega$ is small.

To every map $A \in \Hcal$ we associate the following functions, by induction,
\be
\begin{split} 
& \pi_j^{(1)}(A,u) :=  r_j(A,u) \cdot \nabla_u \lam_j(A,u), 
\\
& \pi_j^{(k+1)}(A,u) : = r_j(A,u) \cdot \nabla_u \pi_j^{(k)}(A,u), \quad k = 1,2, \ldots 
\end{split}
\label{DefAlpha}
\ee
We claim that: 

\begin{theorem}
\label{ApproxFlux} {\rm (Density of the set of nondegenerate systems.)} 
Given any matrix-valued map $A \in \Hcal$ there is a sequence $A^l \in \Hcal$ such that
\be
\label{AFCV}
A^l \rightarrow A  \, \text{ in the (strong) $C^\infty$ Whitney topology,}
\ee
and $A^l$ is nondegenerate in the sense that, 
for $l=1,2,\ldots$, $u\in \Omega$, and $j=1, \ldots, N$, 
\be
\big(\pi_j^{(1)}(A^l,u), \pi_j^{(2)}(A^l,u), \ldots, \pi_j^{(N+1)}(A^l,u) \big) 
\neq (0,0, \ldots, 0).  
\label{generic}
\ee
Moreover, in the special case that $A$ is conservative, i.e.~$A=Df$, the sequence can be chosen 
to be conservative, $A^l = Df^l$, and $f^l$ converges to $f$ in the Whitney topology.  
\end{theorem}

Hence, all the results established earlier in this section for nondegenerate flux  
extend to arbitrary, strictly hyperbolic flux, and this leads us to Corollary~\ref{lecorollaire}.  
In the following, a map $A \in \Hcal$ (conservative or not) is fixed, and we consider the set 
$$
\Fcal := \big\{ (u,A(u)), \ u \in \Omega \big\}. 
$$
Clearly, \eqref{Hyperb} is an open condition and, therefore, 
for each $u \in \Omega$ one can find $\delta>0$ such that the ball $\Bcal_\delta(u)$ centered at $u$
is included in $\Omega$ and all matrices $B \in \Bcal_\delta(A(u))$ have $N$ real and distinct eigenvalues:
$$ 
\Lambda_1(B) < \ldots <  \Lambda_N(B),
$$ 
and, therefore, basis of left- and right-eigenvectors
$L_1(B), \ldots, L_N(B)$ and $R_1(B), \ldots, R_N(B),$
which we will normalize so that 
\be
\label{NormalisationLRb}
\aligned 
& |R_j(B) | = 1, \qquad  L_i(B) \cdot R_i(B) = 1, 
\\
&  L_i(B) \cdot R_j(B) = 0 \quad (i \neq j).  
\endaligned 
\ee
Clearly, the maps $\Lambda_j, L_j, R_j$ depend smoothly upon $B$ (at least locally), and represent extension of 
the maps $\lambda_j, l_j, r_j$ to the ball $\Bcal_\delta(A(u))$.  
Reducing $\delta$ if necessary, we can assume that for all $B_1, B_2 \in \Bcal_\delta(A(u))$
$$ 
L_i(B_1) \cdot R_i(B_2) \geq 1/2.
$$

Now, in view of the paracompactness property of $\Fcal$, there exists a locally finite covering of $\Fcal$ 
by open sets finer than $\Bcal_{\delta/2}((u,A(u)))$. 
We define  $\Gcal$ as the unions of all balls in such a covering. 
We define $R_j(B)$ and $L_j(B)$ on $\Gcal$ for $B \in \Bcal_{\delta/2}(A)$ as the value of $R_j(B)$ determined in that ball. 
It follows easily from this construction that this does not depend on the ball that contains $B$. 
Hence, we conclude that in $\Gcal$ the maps $\Gcal \ni A \mapsto \Lambda_j(A), L_j(A), R_j(A)$ are smooth.

To fix the idea we discuss the conservative case, with $A=Df$ for some flux $f$. Minor modifications are needed
to cover the nonconservative case, which we omit.

Let $\Ocal \subset C^\infty(\Omega, \RN)$ be a small neighborhood of the flux $f$, chosen so that
$(u,Dg(u)) \in \Gcal$ for all $g \in \Ocal$ and $u \in \Omega$. Introduce the subset 
$$ 
\aligned 
\Jcal_k(\Ocal) 
& = \big\{ j^k_{g}(u):=\big( u, g(u), g^{(1)}(u), \ldots, g^{(k)}(u)\bigr)
 \, / \,   u \in \Omega, \ g \in \Ocal \big\} 
 \\
 & \subset J^k(\Omega,\RN), 
\endaligned 
$$
consisting of all $k$-th jets of maps from
$\Omega$ into $\RN$. Here, $g^{(j)}$ denotes the $j$-order differential of the map $g$, 
and the mapping $j^{k}_g$ is the {\bf $k$-th jet extension of $g$.}

It follows easily from the definition \eqref{DefAlpha} that the functions $u \mapsto \pi_j^{(k)}(f,u)$ can be expressed 
in terms of jets, that is 
\be
\label{DefPhii}
\aligned 
\pi_j^{(k)}(f,u) & = \varphi_{k}(u, f(u), f^{(1)}(u), \ldots, f^{(k+1)}(u))
\\
&  =  \varphi_{k}(j^{k+1}_f(u)),
\endaligned 
\ee
where the functions $\varphi_k : J^{k+1}(\Omega,\RN) \rightarrow \RR$ are smooth in $\Jcal_{k+1}(\Ocal)$. 
Note in passing that $\pi_j^{(k)}(f,u)$ does not depend upon the first two components of the jet, $u$ and $f(u)$.

The proof of Theorem~\ref{ApproxFlux} follows from Thom's transversality theorem, 
which we recall for readers' convenience (for a proof see, for instance, \cite{Hirsch}), 
and a technical lemma that we establish below. 

\begin{theorem}[Thom]
\label{Thom}
Let $X$ and $Y$ be two smooth manifolds, and $Z$ be a submanifold of the $k$-th jet space $J^k(X,Y)$ for some 
$k \geq 0$. Then, there exists a Baire set of second category $E \subset C^\infty(X,Y)$ 
for the (strong) $C^\infty$ Whitney topology
such that 
the $k$-th jet extension of any $f \in E$ is transverse to $Z$, and consequently
$(j^k_f)^{-1} (Z)$ is either the empty set or a submanifold of $M$ having the same codimension as $Z$.
\end{theorem}

Recall that a Baire set is a set whose complement is the union of at most countably many, nowhere dense sets. 

\begin{lemma}
\label{MatSsVar}
For any $k \geq 2$, the set
$$
\aligned 
Z & := \Big\{ q \in \Jcal_k(\Ocal) \, / \, \varphi_j(q)=0, \, j =1, \ldots, k-1  \Big\}
\\
& \subset J^k(\Omega,\RN)
\endaligned 
$$
is a submanifold of codimension $k-1$. 
\end{lemma}

\medskip
\noindent
{\bf Proof of Theorem~\ref{ApproxFlux}.} 
Applying Lemma~\ref{MatSsVar} with $k=N+2$, we obtain a submanifold $Z$ of $J^{N+2}(\Omega,\RN)$. 
By Theorem~\ref{Thom}, 
the set $E$ of maps $g:\Omega \rightarrow \RN$ such that $g$ is transversal to $Z$ is a Baire set of second category. 
For $g \in E$ we introduce
$$
\Psi_g: u \mapsto ( \pi_j^{(1)}(g,u), \ldots, \pi_j^{(k-1)}(g,u)).
$$
For any $g \in E$, $\Psi_g^{-1}(0)$ is either empty or is a submanifold of codimension $N+1$. 
Clearly, there are no submanifold of codimension $N+1$ in $\Omega \subset \RN$, hence $\Psi_{g}^{-1}(0) = \emptyset$. 
Thus, one can find a sequence $f^k \in E$ converging to $f$ for the Whitney topology.  
\quad \qed 
\medskip

\medskip
\noindent
{\bf Proof of Lemma~\ref{MatSsVar}.} 
It is enough to prove that the differential $d\psi_k$ of the mapping 
$$
\psi_k: J^k(\Omega,\RN) \ni q \mapsto ( \varphi_1(q), \ldots, \varphi_{k-1}(q))
$$
has constant maximal rank $k-1$ for each $q$ in $J^k(\Omega,\RN)$. 
Let us write $(A_0,A_1, \ldots, A_{k-1})$ for a general vector of the tangent space at $q$, 
$T_q J^k(\Omega,\RN)$, 
with $A_0, A_1 \in \RN$, $A_i \in \RR^{d_i}$, $i \geq 1$, where $d_i =(N+i-1)! /[(N-1)! i!]$. 

By a straightforward induction argument we can check that it is sufficient to prove, for every $k$, 
\be
\label{Range}
(0, \ldots, 0, 1) \in \mbox{Range}(D\psi_k).
\ee
Denoting by $f^{(k)}$ the $k$-th order differential of $f$, we first observe that
$$
\aligned 
\pi_j^1(f,u) 
& = D_u[\Lambda_i(f^{(1)}(u))] (R_i(f^{(1)}(u))) 
\\
& = (D_A \Lambda_i)(f^{(1)}(u)) \cdot (f^{(2)})(u) \cdot \big[ R_i( f^{(1)} (u) ) \big] ,
\endaligned
$$ 
where $f^{(2)}$ is regarded as an element of ${\mathcal L}(\RN; {\mathcal L}(\RN,\RN))$. 

It is straightforward to check, by induction, that the term in $\varphi_{k-1}$ containing the highest-order differential of $f$ 
is 
$$ 
(D_A \Lambda_i)({f^{(1)}(u)}) \cdot f^{(k)}(u) \cdot [ \otimes^{k-1} R_i(f^{(1)}(u))],
$$
where $f^{(k)}$ is regarded as an element of  ${{\mathcal L}((\RN)^{\otimes k-1}; {\mathcal L}(\RN,\RN))}$. In other words, 
we have 
$$
\aligned
& \varphi_{k-1}(u, f(u), \ldots, f^{(k)})(u)) 
\\
& = (D_A \Lambda_i)({f^{(1)}(u)}). f^{(k)}(u). [\otimes^{k-1} R_i(f^{(1)}(u))] + \psi_{k-1}(u, f(u), \ldots, f^{(k-1)}),
\endaligned 
$$
where $\psi_{k-1}$ is a smooth function of $j^{k-1}_f(u)$.

Hence, to establish \eqref{Range}, it is enough to prove that for all $A \in \Hcal$ there exists 
$$
A_{k} \in \displaystyle{{\mathcal L}( (\RN)^{\otimes k-1}; {\mathcal L}(\RN,\RN))}
\approx \displaystyle{{\mathcal L}( (\RN)^{\otimes k};\RN)}
$$
such that 
\be
\label{But}
D_A \Lambda_i(A) \cdot A_{k} \cdot (\otimes^{k-1} R_i(A)) \not =0.
\ee
Note that $A$, $R_i(A)$ and $\Lambda_i(A)$ are related together via the identity 
$$
A \, R_i(A) = \Lambda_i(A) R_i(A).
$$
Differentiating this identity with respect to $A$ (in the direction $H$), we obtain 
$$
\aligned 
& H \, R_i(A) + A \cdot D_A R_i(A) \cdot H 
\\
& = (D_A \Lambda_i(A) \cdot H) \, R_i(A) + \Lambda_i(A) \, D_A R_i(A) \cdot H.
\endaligned 
$$
Taking the product by the left eigenvector $L_i(A)$ and using \eqref{NormalisationLR}, this yields
$$
D_A \Lambda_i(A) H = L_i(A) H R_i(A). 
$$ 
Hence, we arrive at \eqref{But} by considering $A_{k}$ sending $\displaystyle{\otimes^{k} R_i(A)}$ to $R_i(A)$ 
and other elements of a basis of $(\RN)^{\otimes k}$ to $0$, and $A_0=0$,\ldots, $A_{k-1}=0$. 
\quad \qed
\medskip 

\begin{rema}
It follows from the above proof, that generically in $f$ (or in $A$ in the non-conservative case), one has 
the stronger property that 
$$
\Big\{ u \in \Omega \ \Big/ \ \big(\pi_j^{(1)}(u), \ldots, \pi_j^{(k)}(u) \big) 
\Big\} 
$$ 
is a (possibly empty) submanifold of codimension $k$. \par
As a consequence, for such a flux and for almost every $u \in \Omega$, the (Lipschitz continuous) 
wave curve starting from $u$ does not meet any point where $(\pi_j^{(1)}(u), \pi_j^{(2)}(u)) =(0,0)$.
\end{rema}


\subsection{The PGNL condition is not generic} 

We have established that any flux can be approached by a sequence of flux satisfying
$$
(\pi_j^{(1)}(u), \ldots,\pi_j^{(N+1)}(u)) \not = (0, \ldots,0)
$$
for all $u$ and $j$. We will now prove that, in general, a flux cannot be approached by a sequence 
of flux satisfying the stronger condition 
$$
 (\pi_j^{(1)}(u), \ldots,\pi_j^{(N)}(u)) \not = (0, \ldots,0) 
$$
for all $u$ and $j$. Hence, the latter condition is not generic. 
In particular, for $N \geq 2$, there exist fluxes that cannot be smoothly approached by PGNL fluxes. 
Note that our counterexample below is ``local'', in the sense that it persists even if we shrink the neighborhood of the base point 
under consideration. Moreover, our result is valid for both conservative or nonconservative hyperbolic systems.

Consider flux $f :\RR^N \rightarrow \RR^N$ of the following triangular form:
$$ 
f(u_1,\ldots,u_N)= \left( \begin{array}{c}
f_1(u_1) \\
f_2(u_1,u_2) \\
\vdots \\
f_N(u_1,\ldots,u_N) 
\end{array} \right), 
$$
whose Jacobian matrix $A(u_1,\ldots,u_N)$ is given by 
$$ 
\left( \begin{array}{ccccc}
\partial_{u_1} f_1(u_1) & 0 & \ldots & \ldots &  0 \\
\partial_{u_1} f_2(u_1,u_2) & \partial_{u_2} f_2(u_1,u_2)  & 0 & \ldots &  0 \\
\vdots & \vdots & \ddots & & \vdots \\
\partial_{u_1} f_N(u_1,\ldots, u_N) & \partial_{u_2} f_N(u_1,\ldots, u_N)  & \ldots &   & \partial_{u_N} f_N(u_1,\ldots, u_N) \\
\end{array} \right). 
$$
Under the  assumption 
$$
\partial_{u_1} f_1 < \partial_{u_2} f_2 <  \ldots < \partial_{u_N} f_N, 
$$
the system of conservation laws associated with $f$ is strictly hyperbolic, with 
$$ 
\lambda_j(u_1, \ldots,u_N)=\partial_{u_j} f_j(u_1,\ldots,u_N), 
\qquad 
r_N(u_1,\ldots,u_N)= \left( \begin{array}{c}
0 \\
\vdots \\
0 \\
1
\end{array} \right). 
$$

It will be convenient to choose the component $f_N$ to be  the following polynomial expression 
$$ 
f_N(u_1,\ldots,u_N) :=  u_1 \, u_N^2 + u_2 \, u_N^3 + \ldots + u_{N-1} \, u_N^N + u_N^{N+2},
$$
while $f_1, \ldots, f_{N-1}$ are chosen arbitrarily but so that the strict hyperbolicity property holds
for all $u$ in the ball $\Bcal_\delta$. One can compute
\be
\label{Valeurs}
\aligned 
& \lambda_N(u_1, \ldots, u_N) 
\\
& = 2 u_1 u_N + 3 u_2 u_N^2 + \ldots + N u_{N-1}u_N^{N-1} + (N+2) u_N^{N+1},
\endaligned 
\ee
thus, since $r_N \cdot \nabla = \del_{u_N}$, 
\be
\label{Valeurs2}
\aligned 
& \pi^{(1)}_N(u_1, \ldots, u_N)  := r_N \cdot \nabla \lambda_N (u_1, \ldots, u_N)
\\
& =  2 \, u_1 + \ldots + N (N-1) \, u_{N-1} \, u_N^{N-2} + (N+2)(N+1) \,u_N^{N},
\endaligned 
\ee
and so on for all $\pi_j^{(N-1)}$, in particular  
\be
\label{Valeurs3}
\begin{array}{l}
\pi^{(N-1)}_N(u_1, \ldots, u_N)  = N! \ u_{N-1} + \frac{(N+2)!}{2} \, u_N^2, 
\\
\pi^{(N)}_N(u_1, \ldots, u_N)  = (N+2)! \ u_N.
\end{array}
\ee

Now, we claim that: 

\begin{proposition}
\label{PropContrex} Consider the triangular flux $f$ introduced above. 
For any $\delta \in (0,1)$ and for all $g$ sufficiently close to $f$ on $\Bcal_\delta$ (in the $C^{N+1}$-norm), 
there exists a state $u^* \in \Bcal_\delta$ such that 
\be 
\label{DegenN} 
(\pi_N^{(1)}(u^*), \ldots, \pi_N^{(N)}(u^*)) = (0, \ldots,0).
\ee
\end{proposition}

Our proof below remains valid if, instead of flux $g$, we consider mappings 
$A=A(u)$ ($u \in \Bcal_\delta$) which are sufficiently close to $df$.

\

\noindent{\bf Proof .} 
Clearly, any perturbation $g :\Bcal_\delta \to \RR^N$ of $f$ is  still strictly hyperbolic. 
Denote by $R_N[g],\Lambda_N[g]$ the eigenvector and eigenvalues of $Dg$ 
that are associated with the $N$-characteristic family, and define 
$\pi^{(1)}_{N}[g](u), \ldots, \pi^{(N)}_N [g](u)$ in the usual way. 
To exhibit a state satisfying \eqref{DegenN}
we consider the following vector field $W[g]: \overline\Bcal_\delta \to \RR^N$:
$$
W[g](u_1, \ldots,u_N) := \left( \begin{array}{c}
\pi^{(1)}_N [g] (u_1, \ldots,u_N) \\
\vdots \\
\pi^{(N)}_N [g] (u_1, \ldots,u_N) 
\end{array} \right).
$$ 

From \eqref{Valeurs}--\eqref{Valeurs3} one sees that, for any $(u_1, \ldots, u_N) \in \Bcal_\delta$, 
$$
W[f](u_1, \ldots,u_N) \cdot (u_1, \ldots,u_N) 
= \sum_{j=1}^{N-1} (j+1)! \, u_j^2 + (N+2)! \, u_N^2 + P(u_1, \ldots, u_N),
$$
where $P(u_1, \ldots, u_N)$ is a polynomial expression in $u_1, \ldots, u_N$ such that 
which each term is of degree $3$, at least. It follows that for some $\delta_0>0$
and for all $\delta \in (0,\delta_0)$ and $(u_1, \ldots, u_N) \in \Scal_\delta$ (the sphere with radius $\delta$) 
$$
W[f](u_1, \ldots,u_N) \cdot (u_1, \ldots,u_N) > 0. 
$$
In consequence, for $\delta \in (0,\delta_0)$, 
the index of the vector field $W[f]$ around the sphere $\Scal_\delta$ is $1$ 
(by a standard homotopy argument, for instance). Now, it follows that for all $g$ close to $f$, 
the vector field $W[g]$ does not vanish on $\Scal_\delta$ and is of index $1$ around $\Scal_\delta$. 
Hence, $W[g]$ necessarily vanishes inside  $\Bcal_\delta$, and this establishes the claim. \quad 
$\Box$


\section{A second proof of the inner speed variation estimates}
\label{Sec:BB}
 
\subsection{Notation and preliminaries}                                                                                                                                                                                                                                                                                                                                                                                                                                                   

In the present section, we rely on an alternative approach to the Riemann problem due to 
\auth{Bianchini and Bressan} \cite{BianchiniBressan,Bianchini,Bianchini2}
and we provide a second proof of the inner speed variation estimate already derived in the previous section. 
 
Throughout, we consider a strictly hyperbolic system of conservation laws \eqref{Class.1} 
and use the notation of Section~\ref{Riem-0}. Note that the parametrization (often denoted by the letter $s$) 
along the wave curves will be here different from the one in the previous section. The sup norm in $\RN$ is used 
throughout the present section. 

The following framework is based on a prescribed family of ``traveling waves''; for the motivations of which  
we refer the reader to \cite{LeFloch1,DLM,LeFloch5,BianchiniBressan,Bianchini2} 
and to the discussion in Section~\ref{NC-0} below. 
Consider a given family of smooth, vector-valued maps $\tr_j = \tr_j(u,v_j,\sigma_j)$ 
for $(u,v_j,\sigma_j) \in \RN \times \RR \times \RR$ 
(presumably associated with viscous profiles of a regularized version of \eqref{Class.1}). 
Given a base point \emph{which can be assumed to be the origin in $\RN$}, and setting $l_j^0 := l_j(0)$, etc., 
we can normalize the vectors $\tr_j$ so that 
$l_j^0 \cdot \widetilde{r}_j =1$. 

From $\tr_j$ we can determine the speed functions $\tlam_j$ 
$$
\tlam_j(u, v_j, \sigma_j) : = l_j^0 \cdot A(u) \, \tr_j(u, v_j, \sigma_j), 
$$
and we assume as in \cite{BianchiniBressan} that, for some constant $C_0>0$ 
\be
\tr_j(u, 0, \sigma_j) = r_j (u), 
\qquad 
\big| \frac{\del \tr_j}{\del \sigma_j} (u, v_j, \sigma_j) \big| \leq C_0 \, |v_j|, 
\label{PropBianchiniBressan}
\ee  
and 
\be
\label{EstimeeBianchiniBressan}
\big| \frac{\del \tlam_j}{\del v_j} (u,v_j,\sigma_j) \big| \leq C_0 \, |u|, 
\qquad 
\big| \frac{\del \tlam_j}{\del \sigma_j} (u,v_j,\sigma_j) \big| \leq C_0 \, |u| \, |v_j| .
\ee 
In particular, it follows that $\tlam_j(u, 0, \sigma_j) = \lam_j(u)$.

Then, the wave curves associated with the system \eqref{Class.1} are constructed in the following way.
Fix some small $\delta_1 < \delta_0$ and, for $j=1, \ldots, N$, $s \in [0,\delta_1)$, and $u_- \in \Bcal_{\delta_0}$, 
define a family of curves $\Gamma_j(s; u_-)$, issuing from $u_-$ and of ``length'' $s$, as follows: 
$$
\aligned 
\Gamma_j(s; u_-) := \Big\{ 
& \gam(\tau) =(u(\tau),v_j(\tau),\sigma_j(\tau)) \in Lip([0,s]; \RR^{N+2}) \, / \,  
\\
&  u(0)=u_-, \quad \mu_j(u(\tau)) = \mu_j(u_-) + \tau, 
\\  
& v_j(0)=0, \ |v_j(\tau)| \leq \delta_1, 
\quad  |\sigma_j(\tau) -  \lam_j^0| \leq 2C_0 \, \delta_1 
\Big\},
\endaligned
$$ 
where $\mu_j(u) := l_j^0 \cdot u$.
To any curve $\gam \in \Gamma_j(s;u_-)$ we associate its {\bf $j$-flux function}
$$
\tf_j[\gam](\tau) := \int_0^\tau \tlam_j(u, v_j, \sigma_j)(\tau') \, d\tau'. 
$$
Define also the nonlinear operator 
$\Tcal_j : \gam \mapsto (\widehat u, \widehat{v}_j, \widehat{\sigma}_j)$ on $\Gamma_j(s,u_-)$
by 
\be
\label{OperateurBB}
\aligned 
& \widehat{u}(\tau) = u_- + \displaystyle \int_0^\tau \tr_j (u, v_j, \sigma_j)(\tau') \, d\tau', 
\\
&\widehat{v}_j(\tau) = \tf_j[\gam](\tau) - \conv_{[0,s]} \tf_j[\gam](\tau), 
\\ 
&\widehat{\sigma}_j(\tau) = \frac{d}{d\tau} \Big( \conv_{[0,s]}  \tf_j[\gam](\tau) \Big). 
\endaligned 
\ee

Define a distance between $\gam, \gam' \in \Gamma_j(s,u_-)$ by 
\be
\label{DefD}
D(\gam,\gam') := \delta_1 \| u-u' \|_{L^\infty} + \| v_j - v'_j \|_{L^1} 
+ \| v_j \sigma_j - v'_j \sigma'_j \|_{L^1}. 
\ee
It was established in \cite{Bianchini} (Proposition 3.2 therein) that, 
for $\delta_1$ suitably small, the operator $\Tcal_j$ is a contraction with constant $1/2$
in the metric space $\Gamma_j(s;u_-)$ endowed with the distance $D$. To any $s$ and 
$u_-$ we can thus associate a unique curve $\gam^\star= \gam^\star_{s, u_-}$ that satisfies  
$\Tcal_j(\gam^\star) = \gam^\star$. 

This allows us to introduce the {\bf wave curve} $\psi_j= \psi_j(s;u_-)$ by keeping, for every $s>0$,
the end-point of the curve at $\tau =s$, only, i.e. 
$$
\psi_j(s;u_-) := \gam^\star_{s;u_-}(s).
$$
For every $s>0$, the state $u_-$ is connected to the right-hand state $u_+:=\psi_j(s;u_-)$ 
along the curve $\gam^\star = (u^\star,v_j^\star,\sigma_j^\star)$ 
(associated with the given $s,u_-$) 
and we can define the corresponding $j$-wave fan: 
\be
\nonumber  
u(t,x) := 
\begin{cases}
u_-,     & x/t \leq \sigma_j(0), 
\\
u(\tau),   &  x/t = \sigma_j(\tau), 
\\
u_-,       &  x/t \geq \sigma_j(s).
\end{cases}
\ee
When $s<0$ a similar construction is done using the concave hull (denoted below by $\conc$) instead of the convex hull.

The wave curves are Lipschitz continuous and, by the implicit function theorem, 
any {\bf Riemann problem} $(u_l,u_r)$ can be solved uniquely, by combining wave curves together, i.e. 
$$ 
u_r = \Psi(\s; u_l) := \psi_N(s_N; \psi_{N-1}(s_{N-1}, \ldots \psi_1(s_1;u_l))), 
$$
where we use the notation $\s : = (s_1, \ldots, s_N)$.

To estimate the interaction between waves we introduce Bianchini-Bres\-san's {amount of interaction} 
\cite{Bianchini} 
as follows. First of all, consider two waves of the same family 
\be
\label{4etats}
u_+ = \psi_j(s, u_-), \qquad  u'_+ = \psi_j(s';u'_- ),
\ee
together with the corresponding curves $\gam$ and $\gam'$ and the corresponding 
$j$-flux  $\tf_j= \tf_j [\gam]$ 
and $\tf'_j = \tf_j [\gam']$. 
Define the real $\IcalBB = \IcalBB (u_-,u_+; u'_-,u'_+)$, as follows, where 
for definiteness we assume that $s >0$. (When $s<0$ one should replace all convex/concave hulls by concave/convex hulls.) 
\begin{itemize}
\item If $s,s' >0$, we set 
\be
\label{fUf}
\tf_j \cup  \tf'_j := 
\begin{cases} 
\tf_j(\tau),      &  \tau \in [0,s], 
\\
\tf_j(s) + \tf'_j(\tau-s),      & \tau \in [s,s+s'],
\end{cases} 
\ee
and 
$$
\aligned 
\IcalBB(u_-,u_+; & u'_-,u'_+)
\\  
:= & \int_0^s \left| \conv_{[0,s]} (\tf_j\cup \tf'_j)(\xi) - \conv_{[0,s+s']} (\tf_j\cup \tf'_j)(\xi) \right|  
\\
& +
\int_{s}^{s+s'} \left| \conv_{[s,s+s']} (\tf_j\cup \tf'_j)(\xi) - \conv_{[0,s+s']} (\tf_j\cup \tf'_j)(\xi) \right|,
\endaligned
$$ 
\item If $-s \leq s' <0$, we set 
$$
\aligned  
\IcalBB(u_-,u_+;u'_-,u'_+) :=
& \int_0^{s+s'} \left| \conv_{[0,s]} \tf_j(\xi) - \conv_{[0,s+s']} \tf_j(\xi) \right|  
\\
& + \int_{s+s'}^{s} \left| \conv_{[0,s]} \tf_j(\xi) - \conc_{[s+s',s]}  \tf_j(\xi) \right|,
\endaligned 
$$ 
\item If $ s' < - s <0$, we set 
$$
\aligned 
\IcalBB(u_-,u_+; u'_-,u'_+) :=
& \int_{s'}^{-s} \left| \conc_{[s',0]} \tf'_j(\xi) - \conc_{[s',-s]} \tf'_j(\xi) \right|  
\\
& +
\int_{-s}^{0} \left| \conc_{[s',0]} \tf'_j(\xi) - \conv_{[-s,0]}  \tf'_j(\xi) \right|.  
\endaligned 
$$
\end{itemize}

More generally, define the amount of interaction associated with two Riemann problems 
$u_r = \Psi(\s;u_l)$, $u'_r = \Psi(\s';u'_l) $ by 
$$
\IcalBB(u_l, u_r; u'_l,u'_r) := \sum_{i<j} |s_i \, s_j'| + \sum_j \IcalBB(u_-, u_+;u'_-,u'_+), 
$$
where the second sum is over all $j$-waves $(u_-,u_+)$ and $(u'_-,u'_+)$ in the Riemann problems  
$(u_l, u_m)$ and $(u_m, u_r)$, respectively. Now the interaction functional $\QBB$ is defined as  
$$
\QBB(u_l,u_r,u'_l,u'_r) 
= \sum_{i > j}  | s_i  \, s'_j | + \sum_j \int_0^{s_j} \int_0^{s'_j} | \sigma_j(\tau) - \sigma'_j(\tau) |.
$$
It will be convenient also to set $\QBB(u_l, u_m,u_r) := \QBB(u_l, u_m; u_m,u_r)$ for all $u_l,u_m,u_r$.

Finally, we will need another distance between curves. 
Consider two curves $\gam \in \Gamma_j(s_j;u)$ and $\gam' \in \Gamma_j(s'_j;u')$, and 
restrict attention to the case where $s_j$ and $s'_j$ have the same sign. 
If, for instance, both of them are positive we set $\underline{s}_j := \min(s_j,s'_j)$ and 
$$ 
P(\gam,\gam') := D(\gam_{|[0,\underline{s}_j]}, \gam'_{|[0,\underline{s}_j]}) + |s_j - s'_j |,
$$
where the distance $D$ introduced earlier is extended in an obvious way to any two curves with the
same length but not the same base point.  


\subsection{Interactions estimates on the inner speed variation}                                                                                                                                                                                                                                                                                        

Recall from \cite{Bianchini} the following version of Glimm's interaction estimates.  

\begin{theorem} 
\label{TheoBianchini1}
1. Consider any three states $u_l, u_m, u_r$ such that
$u_m = \Psi(\s_{lm};u_l)$, $u_r = \Psi(\s_{mr};u_m)$, and $u_r = \Psi(\s_{lr};u_l)$.  
Then, the corresponding wave strength vectors satisfy 
\be
\label{EstB1}
|\s_{lr} - \s_{lm} - \s_{mr} | \lesssim \IcalBB(u_l, u_m, u_r). 
\ee 

2.  Let $\gam_j^{lm}$, $\gam_j^{mr}$, and $\gam^{lr}_j$ be the curves associated with the $j$-waves 
in the Riemann solutions under consideration. Define $\hgam_j$ by 
$$
\hgam_j := 
\begin{cases}
\gam_j^{lm} \cup \gam_j^{mr}  & \text{ on } [0,s^{lm}_j+s_j^{mr}],Ê\text{ if } s^{lm}_j, s^{mr}_j \geq 0, 
\\
\gam_{j |[0,s^{lm}_j+s^{mr}_j]}^{lm} & \text{ on } [0,s^{lm}_j+s^{mr}_j],Ê\text{ if }  -s^{lm}_j  < s^{mr}_j < 0, 
\\
\gam_{j|[s^{mr}_j,-s^{lm}_j]}^{mr} (\cdot - s^{lm}_j) & \text{ on } [s_j^{lm}+s_j^{mr},0], \text{ if }  s^{mr}_j  < -s^{lm}_j < 0,
\end{cases} 
$$
where $\gam_j^{lm} \cup \gam_j^{mr}$ is defined by 
\be
\label{gammaUgamma}
\gam_j^{lm} \cup \gam_j^{mr} := 
\begin{cases} 
\gam_j^{lm}(\tau),        & 0 \leq \tau \leq s_j^{lm}, 
\\
\gam_j^{lm}(s_j^{lm}) + \gam_j^{mr}(\tau - s_j^{lm}),     &   s_j^{lm} \leq \tau \leq s_j^{lm} + s^{mr}_j,
\end{cases} 
\ee
(the definition when $s^{lm}_j, s^{mr}_j$ are negative being similar).
Then, the curves satisfy the following interaction estimates:  
$$
\sum_j P(\hgam_j, \gam^{lr}_j ) \lesssim \IcalBB(u_l, u_m, u_r). 
$$  

3. Consider a piecewise constant function with small total variation $u:\RR \rightarrow \RR^N$. 
Let $u'$ be the function obtained from $u$ by replacing two consecutive Riemann problems 
$(u_k,u_{k+1})$ and $(u_{k+1},u_{k+2})$ by $(u_k,u_{k+2})$. Then, for some $c>0$ we have the estimate 
\be
\label{EstGlimmBianchini}
\QBB(u') \leq \QBB(u) - c \, \IcalBB(u_k,u_{k+1},u_{k+2}).
\ee

4. In the situation above, let $Q'$ be the interaction amount associated with $u$ but 
with the $i$-th curve $\gamma^i_k$ in the Riemann problem $(u_k,u_{k+1})$ being replaced 
by $\gamma'_i$ corresponding to another problem $u_1'=\psi_i(u_0')$. This modifies the $Q$ in the following way:
\be
\label{EstGlimmBianchini2}
|Q' - \QBB(u)| \lesssim P(\gamma_i^k,\gamma'_i).
\ee

\end{theorem}


We now turn to the investigation of the inner speed variation $\vartheta_j$ (defined in the previous section) which here takes the following form: 
given $u_r=\psi_j(s,u_l)$ with corresponding curve $\gamma=(u,v_j,\sigma_j)$:
$$
\vartheta_j(u_l;u_r)= \sigma_j(s)-\sigma_j(0).
$$
We provide a new proof of the estimate already derived in the previous section.

\begin{proposition} {\rm (Inner speed variation -- same family.)} 
\label{Prop:Interaction1} 
Consider two $j$-wave fans 
$$ 
u_m = \psi_j(s_{lm}; u_l),   \qquad   u_r = \psi_j(s_{mr}; u_m), 
$$
and denote by $\gam_{lm} = (u^{lm}, v_i ^{lm}, \sigma_i^{lm})$ and 
$\gam_{mr} = (u^{mr}, v_i ^{mr}, \sigma_i^{mr})$ the associated curves. 
Set $u_r' := \psi_j(s_{lm} + s_{mr}; u_l)$ and denote by $\gam' := (u', v_j', \sigma_j')$ the 
associated curve. Then, the inner speed variation of the $j$-wave fan $(u_l,\tu_r)$ 
satisfies the following estimates (assuming $s_{lm} \geq 0$): 
\begin{itemize}
\item {\it Monotone case : $0 \leq s_{lm}, s_{mr}$.} 
\be
\aligned 
\vartheta_j(u_l,u_r') 
\leq & 
\max\big( \vartheta_j(u_l,u_m), \vartheta_j(u_m,u_r) \big)  
\\
& + \big(\sigma^{mr}_j(0) - \sigma^{lm}_j(0) \big)_+ + \GO(1) \, \IcalBB(u_l,u_m,u_r).
\endaligned 
\label{OscDiminue1}
\ee
\item {\it Non-monotone case (I) : $0 \leq s_{lm}+  s_{mr} \leq s_{lm}$.} 
\be
\label{OscDiminue2}
\vartheta_j(u_l,u_r')
\leq \vartheta_j(u_l,u_m) + \GO(1) \, |s_{mr}|. 
\ee
\item {\it Non-monotone case (II) : $s_{lm}+  s_{mr} \leq 0 \leq s_{lm}$.} 
\be
\label{OscDiminue3}
\vartheta_j(u_l,u_r')
\leq \vartheta_j(u_m,u_r) + \GO(1) \, s_{lm}. 
\ee
\end{itemize}
Moreover, completely similar estimates hold in the case $s_{lm} \leq 0$.
\end{proposition}


\begin{proposition}  {\rm (Inner speed variation -- artificial wave fronts.)} 
\label{Prop:Interaction3}
Consider three states $u_l, u_m,u_r$ satisfying solely 
$$
u_r = \psi_j(s; u_m), 
$$
and let $\gam=(u, v_j, \sigma_j)$ be the curve associated with the wave $(u_m,u_r)$. Introduce $u_m' := \psi_j(s;u_l)$
and the corresponding curve $\gam' := (u', v'_j, \sigma'_j)$. Then, it holds 
$$
\vartheta_j(u_l,u_m')  = \vartheta_j(u_m,u_r) + \GO(s) \, |u_m-u_l|. 
$$
A completely similar statement holds when $(u_l, u_m)$ is a $j$-wave and $(u_m,u_r)$ is an arbitrary jump discontinuity. 
\end{proposition}


In fact, from the above statement, we can also deduce the desired estimate for the interaction of waves of different families. 

\begin{proposition}  {\rm (Inner speed variation -- different families.)} 
\label{Prop:Interaction2}
Consider three states $u_l, u_m,u_r$ such that, for $i \not= j$, 
$$ 
u_m = \psi_j(s_{lm};u_l) , \qquad u_r = \psi_i(s_{mr};u_m),
$$
and denote by $\gam^{lm}=(u^{lm}, v_j^{lm}, \sigma_j^{lm})$ and $\gam^{mr}=(u^{mr}, v_i^{mr}, \sigma_i^{mr})$ 
the curves associated with the wave $(u_l, u_m)$ and $(u_m,u_r)$, respectively. 
Set also $\tu_m:= \psi_i(s_{mr};u_l)$ and $\tu_r:= \psi_i(s_{lm}; \tu_m)$,
and denote by ${\gam}'^{lm}:=({u}'^{lm},{v'^{lm}_j}, {\sigma'^{lm}_j})$, 
${\gam}'^{mr}:=({u}'^{mr},{v'^{mr}_i}, {\sigma'^{mr}_j})$
the corresponding curves. Then, one has 
$$
\aligned 
& \vartheta_j(\tu_m,\tu_r)  = \vartheta_j(u_l,u_m) + \GO(s_{lm} \, s_{mr}), 
\\
& \vartheta_j(u_l,\tu_m)  = \vartheta_i(u_m,u_r)  + \GO(s_{lm} \, s_{mr}).
\endaligned
$$   
\end{proposition}

\begin{corollary}
\label{CorRegTheta} 
Consider four states $u_0,u'_0, u_1, u'_1$ such that
$$ 
u_1 = \psi_j(s;u_0),    \qquad  u'_1 = \psi_j(s';u'_0) 
$$
with $0<s< s'$. Denote by $\gam=(u, v_j, \sigma_j)$ and $\gam'=(u', v'_j, \sigma'_j)$ the 
curves associated with the waves $(u_0,u_1)$ and $(u'_0, u'_1)$, respectively. Then, one has 
\be
\label{RegTheta}
| {\sigma'_j}(s') - {\sigma'_j}(0) - \sigma_j(s) + \sigma_i(0) | \lesssim
|u_0 - u'_0 | + |s'-s|.
\ee
\end{corollary}

\

Clearly, Proposition~\ref{Prop:Interaction2} follows from the fact that 
$$
s \lesssim | \psi_j(s; u) - u| \lesssim s.
$$ 
We postpone the proof of the other statements to Subsection~\ref{PreuvesEI}. In the next subsection, 
we derive some additional properties of the Riemann solver constructed via \eqref{OperateurBB}, 
which will be useful in our proofs but are also of independent interest. 


\subsection{Additional properties of the Riemann solver} 

To derive the inner speed variation estimates stated earlier, it will be convenient 
to decompose the problem of finding a fixed point of \eqref{OperateurBB} in two steps: on one hand, 
finding the geometric wave curves
(that is the component $u$) and, on the other hand, finding the speed of propagation of the wave fan
(that is the components $(v_j, \sigma_j)$). This is the subject of the following proposition. 
Consider the set 
$$
\aligned 
\Upsilon_j(s; u_-) := \Big\{ 
& (v_j(\cdot),\sigma_j(\cdot)) \in Lip([0,s]; \RR^2) \, \big/
\\
& v_j(0) = 0, \ |v_j(\tau)| \leq \delta_1, \ |\sigma_j(\tau) -  \lam_j^0| \leq 2C_0 \,\delta_1 \Big\},
\endaligned 
$$
endowed with the sup norm $\| (v_j,\sigma_j) \|_\infty := \| v_j \|_\infty + \| \sigma_j \|_\infty$. 
Given $u^* \in Lip([0,s];\RN)$ satisfying $u^*_j(\tau)  = u_- + \tau$ and $|u^*_j(\tau) - u_- | \leq \tau$, 
consider the operator 
$$
\Omega^{u^*}_s : \ (v_j,\sigma_j) \mapsto
\begin{cases}
&   \hat{v}_j(\tau) : = \tf_j [u^*,v_j,\sigma_j](\tau) - \conv_{[0,s]} \tf_j[u^*,v_j,\sigma_j](\tau), 
\\
&   \hat{\sigma}_j(\tau) = \frac{d}{d\tau} \Big( \conv_{[0,s]} \tf_j [u^*,v_j,\sigma_j] \Big)(\tau).
\end{cases} 
$$

\begin{lemma}
\label{Prop:AutreOperateur}
Fix $u^* \in Lip([0,s];\RN)$ satisfying $u^*_j(\tau)  = u_- + \tau$ and $|u^*_j(\tau) - u_- | \leq \tau$. 
Provided the range of $u^*$ is included in a sufficiently small neighborhood of $0$, 
the operator $\Omega^{u^*}_s$ admits a unique fixed point in $\Upsilon_j(s;u_-)$, which will be 
denoted by $(V^{u^*}_s, \Sigma^{u^*}_s)$. Moreover, for some $C_1>0$, 
$$
\| (V^{u^1}_s, \Sigma^{u^1}_s) - (V^{u^2}_s, \Sigma^{u^2}_s) \|_\infty 
\leq C_1 \, \| u^1 - u^2 \|_\infty.
$$
\end{lemma}

\begin{proof}  
We will prove that if the range of $u^*$ is included  in a sufficiently small 
neighborhood of $0$, then the map $\Omega^{u^*}_s$ is a contraction with constant $1/2$. Indeed, one has 
\begin{eqnarray*}
\| \Omega^{u^*}_s(v_1,\sigma_1) -  \Omega^{u^*}_s(v_2,\sigma_2) \|_\infty & \lesssim &
\| \frac{d}{d\tau} \tf_j(\cdot,u^*, v_1,\sigma_1) - \frac{d}{d\tau} \tf_j(\cdot,u^*, v_2,\sigma_2) \|_\infty 
\\
 & \lesssim &
\| \tlam_j(u^*, v_1,\sigma_1) - \tlam_j (u^*, v_2,\sigma_2) \|_\infty.
\end{eqnarray*}
With \eqref{EstimeeBianchiniBressan} this yields
\begin{multline*}
\| \Omega^{u^*}_s(v_1,\sigma_1) -  \Omega^{u^*}_s(v_2,\sigma_2) \|_\infty  \\
\lesssim  C_0  \|u^* - u_0\|_\infty (1 + \|v_j\|_\infty) \| (v_1,\sigma_1) - (v_2,\sigma_2) \|_\infty,
\end{multline*}
 This proves the first part of the lemma.

To establish the second part, we use the above contraction property and
write 
\begin{eqnarray*}
\| (V_s^{u_1},\Sigma_s^{u_1}) -  (V_s^{u_2},\Sigma_s^{u_2}) \|_\infty & \lesssim &
\| \Omega^{u_1}_s (V_s^{u_2},\Sigma_s^{u_2}) -  (V_s^{u_2},\Sigma_s^{u_2}) \|_\infty \\
& \lesssim & \| \tlam_j (u_1,V_s^{u_2},\Sigma_s^{u_2}) - \tlam_j (u_2,V_s^{u_2},\Sigma_s^{u_2}) \|_\infty \\
& \lesssim &  \| u_1 - u_2 \|_\infty,
\end{eqnarray*}
where  the regularity of the function $\tlam_j$ has been used. 
This completes the proof. 
\quad \end{proof}

We now state some properties of the curve describing the wave fan, under the assumption that the $u$-component
is already known.  Introduce the translation operator $\tau_s : f \mapsto \tau_s f$, defined by 
$$
\tau_s(f)(y) := f(y-s). 
$$
When a function $u^*$ is defined on an interval larger than $[0,s]$, in order 
to simplify the notation, we simply write $(V_{s}^{u^*},\Sigma_{s}^{u^*})$ instead 
of  $(V_{s}^{u^*_{|[0,s]}},\Sigma_{s}^{u^*_{|[0,s]}})$.

\begin{lemma} {\rm (A splitting property.)} 
\label{Lem:AutreOperateur}
Let $0<s_1<s_2$ and $u^* \in Lip([0,s_2],\RN)$ be as in Lemma~\ref{Prop:AutreOperateur}.
If
$$
V_{s_2}^{u^*}(s_1)=0,
$$ 
then $(V_{s_2}^{u^*},\Sigma_{s_2}^{u^*})$ coincides with $(V_{s_1}^{u^*}, \Sigma_{s_1}^{u^*})$ 
on the interval $[0,s_1]$, and coincides with the translate 
$\tau_{s_1}(V_{s_2- s_1}^{\tau_{-s_1}u^*}, \Sigma_{s_2-s_1}^{\tau_{-s_1}u^*})$ on 
the interval $[s_1,s_2]$.
\end{lemma}

\begin{proof}   
It suffices to check that $(V_{s_2}^{u^*},\Sigma_{s_2}^{u^*})_{|[0,s_1]}$ 
and $\tau_{-s_1} (V_{s_2}^{u^*},\Sigma_{s_2}^{u^*})_{|[s_1,s_2]}$ are fixed points of 
the operators $\Omega^{u^*_{|[0,s_1]}}_{s_1}$ and
$\Omega^{\tau_{-s_1}u^*_{|[s_1,s_2]}}_{s_2-s_1}$, respectively. 
But this property is a direct consequence of the following fact: 
for every Lipschitz continuous function $g:[0,s_2] \rightarrow \RR$ and $s_1 \in [0,s_2]$ 
$$ 
(\conv_{[0,s_2]} g)(s_1)=g(s_1) \ \Longrightarrow \conv_{[0,s_2]} g \equiv 
\left\{ \begin{array}{l}
\conv_{[0,s_1]} g \text{ on } [0,s_1], \\
\conv_{[s_1,s_2]} g \quad \text{ on } [s_1,s_2].
\end{array} \right.
$$  \quad \qed  
\end{proof}

Now, we prove (compare with \cite[Lemma 3.3]{Bianchini2}):

\begin{lemma} {\rm (A monotonicity property.)} 
\label{Prop:AutreOperateur2} 
For any $0<s_1<s_2$ and $u^* \in Lip([0,s_2],\RN)$ as in Lemma~\ref{Prop:AutreOperateur}, it holds:  
\be
\label{PseudoLiu}
\Sigma_{s_1}^{u^*} \geq \Sigma_{s_2}^{u^*} \text{ on } [0,s_1].
\ee
\end{lemma}

\begin{proof}  Set 
$$
s^\square := \max \Big\{ s \in [0,s_1], \ V_{s_2}^{u^*}(s) =0 \Big\}, 
$$
and note that the function $\conv_{[0,s_2]} f_j(\tau,u^*,V_{s_2}^{u^*},\Sigma_{s_2}^{u^*})$ is affine 
on the interval $[s^\square,s_1]$. 
Suppose first that $s^\square >0$. We have $V_{s_2}^{u^*}(s^\square)=0$  and, by 
Lemma~\ref{Lem:AutreOperateur}, the curves 
$(V_{s_1}^{u^*}, \Sigma_{s_1}^{u^*})$ and $(V_{s_2}^{u^*}, \Sigma_{s_2}^{u^*})$ coincide on 
the interval $[0,s^\square]$.

Hence, it is clear that \eqref{PseudoLiu} is valid in $[0,s^\square]$, so that on the interval $[s^\square,s_1]$ we 
obtain 
$$ 
\Sigma_{s_1}^{u^*}(s) 
\geq \Sigma_{s_1}^{u^*}(s^\square) 
 = \Sigma_{s_2}^{u^*}(s^\square) = \Sigma_{s_2}^{u^*}(s).
$$ 
Recall here that, since $\gam$ is Lipschitz continuous, the function $\widetilde{f}_j(\gam)$ is $W^{2,\infty}$.

Suppose next that $s^\square =0$. Because $\conv_{[0,s_2]} f_j(\tau,u^*,V_{s_2}^{u^*},\Sigma_{s_2}^{u^*})$ is 
affine in $[0,s_1]$, it suffices to establish \eqref{PseudoLiu} at $s=0$. 
We distinguish between two cases:
\begin{itemize}
\item either $0$ is an accumulation point of values $s$ such that $V_{s_1}^{u^*}(s)=0$. In this case, 
we have 
$$
\aligned 
\Sigma_{s_1}^{u^*}(0) 
& = \frac{d}{d\tau} \widetilde{f}_j (0,u^*, V_{s_1}^{u^*}, \Sigma_{s_1}^{u^*})
& = \lam_j(u^*(0)),
\endaligned 
$$
and, on the other hand, 
$$
\aligned 
\Sigma_{s_2}^{u^*}(0) 
 \leq \lam_j(u^*(0)) 
 =\frac{d}{d\tau} \widetilde{f}_j (0,u^*, V_{s_2}^{u^*}, \Sigma_{s_2}^{u^*}).
\endaligned 
$$
\item or $s^\triangle$ is the smallest positive $\tau$ such that $V_{s_1}^{u^*}(\tau)=0$. 
In that case, 
both $\Sigma_{s_1}^{u^*}$ and $\Sigma_{s_2}^{u^*}$ are constant in the interval $[0, s^\triangle]$, 
the former coinciding with $\Sigma_{s^\triangle}^{u^*}$ in $[0, s^\triangle]$ thanks to Lemma~\ref{Lem:AutreOperateur}.
Now, let us use \eqref{EstimeeBianchiniBressan} and recall that both $\Sigma_{s^\triangle}^{u^*}$ and $\Sigma_{s_2}^{u^*}$ 
are constant in $[0,s^\triangle]$. Recall also that both $V_{s^\triangle}^{u^*}$ and $V_{s^\triangle}^{u^*}$ 
are linear in $[0,s^\triangle]$ and $V_{s^\triangle}^{u^*}(s^\triangle)=V_{s^\triangle}^{u^*}(0)=V_{s_2}^{u^*}(0)=0$. 
We then find: 
$$
\aligned 
& \| \widetilde{f}_j(u^*,V_{s^\triangle}^{u^*}, \Sigma_{s^\triangle}^{u^*})
 - \widetilde{f}_j(u^*,V_{s_2}^{u^*}, \Sigma_{s_2}^{u^*}) \|_{L^\infty(0,s^\triangle)}  \\
 & \leq \| \tlam_j(u^*,V_{s^\triangle}^{u^*}, \Sigma_{s^\triangle}^{u^*})
 - \tlam_j(u^*,V_{s_2}^{u^*}, \Sigma_{s_2}^{u^*}) \|_{L^1(0,s^\triangle)}  \\
& \lesssim \| u^*\|_\infty \| (V_{s^\triangle}^{u^*}, \Sigma_{s^\triangle}^{u^*})
 - (V_{s_2}^{u^*}, \Sigma_{s_2}^{u^*}) \|_{L^1(0,s^\triangle)}  
 \leq  \frac{1}{2} V_{s_2}^{u^*}(s^\triangle), 
\endaligned 
$$
at least if $u^*$ remains within a small neighborhood of the base point $0$. This yields us 
$$ 
 \widetilde{f}_j(u^*,V_{s^\triangle}^{u^*}, \Sigma_{s^\triangle}^{u^*})(s^\triangle) 
 \geq \conv_{[0,s_2]} \widetilde{f}_j(u^*,V_{s_2}^{u^*}, \Sigma_{s_2}^{u^*})  (s) + \frac{1}{2} V_{s_2}^{u^*}(s^\triangle), 
$$
which leads to \eqref{PseudoLiu} and completes the proof of Lemma~\ref{Prop:AutreOperateur2}. 
\quad 
\end{itemize} 
\end{proof}

From Lemma~\ref{Prop:AutreOperateur2} we deduce:

\begin{lemma}
\label{Prop:AutreOperateur3} 
Let $s_1,s_2>0$, $\ts \in (0,s_1)$, and $u^* \in Lip([0,s_1+s_2];\RN)$. If
\be
\label{AnnuleMilieu}
V^{u^*}_{s_1+s_2}(\ts)=0,
\ee
then 
\be
\label{EntreDeux}
V^{u^*}_{s_1}(\ts)=0. 
\ee
Similarly, if \eqref{AnnuleMilieu} holds for some $\ts \in (s_1,s_1+s_2)$, then 
$V^{\tau_{-s_1}u^*}_{s_2}(\ts-s_1)=0$.
\end{lemma}

\begin{proof}  
We only prove the first statement since the second one is similar. 
We introduce 
$$
\check{u}:= \tau_{-\ts} u^*_{|[\ts,s_1+s_2]}, 
\qquad 
V^{\check{u}}_{s_1+s_2-\ts}, 
\qquad 
\Sigma^{\check{u}}_{s_1+s_2-\ts}.
$$  
Using Lemma~\ref{Lem:AutreOperateur}, we obtain 
$$
\aligned 
\Sigma^{u^*}_{\ts} (\ts) 
& = \Sigma^{u^*}_{s_1+s_2} (\ts)
\\
& = \frac{d}{d\tau} \conv_{[0,s_1+s_2]} f_j(\ts, u^*, \Sigma^{u^*}_{s_1+s_2}, V^{u^*}_{s_1+s_2}).  
\endaligned
$$
Using $V^{u^*}_{s_1+s_2}(\ts)=0$, we then have 
$$
\aligned  
\Sigma^{u^*}_{\ts} (\ts) 
& = \frac{d}{d\tau}  f_j(\ts, u^*, V^{u^*}_{s_1+s_2}, \Sigma^{u^*}_{s_1+s_2})
\\
& = \tlam_j(u^*(\ts), 0, \Sigma^{u^*}_{s_1+s_2}(\ts)) = \lam_j(u^*(\ts)).
\endaligned
$$
Now, from Lemma~\ref{Prop:AutreOperateur2} and \ref{Lem:AutreOperateur}, we deduce
\begin{eqnarray*}
\tau_{\ts}\Sigma_{s_1-\ts}^{\check{u}} (\ts)  
\geq  \tau_{\ts}\Sigma_{s_1+s_2 -\ts}^{\check{u}} (\ts)
=  \Sigma_{s_1+s_2 }^{u^*} (\ts),
\end{eqnarray*}
It follows that
\begin{eqnarray*}
\tau_{\ts}\Sigma_{s_1-\ts}^{\check{u}} (\ts) \geq \lam_1(u^*(\ts)).
\end{eqnarray*}

But, on the other hand,
$$
\aligned 
\tau_{\ts}\Sigma_{s_1-\ts}^{\check{u}} (\ts) 
& = \frac{d}{d\tau} \conv_{[0,s_1-\ts]} f_j (\check{u}, V_{s_1-\ts}^{\check{u}}, \Sigma_{s_1-\ts}^{\check{u}}) (0) 
\\ 
& \leq \frac{d}{d\tau} f_j (\check{u}, V_{s_1-\ts}^{\check{u}}, \Sigma_{s_1-\ts}^{\check{u}}) (0) 
\\
& = \tlam_1(u^*(\ts),0,\Sigma_{s_1-\ts}^{\check{u}}) = \lam_1(u^*(\ts)).
\endaligned
$$ 
Hence, we obtain  
\begin{eqnarray*}
\tau_{\ts}\Sigma_{s_1-\ts}^{\check{u}_{|[0,s_1-\ts]}} (\ts) = \lam_1(u^*(\ts)) = 
\Sigma_{\ts}^{u^*_{|[0,\ts]}}(\ts).
\end{eqnarray*}
But, we have also
\begin{eqnarray*}
\tau_{\ts} V_{s_1-\ts}^{\check{u}_{|[0,s_1-\ts]}} (\ts) = 0 = 
V_{\ts}^{u^*_{|[0,\ts]}}(\ts),
\text{ and }
\tau_{\ts}\check{u}_{|[0,s_1-\ts]} (\ts) = u^*(\ts).
\end{eqnarray*}

We introduce the curve $(\Vb,\Sigmab): [0,s_1] \rightarrow \RR^{2}$ given by 
\be
\label{EntreDeuxFort}
(\Vb,\Sigmab) 
\equiv
\begin{cases} 
(V_{\ts}^{u^*}, \Sigma_{\ts}^{u^*}),     &  [0,\ts], 
\\
\tau_{\ts}(V_{s_1-\ts}^{\check{u}}, \Sigma_{s_1-\ts}^{\check{u}}), 
        & [\ts,s_1],
\end{cases} 
\ee
Clearly, $\overline{\gam}$ is Lipschitz continuous. 
To prove \eqref{EntreDeux}, it suffices to prove that $(V_{s_1}^{u^*}, \Sigma_{s_1}^{u^*})=(\Vb,\Sigmab)$,
and hence that $\overline{\gam}$ is a fixed point for $\Omega^{\check{u}}_{s_1}$.

Now the following property is satisfied by any Lipschitz continuous function $g:[0,s_1] \rightarrow \RR$: 
if its convex hulls $h_1:=\conv_{[0,\ts]} g$ and $h_2:=\conv_{[\ts,s_1]} g$ satisfy 
$$
h_1(\ts) = h_2(\ts), \qquad h'_1(\ts)= h'_2(\ts), 
$$
then $\conv_{[0,s_1]} g$ coincides with $h_1$ in $[0,\ts]$ and with $h_2$ in $[\ts,s_1]$. 
Using the definitions of $(V_{\ts}^{u^*}, \Sigma_{\ts}^{u^*})$ and $(V_{s_1-\ts}^{\check{u}}, \Sigma_{s_1-\ts}^{\check{u}})$, 
one deduces $(V_{s_1}^{u^*}, \Sigma_{s_1}^{u^*})=(\Vb,\Sigmab)$ and hence \eqref{EntreDeuxFort}. 
This completes the proof of Lemma~\ref{Prop:AutreOperateur3}. 
\quad  \end{proof}

We conclude this section with a technical observation which will be useful in the course of the proofs
of the above statements. 

\begin{lemma}
\label{ConstLip}
All of the curves $\gam$ obtained by a fixed point argument with the operator \eqref{OperateurBB}
are Lipschitz continuous with a uniform Lipschitz constant.
\end{lemma}

\noindent
{\bf Proof.} 
Consider first the component $(u,v_j)$ part; the desired estimate is  elementary in view of 
the property $\| d (\conv f) \|_\infty \leq \| d f\|_\infty$. 
It remains to consider the $\sigma_j$-part and, precisely, to prove
$$
| \sigma_j(\tau) - \sigma_j(\tau') | \leq C |\tau - \tau'|.
$$

First of all, observe that for every $f \in W^{2,\infty}$ and $0 \leq x \leq y \leq s$ the convex hull $h:=\conv_{[0,s]} f$
satisfies
$$
\osc_{[x,y]} (h') \leq \osc_{[x,y]} (f').
$$
(This is easily check by considering the contact points between $f$ and $h$ in $[x,y]$). Hence,
we have 
$$
| \sigma_j(\tau) - \sigma_j(\tau') |  \lesssim  \osc_{[\tau,\tau']} (\tlam_j (\gam(\tau)), 
$$
and it follows that
$$
\osc_{[\tau, \tau']}(\sigma_j) \leq  \| \del_u \tlam_j \|_\infty \osc_{[\tau, \tau']}(u) 
+ \| \del_v \tlam_j \|_\infty \osc_{[\tau, \tau']}(v_j) 
+ \| \del_\sigma \tlam_j \|_\infty \osc_{[\tau, \tau']}(\sigma_j).
$$
Now, we already have uniform Lipschitz estimates on $(u,v_j)$, and it follows from \eqref{EstimeeBianchiniBressan}
that the coefficient in front of $\osc(\sigma_j)$ is less than $1/2$, which allows us to conclude.
\quad \qed \medskip


\subsection{Proof of the inner speed variation estimates}
\label{PreuvesEI}
 
 We begin with an elementary lemma.

\begin{lemma} 
\label{LemmeEC}
Given $g \in W^{2,\infty}([-M,M])$, there exists a positive constant $C_1=C_1(g)$
(depending upon $M$ and $\sup | {d^2 g/du^2} |$) such that, for all $-M \leq a \leq b \leq c \leq M$, one has
\be
\label{LEC1}
\| (\conv_{[a,c]} g)'_{|[a,b]} -  (\conv_{[a,b]} g)' \|_\infty \leq C_1 \, (c-b), \\
\ee
\be
\label{LEC2}
\| (\conv_{[a,c]} g)'_{|[b,c]} -  (\conv_{[b,c]}g)' \|_\infty \leq C_1 \, (b-a).
\ee
A similar statement stands for the concave hull.
\end{lemma}

\medskip
\noindent
{\bf Proof.} 
We limit our attention to the proof of \eqref{LEC1} in the case of the convex hull since the other inequalities are proven in the same way.  
We use the notation $h_{[a,b]} := \conv_{[a,b]} g$, etc. In view of this definition we note that 
$$ 
h_{[a,c]} \leq h_{[a,b]} \leq g \quad \text{ in } [a,b],
$$
hence we can restrict attention to the case where $a$ and $c$ are the
only contact points between $g$ and $h_{[a,c]}$. Hence one gets
$$
h_{[a,c]}' \equiv \frac{g(c) - g(a)}{c - a}.
$$ 
As $h_{[a,c]}'(a) \leq h_{[a,b]}'(a)$ it is sufficient to prove that
$$ 
h_{[a,b]}'(b) - \frac{g(c) - g(a)}{c-a} \leq \Ocal(c-b).
$$

Suppose that $h_{[a,b]}$ and $g$ are not tangent at $b$, and denote by $d < b$ the ``last'' contact point 
of $h_{[a,b]}$ with $g$. It follows that 
\begin{eqnarray*}
h_{[a,b]}'(b) -  \frac{g(c) - g(a)}{c-a} 
&=&  \frac{g(d) - g(b)}{d-b} -  \frac{g(a) - g(c)}{a-c} 
\\ 
& \leq & \frac{g(d) - g(b)}{d-b} -  \frac{g(d) - g(c)}{d-c}.
\end{eqnarray*}
In the latter inequality we used the fact that the point $(d,g(d))$ is
above the line connecting $(a,g(a))$ to $(c,g(c))$. The latter term in the inequality above is
obviously $\Ocal(c-b)$. 

Now, in the case where $h_{[a,b]}$ and $f$ are tangent at $b$, the difference is
\begin{eqnarray*}
h'(b) -  \frac{g(c)-g(a)}{c-a} &=& g'(b) -  \frac{g(a) - g(c)}{a-c}  
\\
&\leq &  g'(b) -  \frac{g(b) - g(c)}{b-c},
\end{eqnarray*}
where we have used that $(b,g(b))$ is above the line connecting $(a,g(a))$ to
$(c,g(c))$. The right-hand side is also clearly $\Ocal(c-b)$ and this completes the proof of \eqref{LEC1}. 
\quad \qed 
\medskip

\noindent
{\bf Proof of Proposition~\ref{Prop:Interaction1}.} 
To simplify the notation, we use the index $1$ (resp. $2$) to denote the objects ($s$, $\gam=(u,v_i,\sigma_i)$) 
relative to the left-middle pattern (that is, $lm$)  
(resp. the middle-right pattern $mr$). We drop the index $i$ corresponding to the wave family 
(which is the same for all objects considered here). Recall that a prime indicates an object corresponding the outgoing $i$-th wave.

We distinguish between three cases: $s_{1},s_{2}>0$, $-s_{1}<s_{2}<0$ and $s_{2}<-s_{1} < 0$. 
The cases where $s_{1} <0$ are treated similarly using the concave hull. \par
\

\noindent
{\bf Case 1 :} Assume $s_{1},s_{2}>0$. \\
Let us define $\gam_{1} \cup \gam_{2}$ by \eqref{gammaUgamma}, let $u^{\#}$ be the $u$-part of it.
It follows from Theorem~\ref{TheoBianchini1} and Corollary~\ref{CorRegTheta} 
that 
$$ 
D({\gam}',\gam_{1} \cup \gam_{2} ) = \GO(1) \, \IcalBB(u_l,u_m,u_m,u_r),
$$
where the distance $D$ is given by \eqref{DefD}. In particular, 
$$
\| u^{\#} - {u}' \|_\infty = \GO(1) \, \IcalBB,
$$
which in turn yields, using Proposition~\ref{Prop:AutreOperateur},
$$
\| (V_{s_1+s_2}^{u^{\#}},\Sigma_{s_1+s_2}^{u^{\#}}) - (V_{s_1+s_2}^{{u}'},\Sigma_{s_1+s_2}^{{u}'}) \|_\infty 
= \GO(1) \, \IcalBB.
$$
Hence, to establish \eqref{OscDiminue1} it is sufficient to prove
\be
\label{OscDiminueApprox}
\aligned 
& \Sigma_{s_1 + s_2}^{u^{\#}} (s_1 + s_2) - \Sigma_{s_1 + s_2}^{u^{\#}}(0) 
\\
& \leq \max \big( \sigma^{1}(s_1) - \sigma^{1}(0), \sigma^{2}(s_{2}) 
      - \sigma^{2}(0) \big) + \big( \sigma^{2}(0) - \sigma^{1}(0) \big)_+.
\endaligned 
\ee

Let us introduce 
$$
\aligned 
& s_a := \max \Big\{ s \in [0,s_{1}], \  V_{s_1 + s_2}^{u^{\#}}(s)=0 \Big\}, \\
& s_b := \min \Big\{ s \in [s_{1},s_{1}+s_{2}], \  V_{s_1 + s_2}^{u^{\#}}(s)=0 \Big\}. 
\endaligned
$$ 
From Lemmas~\ref{Lem:AutreOperateur} and~\ref{Prop:AutreOperateur3}, we deduce that 
$\Sigma_{s_{1}+s_{2}}^{u^{\#}}$ coincides with $\Sigma^{u_{1}}_{s_{1}}$ in $[0,s_a]$,
and with $\tau_{s_{1}}\Sigma^{u_{2}}_{s_{2}}$ in $[s_b,s_{1}+s_{2}]$. Moreover, $\Sigma_{s_1+s_2}^{u^{\#}}$ 
is constant in $[s_a,s_b]$. Therefore, we see that 
\begin{itemize}
\item if $s_a > 0$ and $s_b < s_1+s_2$, then
\begin{eqnarray*}
\Sigma_{s_1+s_2}^{u^{\#}}(s_1 + s_2) -\Sigma_{s_1+s_2}^{u^{\#}}(0)
&=& \Sigma_{s_2}^{u_2}(s_2) - \Sigma_{s_1}^{u_1}(0) \\
&=& \Sigma_{s_2}^{u_2}(s_2) - \Sigma_{s_2}^{u_2}(0) + \Sigma_{s_2}^{u_2}(0) - \Sigma_{s_1}^{u_1}(0),
\end{eqnarray*}
which yields \eqref{OscDiminueApprox}.

\item if $s_a > 0$ and $s_b = s_1+s_2$ (the equivalent could be done in the case $s_a = 0$ and $s_b < s_1+s_2$), 
then $\Sigma_{s_1+s_2}$ is constant in $[s_a,s_1+s_2]$, hence
\begin{eqnarray*}
\Sigma_{s_1+s_2}^{u^{\#}}(s_1 + s_2) -\Sigma_{s_1+s_2}^{u^{\#}}(0)
&=& \Sigma_{s_1}^{u_1}(s_a) - \Sigma_{s_1}^{u_1}(0),
\end{eqnarray*}
which yields again \eqref{OscDiminueApprox}.
\item if $s_a = 0$ and $s_b = s_1+s_2$ then $\Sigma^{u^{\#}}_{s_1+s_2}$ is constant in $[0,s_1+s_2]$, hence the result is satisfied.
\end{itemize}

{\bf Case 2 :} Assume $-s_1<s_2<0$. \\
Let us write:
$$
\gam^{\#} := \gam_{1|[0,s_1+s_2]}.
$$
It follows from Theorem~\ref{TheoBianchini1} and Corollary~\ref{CorRegTheta} 
that 
$$ 
D(\tilde{\gam},\gam^\#) = \GO(1) \, \IcalBB,
$$
hence
$$\| (V_{s_1+s_2}^{u^{\#}},\Sigma_{s_1+s_2}^{u^{\#}}) - (V_{s_1+s_2}^\tu,\Sigma_{s_1+s_2}^\tu) \|_\infty 
= \GO(1) \, \IcalBB.
$$
We note that, as all the speeds are bounded, $\IcalBB=\GO(s_2)$. Hence, to establish 
\eqref{OscDiminue2} it suffices to show 
$$ 
\sigma^{\#}(s_1+s_2) -  \sigma^{\#}(0) \leq \sigma_1(s_1) -  \sigma_1(0) + \GO(s_2).
$$
As previously, using the contraction property, it is sufficient to prove that
$$ 
\| \Omega^{u^{\#}}_{s_1+s_2} (v_1,\sigma_1)_{|[0,s_1+s_2]} -  (v_1,\sigma_1)_{|[0,s_1+s_2]} \|_\infty = \GO(s_2).
$$
Computing the difference, we are led to establish that
$$ 
\| [\conv_{[0,s_1]} \tf(\gam^1) ]_{|[0,s_1+s_2]} - \conv_{|[0,s_1 +s_2]} \tf(\gam^1) \|_{\Lip} = \GO(s_2).
$$
Using Lemma~\ref{LemmeEC},  
the expression of $\tf_i$ and the smoothness of $\tlam_i$, we see that it is sufficient to prove 
a uniform Lipschitz 
continuous estimate for the curves $\gam$. This is precisely given by Lemma~\ref{ConstLip}. 

\

\noindent
{\bf Case 3 :} Assume $s_2<-s_1<0$. \\
This case is done similarly as the previous one, and  
this completes the proof of Proposition~\ref{Prop:Interaction1}. 
\quad\qed 
\medskip 


\medskip
\noindent
{\bf Proof of Proposition~\ref{Prop:Interaction3}.} 
In the sequel, all objects ($s$, $\gam=(u,v_j,\sigma_j)$) 
without a prime are relative to the incoming $j$-wave, the objects with a prime to the outgoing one. 
Without loss of generality, we assume that $s \geq 0$ and use convex hulls.  
We treat only the case that the given $j$-wave is $(u_m,u_r)$, since the other case is similar. 
Define 
$$
\tu(\cdot):= u(\cdot) - {u_m} + u_l, \qquad 
\tilde{\gam}(\cdot):= (\tu, v, \sigma)(\cdot),
$$ 
and consider $d={u_m}-u_l$. For $0 \leq \tau \leq s$ we have  
$$
\tlam_j({\gam})(s) - \tlam_j(\tilde{\gam})(\tau) =
\int_0^{|d|} \frac{\del}{\del u} \tlam_j (u(s) + \varsigma \frac{d}{|d|}, v(s), \sigma(s)) d\varsigma.
$$
In consequence, using the uniform bound available on the second-order derivatives of $\tlam_j$, we find 
$$
\| \tlam_j(\tilde{\gam})(\cdot) - \tlam_j({\gam})(\cdot) 
+\{ \tlam_j(\tilde{\gam})(0) - \tlam_j({\gam})(0) \} \|_{L^\infty([0,s])}
=\GO(s |d|),
$$ 
and, in consequence, 
\be
\label{PerturbationVitesse1}
\| \tilde{f}_j(\tilde{\gam})(\tau) - \tilde{f}_j({\gam})(\tau) - \tau \{ \tlam_j(\tilde{\gam})(0) - \tlam_j({\gam})(0) \}\|_{W^{1,\infty}} = \GO(s |d|).
\ee
Next, since $\| \conv f - \conv g\|_{W^{1,\infty}} \leq \| f - g \|_{W^{1,\infty}}$ we deduce that 
\be
\label{PerturbationVitesse2}
\| \conv_{[0,s]} \tilde{f}_j(\tilde{\gam}) - \conv_{[0,s]} \tilde{f}_j({\gam}) -\tau \{ \tlam_j(\tilde{\gam})(0) - \tlam_j({\gam})(0) \}\|_{W^{1,\infty}} = \GO(s|d|).
\ee 

Let us now show that
\be
\label{Eq:EvaluationNouvelleVitesse}
\| (V^\tu, \Sigma ^\tu) - (v, \sigma +  \{ \tlam_j (\tilde{\gam}(0)) - \tlam_j({\gam}(0)) \} ) \|_\infty 
= \GO(1) \, s|d|.
\ee
In view of \eqref{PropBianchiniBressan} we have
$$
\tlam_j ({\gam}(0)) = \lam_j(u_m), 
\qquad 
\tlam_j (\tilde{\gam}(0)) = \lam_j({u}_l).
$$
We denote $\tilde{\sigma}:=\sigma + \tlam_j (\tilde{\gam}(0)) - \tlam_j ({\gam_j}(0))$ and $\check{\gam} =(\tu, v,\tilde{\sigma})$. 
Using the contraction property, to establish \eqref{Eq:EvaluationNouvelleVitesse} it suffices to show
 $$
\| (v, \tilde{\sigma}) - \Omega_{s}^\tu(v, \tilde{\sigma}) \|_\infty = \GO(1) \, s|d|.
$$
Recall that $\| d\conv f - d\conv g\|_{\infty} \leq \|d f - dg \|_{\infty}$, and observe that $g - \conv g$ remains unchanged  
if one adds an affine function to $g$. Recalling also the definition of the operator $\Omega_{s}^\tu$, we 
obtain 
$$
\aligned 
& \| (v, \tilde{\sigma}) - \Omega_{s}^\tu(v, \tilde{\sigma}) \|_\infty  
\\
& = \| \tilde{f}_j(\gam) -  \conv_{[0,s]} \tilde{f}_j(\gam) - \tilde{f}_j(\check{\gam}) + \conv_{[0,s]} \tilde{f}_j(\check{\gam}) \|_\infty
 \\
& \quad +  \| \frac{d}{d\tau} \conv_{[0,s]} \tilde{f}_j(\gam) + \tlam_j (\tilde{\gam}(0)) - \tlam_j ({\gam_j}(0)) - \frac{d}{d\tau} \conv_{[0,s]} \tilde{f}_j(\check{\gam}) \|_\infty
\\
& \lesssim \| \tilde{f}_j(\gam) -  \conv_{[0,s]} \tilde{f}_j(\gam) - \tilde{f}_j(\tilde{\gam}) + \conv_{[0,s]} \tilde{f}_j(\tilde{\gam}) \|_\infty 
\\ 
& \quad + \| \tilde{f}_j(\tilde{\gam}) -  \conv_{[0,s]} \tilde{f}_j(\tilde{\gam}) - \tilde{f}_j(\check{\gam}) + \conv_{[0,s]} \tilde{f}_j(\check{\gam}) \|_\infty 
\\ 
& \quad +  \| \tlam_j (\tilde{\gam}(0)) - \tlam_j ({\gam}(0)) + \frac{d}{d\tau} \tilde{f}_j(\tau,\gamma) - \frac{d}{d\tau} \tilde{f}_j(\tau,\tilde{\gamma}) \|_\infty 
\\
& \quad +  \| \frac{d}{d\tau} \tilde{f_j}(\tau,\tilde{\gamma}) - \frac{d}{d\tau} \tilde{f}_j(\tau,\check{\gamma}) \|_\infty. 
\endaligned
$$

Taking \eqref{PerturbationVitesse1}-\eqref{PerturbationVitesse2} into account, this yields 
\begin{eqnarray*}
\| (v, \tilde{\sigma}) - \Omega_{s}^\tu(v, \tilde{\sigma}) \|_\infty
\lesssim  \| \tlam_j(\tau,\tu,v,\sigma) - \tlam_j(\tau,\tu,v,\tilde{\sigma}) \|_\infty + \GO(1) \, s|d|, \\
\end{eqnarray*}
and, using \eqref{EstimeeBianchiniBressan}, 
 \begin{eqnarray*}
\| (v, \tilde{\sigma}) - \Omega_{s}^\tu(v, \tilde{\sigma}) \|_\infty 
& \lesssim & \GO(1) \, \|v\|_\infty \|\sigma - \tilde{\sigma} \|_\infty + \GO(1) \, s |d| \\
& \lesssim & \GO(1) \, \|v\|_\infty |{u_m}-u_l | + \GO(1) \, s |d|. 
\end{eqnarray*}

Since $v = \GO(1) \, s$ (thanks to the bound on $\tlam_j$ in the domain under consideration) we deduce 
\eqref{Eq:EvaluationNouvelleVitesse}, and it follows that 
$$
\Sigma^\tu(s) - \Sigma^\tu(0) = \Sigma^{{u}}(s) - \Sigma^{{u}}(0) 
+ \GO(1) \, s |d|.
$$
In view of the proof in Bianchini \cite{Bianchini} (cf.~Lemma 3.8 therein) we have
$$ 
D(\tilde{\gam}, \gam') = \GO(s) \, |u_m -u_l|.
$$
(This follows by the same procedure, it is sufficient to prove that 
$
D(\tilde{\gam}, {\mathcal T}^i(\tilde{\gam}))$ is of order $\GO(1) \, |s| |u_m -u_l|. )
$
From Proposition~\ref{Prop:AutreOperateur} we deduce
$$
\| (V^{u'},\Sigma^{u'}) - (V^\tu,\Sigma^\tu) \|_\infty = \GO(1) \, s |d|.
$$
which yields the desired conclusion. This completes the proof of Proposition~\ref{Prop:Interaction3}
\quad\qed
\medskip

\medskip
\noindent
{\bf Proof of Corollary~\ref{CorRegTheta}.} 
In view of Proposition~\ref{Prop:Interaction3}, we see that it suffices to treat the case  $u_0 = u'_0$. 
One can estimate $\| u^1 - u^2_{|[0,s]} \|_{L^\infty}$ by
$$
\| u - u'_{|[0,s]} \|_{L^\infty} \lesssim D( {\mathcal T}^i_s (\gam'_{|[0,s]}), \gam'_{|[0,s]} ). 
$$
Recalling Lemma~\ref{LemmeEC} and the uniform $W^{2,\infty}$ estimates on $\tilde{f}_i(\gam^2)$ 
(by the regularity of $\tlam_i$ and Lemma~\ref{ConstLip}), we deduce that 
$$
\| u - u'_{|[0,s]} \|_{L^\infty} \lesssim |s'-s|.
$$
In turn, with Proposition~\ref{Prop:AutreOperateur}, this yields
$
\| \Sigma^{u}_s - \Sigma^{u'}_s \|_{L^\infty} \lesssim |s'-s|.
$
Using the uniform Lipschitz continuity estimates on $\sigma_i$, to prove \eqref{RegTheta} it suffices to check 
$$ 
\| (\Sigma^{u'}_{s'})_{|[0,s]} - \Sigma^{u'}_s \|_{L^\infty} \lesssim |s'-s|,
$$
which, again, follows from the contraction property and Lem\-ma~\ref{LemmeEC}.
\quad \qed 
\medskip


\section{A generalization to nonconservative hyperbolic systems}   
\label{NC-0}

We now turn our attention to nonlinear hyperbolic systems in nonconservative form, i.e.
\be 
\del_t u + A(u) \, \del_x u = 0, \qquad u=u(t,x) \in \RN, \quad  t \geq 0, \, x \in \RR. 
\label{NC.system}
\ee
As usual, we assume that the matrix $A(u)$ has real and distinct eigenvalues $\lam_j(u)$ and a basis 
of left- and right-eigenvectors $l_j(u), r_j(u)$, normalized as in Section~\ref{Riem-0}.  
Our aim is to extend to the nonconservative system \eqref{NC.system}
the theory of the Riemann problem discussed in previous sections. Recall that 
the distributional definition of solution does not make sense for nonconservative systems, 
and that a suitable notion of weak solution
for such systems was introduced in \cite{LeFloch0,LeFloch1,DLM,LeFloch5} which is based  on a 
{\sl prescribed family of Lipschitz continuous paths.} 
The Riemann problem for genuinely nonlinear and nonconservative systems was solved in \cite{LeFloch0,DLM}. 
We are here interested here in the generalization to systems that need not be genuinely nonlinear. 

Following LeFloch \cite{LeFloch1,LeFloch5} we can define the family of paths and, therefore, 
the generalized Rankine-Hugoniot relation for \eqref{NC.system} from the family of
traveling wave solutions associated with a regularization with small diffusion 
\be 
\del_t u + A(u) \, \del_x u = \eps \, \del_x \big (B(u) \, \del_x u\big), 
\label{NC.viscous.system}
\ee
where the diffusion matrix $B=B(u)$ is, say, positive-definite. 

We first observe that Bianchini-Bressan's arguments generalize straightforwardly to nonconservative systems. 
More precisely, the Riemann problem admits a unique entropy solution if 
the maps $\tr_j$ satisfy the conditions imposed in the beginning of Section~\ref{Sec:BB}. 
On the other hand, the interactions estimates (stated in Theorem \ref{TheoBianchini1}) were established 
in \cite{Bianchini2} {\sl in the case} that the nonconservative system is 
regularized by an identity diffusion matrix $B(u) = Id$. 
The rest of this section will thus be focused on extending  Iguchi-LeFloch's method. 

Using the notation in Section \ref{Sec:BB} introduced within Bianchini-Bressan's method, 
this amounts 
to fix vectors $\tr_j$ and to 
introduce the {\bf $j$-Rankine--Hugoniot curve} $\Hcal_j(s;u_-) = u(s,s; u_-)$ together with the 
{\bf shock speed} $\lamb_j(s,u_-)=\sigma_j(s,s; u_-)$, 
where we define the map $\tau \mapsto (u,v_j, \sigma_j) (\tau,s; u_-)$ by the differential system: 
\be
\label{NC.RHcurve} 
\aligned
& \del_\tau u (\tau) = \tr_j (u, v_j, \sigma_j)(\tau), 
\\
& \del_\tau v_j (\tau) = \tlam_j(u, v_j, \sigma_j)(\tau) - \sigma_j (0), 
\\
& \sigma_j (\tau) =  \sigma_j(0) = \frac{1}{s} \int_0^s \tlam_j(u, v_j, \sigma_j)(t) dt, 
\endaligned
\ee
with initial condition
\be
\label{NC.RHcurve.CI} 
 u (0,s) =u_-, \qquad  v_j (0,s) = 0.
\ee

\begin{theorem} 
\label{NC-2}
When the shock curves are defined by \eqref{NC.RHcurve}, 
all of the results in Section~\ref{Riem-0}, i.e.~the Riemann problem, 
regularity of the wave curves, and interaction estimates for wave strengths and inner speed variation, 
remain valid for the nonconservative system \eqref{NC.system}. 
In particular, the Hugoniot curves defined by \eqref{NC.RHcurve} satisfy 
\be  
\aligned 
& s \del_s \lamb_j(s;u_-) =  \kappa_j(s;u_-) \, \big(\lam_j(\Hcal_j(s;u_-)) - \lamb_j(s;u_-) \big),
\\
&  \del_s \Hcal_j (s, u_-) 
   =  r_j(\Hcal_j(s, u_-)) + \kappah_j(s;u_-) (\lam_j(\Hcal_j(s;u_-)) - \lamb_j(s;u_-)), 
\endaligned 
\label{NC.speedrelations}
\ee  
where the function $\kappa_j =\kappa_j(s;u_-)>0$ is smooth, bounded, and bounded away from zero, and satisfies 
$\kappa_j(\overline{s};u_-)=1$, and $\kappah_j$ is a smooth vector-valued map.   
\end{theorem} 

Observe that no assumption (of genuinely nonlinearity) need be imposed on the matrix $A$. 
Our result extends the construction of the Riemann problem given by Dal~Maso, LeFloch and Murat \cite{DLM}
for genuinely nonlinear systems. Dealing with non-genuinely nonlinear systems is more involved. 

Observe also that no reference to the vectors $\tr_j (u, v_j, \sigma_j)$ 
is made in \eqref{NC.speedrelations}, so that some further generalization of the above theorem
is possible; see Remark~\ref{plusgeneral} below.

\begin{proof} 
The analysis in Section~\ref{Riem-0} relied 
on the key property \eqref{Pro.1} 
of the shock speed and characteristic speed along the Rankine-Hugoniot curve. 
It is remarkable that this property remains
true when the shock waves for nonconservative systems are defined from traveling waves. 
Differentiating \eqref{NC.RHcurve} with respect to $s$, we obtain a linear differential system 
for the ``unknowns''  $u' := \del_s u$, $v' := \del_s v$, and $\sigma' := \del_s \sigma$:  
\begin{subequations}
\label{NC.RHcurve2}
\be
\del_\tau u'= \nabla_u \tr_j (u, v_j, \sigma_j) \cdot u'  
   				+ \del_v \tr_j (u, v_j, \sigma_j) \, v' +  \del_\sigma \tr_j (u, v_j, \sigma_j) \, \sigma', 
\label{NC.RHcurve2.a}
\ee
\be
\del_\tau v' = \nabla_u \tlam_j (u, v_j, \sigma_j) \cdot u' 
   				+ \del_v \tlam_j (u, v_j, \sigma_j) \, v'  
					+  \del_\sigma \tlam_j (u, v_j, \sigma_j) \, \sigma' - \sigma', 
\label{NC.RHcurve2.b}
\ee
\begin{multline}
\sigma' =  \frac{1}{s}(\lam_j(u(s,s;u_-)) - \sigma(s)) 
+ \frac{1}{s} \int_0^s \Big\{ \nabla_u \tlam_j (u, v_j, \sigma_j) \cdot u' (t,s)\\
   				+ \del_v \tlam_j (u, v_j, \sigma_j) \, v' (t,s)
					+  \del_\sigma \tlam_j (u, v_j, \sigma_j) \sigma'(s) \Big\} dt. 
\label{NC.RHcurve2.c}
\end{multline}
\end{subequations} 
The term $\lam_j(u(s,s;u_-))$ in \eqref{NC.RHcurve2.c} arises from the fact that
$v_j(s,s,u_-)=0$ and $\tlam_j(u,0,\cdot) = \lam_j(u)$.

In view of \eqref{NC.RHcurve.CI},  we have $u' (0,s) =0$ and $v' (0,s) = 0$. Hence, 
from \eqref{NC.RHcurve2.a}-\eqref{NC.RHcurve2.b} it follows that, given any $s$,
\be
\label{NC.maj.uv}
\| u' \|_{L^\infty(0,\tau)}  + \| v' \|_{L^\infty(0,\tau)} \lesssim \tau (\| u' \|_{L^\infty(0,\tau)}  
+ \| v' \|_{L^\infty(0,\tau)} + |\sigma'|).
\ee
Hence, using \eqref{NC.maj.uv} within \eqref{NC.RHcurve2.c} and noting that, 
by \eqref{PropBianchiniBressan}-\eqref{EstimeeBianchiniBressan}, 
we have $|\del_\sigma \tlam_j (u, v_j, \sigma_j)| <<1$, we can deduce (for sufficiently small $\tau,s$) 
\be
\label{NC.prem.estimee}
| \del_s \sigma(s)| \lesssim \left| \frac{1}{s}(\tlam_j(u(s,s;u_-),v_j(s,s;u_-),\sigma(s;u_-)) - \sigma(s)) \right|. 
\ee
In other words, we have established the first equation in \eqref{NC.speedrelations} with $\kappa_j$ {\sl bounded.} 

Considering next the map $\Hcal_j$ we see that 
\begin{eqnarray*}
\del_s \Hcal_j(s;u_-) &=& \del_\tau u(s,s;u_-) + \del_s u(s,s;u_-) 
\\
&=& \tilde{r}_j(u,v,\sigma)(s) + \del_s u(s,s;u_-).
\end{eqnarray*}
To handle the first term, we observe that $v_j(s,s;u_-)=0$ and $\tilde{r}_j(u,0;\sigma)=r_j(u)$.
For the second one, we use again \eqref{NC.maj.uv} which, together with \eqref{NC.prem.estimee},
yields the second equation in \eqref{NC.speedrelations} with $\kappah_j$ {\sl bounded.}

We next establish some regularity on $\kappa_j, \kappah_j$. 
First, from \eqref{NC.RHcurve} it follows that, for every fixed $s$, the mapping $(u,v,\sigma)$ is of class 
$C^\infty$ with respect to the variable $\tau$. 
On the other hand, it was proven in \cite{Bianchini2} (see (3.26) therein) 
that this mapping is Lipschitz continuous with respect to the variable $s$. 
Now, consider \eqref{NC.RHcurve2.a}-\eqref{NC.RHcurve2.b} as a linear differential system, with Lipschitz continuous 
coefficients, the terms cointaining $\del_s \sigma$ being viewed as a source-term. 
(Recall, moreover, that $\del_s \sigma$ is independent of $\tau$.)  
Then, it follows that $u'$ and $v'$ can be written in the form 
\be
\label{NC.uv.wrt.sigma}
(u',v')(\tau,s) =  \del_s \sigma(s) \, (h^u(\tau,s),h^v(\tau,s)),
\ee
where $h^u,h^v$ are Lipschitz continuous in both variables and $\GO(\tau)$ 
at most 
in the sup norm. Taking \eqref{NC.uv.wrt.sigma} into account into \eqref{NC.RHcurve2.c} yields 
$$
\aligned 
& \sigma' =  \frac{\lam_j(u(s,s;u_-)) - \sigma(s)}{s \, ( 1- \rho)},
\\
& \rho :=  \frac{1}{s} \int_0^s \Big( 
\nabla_u \tlam_j (u, v_j, \sigma_j) \cdot h^u(t,s) 
\\
& \hskip2.5cm 
+ \del_v \tlam_j (u, v_j, \sigma_j) \, h^v (t,s) +\del_\sigma \tlam_j (u, v_j, \sigma_j) \Big)  \, dt, 
\endaligned
$$ 
where the coefficient $\rho$ can be assumed to be sufficiently small. In consequence, $\kappa_j$ is Lipschitz continuous. 
We can derive further regularity on $(u,v)$ by \eqref{NC.RHcurve2.a}-\eqref{NC.RHcurve2.b}, 
and one concludes by a bootstrapping argument. 

Finally, the condition $\kappa_j(0,u_-)=1$ follows 
from \eqref{NC.RHcurve2.c}. Note here that the second (integral) term is of order 
$\GO(s) \, |\tlam_j(u(s,s;u_-)) - \sigma(s)|$ since $u'$, $v'$ and $\del_\sigma \tlam_j$ 
are precisely of this order for $t \in [0,s]$. 
\quad 
\end{proof} 


\begin{rema} 
\label{plusgeneral} 
1. Dal~Maso, LeFloch, and Murat's definition allows for more general jump relations that need not be related 
to a specific regularized model associated with \eqref{NC.system}. In view of the discussion above, especially 
Proposition~\ref{NC-2} it seems natural to fix a vector valued map $\kappah_j=\kappah_j(s;u_-)$
and a function $\kappa_j=\kappa_j(s;u_-)$ satisfying ($\kappa_j(0;u_-)=1$) 
and to prescribe the Hugoniot curve via the following differential 
system 
\be 
\aligned 
& \del_s \Hcal_j (s; u_-) =  r_j(\Hcal_j(s; u_-)) + \kappah_j(s;u_-) (\lam_j(\Hcal_j(s;u_-)) - \lamb_j(s;u_-)). 
\\ 
& s \, \del_s \lamb_j(s;u_-) =  \kappa_j(s;u_-) \, \big(\lam_j(\Hcal_j(s;u_-)) - \lamb_j(s;u_-) \big),
\endaligned 
\label{NC.jumprelations}
\ee  
with the initial conditions
\be
\Hcal_j (0; u_-)  = u_-, \qquad  \lamb_j(0;u_-) = \lam_j(u_-). 
\label{NC.initial}
\ee
Note the system is singular at $s=0$. It is however not difficult to check that \eqref{NC.jumprelations}-\eqref{NC.initial}
single out a unique (generalized) Hugoniot curve. 

2. To close this section, we point out that, 
more generally, it would be interesting to identify, within the DLM framework, suitable conditions on the family 
of paths ensuring that the Riemann problem admits a solution that enjoys the properties exhibited in the present paper.
Recall that, for genuinely nonlinear systems, it is known that indeed
the Riemann problem has a unique solution  \cite{LeFloch1,DLM} and that the Glimm scheme converges \cite{LeFlochLiu}. 
\end{rema}

  
\section{Existence theorem and approximation scheme} 
\label{WF-0}

\subsection{Existence result} 
In the present section we consider the Cauchy problem for a conservative or nonconservative,  
strictly hyperbolic system 
\be 
\del_t u + A(u) \, \del_x u = 0, \qquad u=u(t,x) \in \RN, 
\label{NC.system33}
\ee 
with 
\be
u(0,x) = u_0(x), \quad x \in \RR, 
\label{WF.initial}
\ee
where the initial data $u_0: \RR \to \RN$ has sufficiently small, total variation $TV(u_0)$. All wave speeds under consideration
will remain within disjoint intervals $\bigl(\lam_j^\text{min}, \lam_j^\text{max}\bigr)$. 
We establish here the convergence of a front tracking scheme for the approximation 
of  \eqref{NC.system33}-\eqref{WF.initial}. Our original contribution in this section 
is the introduction of a wave splitting strategy and the derivation of new interaction estimates on wave speeds. 
In turn, this provides us with a new existence theorem.

Following \auth{Dafermos, DiPerna, Bressan, and Risebro} we 
construct  approximate solutions to the Cauchy problem \eqref{NC.system33}-\eqref{WF.initial}, 
which are piecewise constant and consist of finitely many propagating fronts. 
In addition to the $j$-fronts  ($1 \leq j \leq N$) associated with one of the wave family of the system we will
introduce {\bf artificial fronts} of small total strength. 
More precisely, given $\varepsilon>0$ we are going to construct a piecewise constant approximate solution
 $u^\eps= u^\eps(t,x)$ satisfying the following properties:
\begin{itemize}
\item The function $u^\varepsilon(t,\cdot)$ admits a finite number of discontinuities for each time $t$, 
and the fronts meet at finitely many interaction points,
\item The propagating discontinuities in $u^\varepsilon$ are of two types:
\begin{itemize}
\item $j$-fronts $(u_-,u_+)$ associated with a family $j \in \{ 1, \ldots, N\}$ and such that $u_+ = \psi_j(s; u_-)$ for some $s$,  
and 
\item artificial fronts $(u_-,u_+)$ propagating with a fixed speed $\hlam$
larger than $\lam_N^{\text{max}}$. We sometimes refer to such a front as an $(N+1)$-front. 
No condition is imposed on the jump, and the strength of such a front is defined as 
$$
\eps_{N+1}(u_-,u_+) := |u_+ - u_-|. 
$$ 
\end{itemize}
\item The inner speed variation of or each $j$-front $(u_-,u_+)$ is uniformly small: 
$$
\vartheta_j(u_-,u_+) =  \GO(\varepsilon).
$$

\item The speed of any $j$-front $(u_-,u_+)$ is uniformly close to its correct speed, that is, calling $\lam$ the speed of the 
front, 
$$
\lam = \lamb_j(0; s; u_-) + \GO(\varepsilon).
$$ 

\item The total strength of artificial fronts is uniformly small: 
$$
\sum_{\text{artificial fronts}} |u_+ - u_-| = \GO(\varepsilon). 
$$ 
\end{itemize}

We will refer to a sequence of functions $u^\eps=u^\eps(t,x)$ satisfying the above properties 
as an {\bf $\eps$-approximate front tracking solution.}

\begin{theorem} Consider a general, nonlinear, 
strictly hyperbolic system in conservative or nonconservative form \eqref{NC.system33}. 
Then there exists a constant ${c>0}$ such that for every initial data of bounded
variation $u_0$ with ${TV(u_0) < c}$ there exists a sequence of 
$\eps$-approximate front tracking solutions which, as $\eps \to 0$, converges to the entropy solution 
of \eqref{NC.system33}-\eqref{WF.initial}. 
\end{theorem}


\subsection{Exact and approximate Riemann solvers} 

We summarize the properties of the Riemann solver that were established earlier via two different techniques 
(Iguchi-LeFloch's explicit construction, Bianchini-Bressan's vanishing viscosity approach). 
The {\bf Riemann solver} is an application which associates, to any two states 
$u_l,u_r \in \Bone$ (for some small $\delta_1>0$), the self-similar solution $u=u(\xi)$ 
($\xi =x/t$) of the corresponding Riemann problem. 
We distinguish between the three sets of conditions which we now describe. 

\subsubsection*{Wave curves.}  
\begin{itemize}
 
\item The Riemann solution is determined 
from $N$ {\sl wave curves} $\psi_j=\psi_j(s;u)$ ($1 \leq j \leq N$),  
which are globally Lipschitz continuous and locally differentiable at the origin $s=0$, with 
$$
\psi_j(0;u) =u, \qquad \del_s \psi_j (0;u) = r_j(u). 
$$

\item It consists of $N+1$ states $u_j$ separated by $N$ wave fans $(u_j,u_{j+1})$,
determined by 
$$
\aligned 
& u_r = \psi_N(s_N; \psi_{N-1}(s_{N-1}, \cdots, \psi_1(s_1;u_l) \cdots )), 
\\
& u_0 := u_l, \quad \ldots, \quad u_{j+1} :=  \psi_j(s_j; u_j), \quad \ldots, \quad u_N := u_r. 
\endaligned
$$  
The parameter values $s_j =: \eps_j(u_j, u_{j+1}) = \eps_{j}(u_l,u_r)$ are refered to as the {\bf wave strengths,} 
and $\sbold =(s_1, \cdots, s_N)$ is refered to as the {\bf strength vector.} We also use the notation 
$$
\Psi(\sbold; u_l) := \psi_N(s_N; \psi_{N-1}(s_{N-1}, \cdots, \psi_1(s_1;u_l) \cdots )).
$$

\item Each $j$-wave packet  $(u_j, u_{j+1}) = (u_j, \psi_j(s_j;u_j))$ 
consists of (finitely or countably) many shock and rarefaction waves, defined as follows.  Suppose for definiteness that $s_j>0$. 
The range of the Riemann solution is a subset of the wave curve $s \in [0,s_j] \mapsto \psi_j(s;u_j)$, 
determined by a continuous, non-decreasing  {\bf wave speed} 
$$
s \in [0,s_j] \mapsto \lamb_j(s,s_j;u_j) \in [\lam^\text{min}_j, \lam^\text{max}_j]. 
$$
Its generalized inverse $\xi \mapsto \big( \lamb_j\big)^{-1}(\xi,s_j; u_j)$ is a (possibly discontinuous) function with 
bounded variation (which could be normalized to be right- or left-continuous). 
The Riemann solution is given by   
\be
\label{WF.PbRie}
u(\xi) = \begin{cases} 
u_j,                  &    \xi \leq \lamb_j(0,s_j;u_j), 
\\
\psi_j(s;u_j),     &     s = \lamb_j^{-1}(\xi,s_j; u_j),  
\\
u_{j+1} = \psi_j(s_j;u_j),      &    \xi \geq \lamb_j(s_j,s_j;u_j).
\end{cases}
\ee

\item Furthermore, the following {\bf superposition/decomposition property} holds: 
given any speed $\xi \in \bigl( \lamb_j(0,s_j;u_j), \lamb_j(s_j,s_j;u_j))$ and introducing 
$s' := \lamb_j^{-1}(\xi,s_j; u_j)$, 
the wave packet $(u_j, \psi_j(s_j; u_j))$ can be obtained by simply 
patching together the wave packet  $(u_j, \psi_j(s'; u_j))$ and the wave packet $(\psi_j(s'; u_j), \psi_j(s_j; u_j))$.   
\end{itemize}


\subsubsection*{Interaction estimates on wave strengths.} 
To any two wave packets $(u_-, u_+)$ with $u_+= \psi_j(s;u_-)$ and $(v_-,v_+)$ with $v_+ =\psi_k(r;v_-)$, 
we can associate the {\bf potential of interaction} and the {\bf amount of interaction} 
$$
Q(u_-,u_+;v_-,v_+) \geq 0, \ \ I(u_-,u_+;v_-,v_+) \geq 0. 
$$ 
In turn, 
the potential of interaction between two Riemann problems $(u_l, u_r)=(u_0, \cdots, u_N)$ 
and $(v_l,v_r)=(v_0, \cdots, v_N)$ is given by 
\begin{gather*}
Q(u_l,u_r; v_l,v_r) := \sum_{i \geq j} Q(u_{i-1},u_i; v_{j-1}, v_j), \\
I(u_l,u_r; v_l,v_r) := \sum_{i \geq j} I(u_{i-1},u_i; v_{j-1}, v_j),
\end{gather*}
 and the interaction potential for a piecewise constant function is obtained by summing $Q(u_l,u_r; v_l,v_r)$ over all the discontinuities. \par
The functions $Q$ and $I$ should satisfy the following properties: 
\begin{enumerate}
\item 
$$
Q(u_{i},u_{i+1}; v_{j}, v_{j+1}) = | \eps_i(u_{i},u_{i+1})| \, | \eps_j(u_{j},u_{j+1})|, \qquad i \neq j, 
$$ 
$$
Q(u_{i},u_{i+1}; v_{j}, v_{j+1}) \lesssim | \eps_i(u_{i},u_{i+1})| \, | \eps_j(u_{j},u_{j+1})|,
$$
$$
Q(u_{0},\psi_i(s,u_{0}); v_{l}, v_{r}) = Q(u'_{0},\psi_i(s',u'_{0}); v_{l}, v_{r}) + \GO(1)(|u_0 - u'_0| + |s-s'|),
$$
$$
\eps_j(u_l,u_r) = \eps_j(u_l,u_m) + \eps_j(u_m,u_r) + \GO(1) \, I_{lmr},
$$  
where $I_{lmr}:=I(u_l,u_m; u_m, u_r)$.
 
\item For $u$ a piecewise constant function $u$ with small total variation, we have for some $c>0$:
$$
Q(u') \leq Q(u) - c \, I_{lmr}, 
$$
where $u'$ is obtained by replacing the Riemann problems $(u_l,u_m)$, $(u_m,u_r)$ with the Riemann problem $(u_l,u_r)$.

\item $Q$ has the following behavior with respect to a decomposition of waves: if the Riemann problem $(u_0,u_2)$ can be decomposed into 
$(u_0,u_1)$ and $(u_1,u_2)$ according to the property above, then the interaction functional remains unchanged 
when splitting $(u_0,u_2)$ as follows: 
\begin{eqnarray*}
Q(u_0,u_2; v_l, v_r) &=& Q(u_0,u_1; v_l, v_r) + Q(u_1,u_2; v_l, v_r),
\\
Q(u_0,u_2; u_0, u_2) &=& Q(u_0,u_1; u_0, u_1) + Q(u_1,u_2; u_1, u_2) \\
&\ & + Q(u_0,u_1; u_1, u_2) + Q(u_1,u_2; u_0, u_1) . 
\end{eqnarray*}
\end{enumerate}

\noindent
Since our construction requires the introduction of artificial fronts, it is necessary to specify how the interaction amount 
$Q(u_-,u_+; v_-,v_+)$ is extended to artificial fronts. {\sl When at least one of the two fronts is artificial,}
we  define
$$
\aligned 
& Q(u_-,u_+; v_-,v_+) 
\\
& := 
\begin{cases}
|u_+ - u_-| \, \eps_j(v_-,v_+),       & \text{artificial front } (u_-,u_+) \text{ and $j$-front } (v_-,v_+), 
\\
0,   & (v_-,v_+) \text{ artificial.}  
\end{cases} 
\endaligned 
$$
This definition is natural since artificial fronts meet all $j$-waves on their right-hand side, 
but do not meet other artificial fronts nor waves on their left-hand side. 


\subsubsection*{Interaction estimates on wave speeds.}  
To each wave packet $(u_{j}, u_{j+1}) = (u_{j}, \psi_j(s_j;u_{j}))$ 
we associate its minimum and maximum speeds 
$$
\lamb_j^\text{min}(u_j, u_{j+1}) := \lamb_j(0,s_j;u_j), 
\qquad 
\lamb_j^\text{max}(u_j, u_{j+1}) := \lamb_j(s_j,s_j; u_j), 
$$
respectively, as well as its {\it inner speed variation}
$$ 
\vartheta_j(u_j,u_{j+1}) := \lam_j^\text{max}(u_j, u_{j+1})  - \lam_j^\text{min}(u_j, u_{j+1}). 
$$ 
More generally, when $(u_j,u_{j+1})$ is the $j$-wave packet in the Riemann solution $(u_l,u_r)$, 
we use the notation 
$$
\aligned 
& \lamb_j^\text{min}(u_l, u_r) := \lamb_j^\text{min}(u_j,u_{j+1}), 
\qquad 
\lamb_j^\text{max}(u_l, u_r) := \lamb_j^\text{max}(u_j,u_{j+1}), 
\\
& \vartheta_j(u_l,u_r) := \vartheta_j(u_j,u_{j+1}).  
\endaligned 
$$

\begin{itemize}
\item When $u_r=\psi_j(s;u_l)$ and $u'_r=\psi_j(s';u'_l)$, 
\be
\label{IWSVForce}
| \vartheta(u_l,u_r) - \vartheta(u'_l,u'_r) | \lesssim |s'-s| + |u'_l - u_l |, 
\ee
so that, in particular, $\vartheta(u,\psi_i(s,u)) \lesssim |s|$.
 
\item When $(u_m,u_r)$ is a $j$-wave packet with $u_r:=\psi_j(s;u_m)$, 
interacting with a Riemann solution $(u_l,u_m)$  
then, setting $u'_m:=\psi_j(s;u_l)$, 
\be  
\vartheta_j(u_l,u'_m) \leq \vartheta_j(u_m,u_r)  + \GO(1) \, |u_l - u_m| \, |s|. 
\label{WF.interactangle3}
\ee 
The same inequality holds if the $j$-wave is located on the right. 
 
\item When $u_l,u_m,u_r$ is an interaction between two wave packets 
(of the same wave family $j$) only: 
 
\begin{itemize}
\item If $\eps_j(u_l, u_m) \, \eps_j(u_m, u_r) \geq 0$, then 
\be  
\aligned 
\vartheta_j(u_l,u_r) & \leq   \max \bigl( \vartheta_j(u_l,u_m), \vartheta_j(u_m,u_r) \bigr) \\
& \! \! \! \! \!  + \bigl( \lamb_j^\text{min}(u_m,u_r) - \lamb_j^\text{min}(u_l,u_m) \bigr)_+  
   + \GO(1) \, I(u_l,u_m; u_m,u_r).
\endaligned 
\label{WF.interactangle1}
\ee 
 
\item If $\eps_j(u_l, u_m) \, \eps_j(u_m, u_r) \leq 0$, then 
\be  
\aligned 
\vartheta_j(u_l,u_r) \leq & \max\bigl( \vartheta_j(u_l,u_m), \vartheta_j(u_m,u_r) \bigr) \\
 & + \GO(1) \, \min(|u_m-u_l|,|u_r -  u_m|).
\label{WF.interactangle2}
\endaligned 
\ee
\end{itemize}
\end{itemize}

\begin{rema}
The Riemann solvers described in Sections~\ref{Riem-0} and \ref{Sec:BB} satisfy the above conditions 
with, in the first case, $I=Q$ and, in the second case, distinct values $\QBB$ and $\IcalBB$.
\end{rema}


\subsubsection*{Approximate Riemann solvers.} 
We now introduce several approximate solvers which will be needed at each wave interaction points. 
The approximate Riemann solutions are piecewise constant functions which we construct in two steps. 
First, we introduce the intermediate states, that is, the $u$-components of the waves, and, second, 
we specify the speeds of each waves. 
Note that the solutions contain $j$-waves propagating with a speed close to one of the characteristic speeds
of the system, as well as artificial waves propagating at the (large, constant) speed $\hlam$.
First of all, we define the intermediate states, as follows. 

\begin{itemize} 
\item {\bf Accurate solver.} The intermediate states here are determined straightforwardly from 
the ones in the exact Riemann solution $(u_l,u_r)$. 

\item {\bf Approximate $ij$-solver} for different families $j>i$. Given a $j$-front and an $i$-front,  
$$
u_m = \psi_j(s_1; u_l), \quad  u_r = \psi_i(s_2;u_m),
$$
the approximate solver consists of the $i$-wave $(u_l,\tildeu_1)$, the $j$-wave $(\tildeu_1,\tildeu_2)$, 
and the artificial front $(\tildeu_2,u_r)$ determined by 
$$
\tildeu_0 := u_l, \quad  \tildeu_1 := \psi_i(s_2;u_l), 
\qquad 
\tildeu_2 := \psi_j(s_1; \tildeu_1), \quad \tildeu_3 = u_r.
$$
 
\item {\bf Approximate $ii$-solver} for a single family $i$. Given three states $u_l, u_m, u_r$ connected by $i$-fronts, 
$$
u_m = \psi_i(s_1; u_l), \quad  u_r = \psi_i(s_2;u_m),
$$
the approximate Riemann solution consists of the $i$-wave $(u_l,\tildeu_r)$ and the artificial front $(\tu_r,u_r)$
determined by 
$$
\tildeu_r := \psi_i(s_1+ s_2;u_l).  
$$
 
\item {\bf Artificial wave solver.} When the left-hand front is an artificial front and the right-hand one is an $i$-front
satisfying
$$
u_r = \psi_i(s;u_m),
$$
we introduce $\tildeu_m:=\psi_i(s;u_l)$ and the approximate Riemann solution consists of the $i$-wave 
$(u_l,\tildeu_m)$ and the artificial front $(\tildeu_m,u_r)$.
\end{itemize}

The second part of our construction consists of prescribing the speeds of propagation of the waves. 
Fronts having a large inner speed variation will need to be split in two or more fronts. 
Fix a threshold $\varepsilon >0$.  First of all, note that artificial fronts are never split and always travel at the speed $\hlam$. 
To handle $j$-wave packets $(u_-, u_+) = (u_-, \psi_j(s_j; u_-))$, we distinguish between two strategies: 
\begin{itemize} 
 
\item{\bf No-splitting strategy.} In this case, we propagate $(u_-, u_+)$ as a single front traveling at the smallest speed of the 
associated wave fan. 
  
\item {\bf Splitting strategy.} In this case, we replace the front by several smaller fronts defined as follows. 
Using the wave speed function $\lamb_j(s') := \lamb_j(s', s_j; u_-)$ we introduce 
$$
\aligned 
& P : = \lfloor \frac{\lamb_j(s_j) - \lamb_j(0)}{\varepsilon} \rfloor +1, 
\\
& \mu_p := \lamb_j(0) + \frac{p}{P} [\lamb_j(s_j) - \lamb_j(0)], \qquad p=0, \ldots, P-1, 
\endaligned 
$$
where $ \lfloor a  \rfloor$ denotes the greatest integer less than or equal to $a \geq 0$.  
Then, along the wave curve from $u_-$, we pick up the states associated with the speeds $\mu_p$, that is
$$
\aligned 
& w_p := \psi_j(s_p'; u_-)
\\
& s_p' := \begin{cases}
\min \big\{ s' \in [0,s_j], \ \,  \lamb_j(s') = \mu_p \big\},    &   s_j \geq 0,
\\
\max \big\{ s' \in [s_j,0], \  \, \lamb_j(s') = \mu_p \big\},    &   s_j < 0.
\end{cases}  
\endaligned 
$$

Then the $i$-wave packet is approximated by $P$ fronts with small strength, as follows:
\be
u(t,x) = 
\begin{cases} 
u_-,    &    x/t < \mu_0, 
\\
w_p,    &    \mu_{p-1} < x/t < \mu_{p}, \quad  (p=1, \ldots, P-1), 
\\
u_+,    &  x/t > \mu_{P-1}, 
\end{cases}
\nonumber 
\ee 
where, for simplicity in the notation, the Riemann solution has been centered at the point $(t,x)=(0,0)$. 
Note that each front corresponds to an exact solution of the system of conservation laws
(in the sense that the right-hand state lies on the wave curve issuing from the left-hand state)
although its propagation speed (in general) only approaches the true speed. 
\end{itemize}


\subsection{Front tracking approximations}

We are now in a position to describe our algorithm, and we fix $\varepsilon, \delta>0$ satisfying 
$$  
\delta \ll \varepsilon^2.
$$ 
Let  $u_0^\varepsilon$ be a piecewise constant approximation of $u_0$ such that
\begin{gather}
\label{CIApprox}
\| u_0^\varepsilon -u_0 \|_{L^1} \leq \varepsilon, \qquad TV(u_0^\varepsilon)  \leq TV(u_0), 
\\
\label{CIApprox2}
u_0^\varepsilon \text{ contains } 1/\varepsilon \text{ discontinuity points, at most. }
\end{gather}

At each discontinuity point of $u_0^\varepsilon$ we use the accurate solver described earlier 
together with the splitting strategy  
and we define the approximate solution $u^\eps=u^\eps(t,x)$ locally in time. The solution is extended until two fronts meet. 
As usual, the speeds of the fronts may need to be modified (cf.~Remark~\ref{RemWF1} below)
in order to ensure that only two fronts meet at every interaction. 
To extend the approximate solutions, we distinguish between several types of interactions
and we use either the accurate or the approximate solver together with the splitting or the no-splitting strategy. 

As a rule, only fronts with inner speed variation less than or equal to $\varepsilon$ will be generated by the scheme.
However, after one or several interactions, the inner speed variation may have increased so much that 
it is greater than $2 \varepsilon$ and need to be split in two fronts. 

\begin{itemize} 

\item{\bf Large interactions.} We call large interaction an interaction that involves 
an $i$-front $(u_l,u_m)$ and a $j$-front $(u_m,u_r)$ such that  
$$ 
I(u_l,u_m; u_m,u_r) \geq \delta.
$$
The solution is defined beyond the interaction time by using the accurate Riemann solver, together with the following 
rule for the speeds of the fronts : 
 
\begin{itemize}
\item For the outgoing $i$- or $j$-waves that have an inner speed variation $\vartheta \leq 2 \varepsilon$, 
we use the no-splitting strategy. 
\item For all other outgoing waves we use the splitting strategy. 
\end{itemize}

\item{\bf Small interactions.} We call small interaction an interaction that involve
an $i$-front $(u_l,u_m)$ 	and a $j$-front $(u_m,u_r)$ such that
$$ 
I(u_l,u_m; u_m, u_r) < \delta.
$$
The solution is defined beyond the interaction time by using the approximate $ij$- or $ii$-solvers. The 
speeds of the fronts are determined exactly as in the case of a large interaction (recall that artificial fronts are never split).

\item{\bf Artificial wave interactions.} At an interaction between a left-hand artificial front
and a right-hand $j$-wave, we use the artificial interaction solver, 
with again the same splitting/no-splitting strategy on the $j$-th outgoing wave.
\end{itemize}

\begin{rema} 
\label{RemWF1}
The speeds of certain fronts may need to be slightly modified in order to avoid
both interaction points involving more than two fronts  
and interaction times involving more than one interaction. 
Starting from $t=0$ or a time $t$ at which an interaction has occured, we let the fronts evolve at the speed determined earlier.
We consider the first time $\tau > t$ where two fronts meet. If there are only two fronts involved, there is nothing to do. 
If there are three (or more) fronts interacting at the time $\tau$ (possibly at different locations), 
we consider the most-left fronts $\alpha$ and $\beta$ (with $\alpha$ on the left of $\beta$). 
Now, if $\alpha$ is a $j$-front ($1 \leq j \leq N$), 
we increase the speed of $\alpha$ by an amount less than $\varepsilon$, so
that the new interaction time is beyond the time $t$. 
If $\alpha$ is an artificial front, then $\beta$ must be a $j$-front, and we diminish the speed of $\beta$ 
by an amount less than $\varepsilon$ so that again that the new interaction time is beyond time $t$. 
By this procedure, in both cases we have only decreased the interaction time. Hence, as $\tau$ was the first interaction time 
greater than $t$, the new interaction time $\tau' \in (t, \tau)$ corresponds to a single interaction of two fronts. 
An important consequence of our construction 
is that {\it two fronts that meet in the process would have met regardless of the speed modification.}  
\end{rema} \par

Based on the properties of the Riemann solver listed at the beginning of the present section, 
it is a standard matter to derive uniform estimates for the following Glimm-type functionals: 
$$
V(t) := \sum_\alpha |\eps_\alpha|, 
\qquad
Q(t) := \sum_{\alpha,\beta} Q_{\alpha\beta},
$$  
where the summations are over all fronts $\alpha, \beta$ in the approximate solution at the time $t$. 
Here, $\eps_\alpha$ denotes the strength of the front $\alpha$ and $Q_{\alpha \beta}$ is the interaction 
amount between two fronts $\alpha,\beta$.

\begin{proposition}
\label{WF-Glimm}
There exist constants $\nu,c, C_0>0$ so that for every initial data satisfying $TV(u_0)<\nu$ and for every sequence $u_0^\varepsilon$
satisfying \eqref{CIApprox}-\eqref{CIApprox2}, (as long as 
the number of waves and the number of interaction remain finite)  the quantity 
$(V + C_0 \, Q)(t)$ is non-increasing in $t$ and, more precisely, at every interaction involving
some fronts $(u_l,u_m)$, $(u_m,u_r)$
\be 
\label{DecrGlimm}
\aligned 
& \big[ V + C_0 \, Q\big] (t) 
\\
& := (V + C_0 \, Q)(t+) - (V + C_0 \, Q)(t-) 
\\
& \leq  
{\begin{cases}
-c  \, \min(|u_m-u_l|, |u_r-u_m|),         & \text{ non-monotone interactions,} 
\\
- c  \, I(u_l,u_m; u_m,u_r),          & \text{ otherwise.} 
\end{cases}
}
\endaligned 
\ee 
\end{proposition}


\subsection{Total number of fronts and interaction points}
\label{Subsec:TNF}

In view of Proposition~\ref{WF-Glimm}, 
to establish that the algorithm is well-behaved we need only show  
that the numbers of fronts and interaction points remain finite.  Define 
$$ 
\aligned 
& \Ncal(t) := \sum_{k=1}^N \Ncal_k(t), 
& \Ncal_j(t) := \# \big\{ \text{ $j$-fronts at time } t \big\}, 
\endaligned 
$$  
$$
\Theta_\varepsilon(t) := \sum_{k=1}^N  \Theta_{\varepsilon,k}(t),   
\qquad 
\Theta_{\varepsilon,j}(t) := \sum_{\text{$j$-fronts}} 
       \frac{\big( \vartheta_\alpha - \varepsilon \big)_+}{\varepsilon},
$$
and introduce the functional 
$$ \Fcal(t):= C_1 \, (V(t) + C_0 \, Q(t)) + 3 \, \Theta_\varepsilon + \Ncal(t),
$$  
where $C_1$ is some positive constant and $C_0$ is the constant already introduced in \eqref{DecrGlimm}.

\begin{lemma}
\label{PropNGF}
For $\delta$ sufficiently small and $C_1$ large enough (depending on $\varepsilon$ and $\delta$) 
the functional $\Fcal$ is non-increasing.
\end{lemma}

\begin{proof}  The function $\Fcal$ is constant away from interaction times; we distinguish between every type of interactions. 
Each interaction involves a left-hand, $i$-front $\alpha=(u_l,u_m)$ with strength $s_\alpha$ ($1 \leq i \leq N+1$) 
and a right-hand, $j$-front $\beta=(u_m,u_r)$ with strength $s_\beta$ ($1 \leq j \leq N$).
We denote by $\theta_\alpha, \lamb_\alpha$ and $\theta_\beta, \lamb_\beta$ the inner speed variation 
and wave speed function associated with the incoming waves.   
Call $\eta_k$ the outgoing $k$-waves ($1 \leq k \leq N+1$), with inner speed variation $\vartheta_{\eta_k}$ and wave speed function $\lamb_{\eta_k}$. 
Finally, call ${\mathcal I}$ the following modified amount of wave interaction:

\begin{itemize}
\item ${\mathcal I} := I$ (the interaction amount at the point), 
unless the interaction is a non-monotone interaction of waves of the same family,
\item ${\mathcal I} =$ the strength of the smallest incoming wave, 
if the interaction is a non-monotone interaction of waves of the same family.
\end{itemize}

We will summarize the evolution of the terms arising in ${\mathcal F}$ in a table (see Table~\ref{Table.1}), 
according to the nature of the interaction, the family of the wave considered, and whether the splitting strategy 
is applied or not. In this table,  ``$i \not =j$'' (respectively ``$i=j$'') refers 
to interactions between waves of different (resp.~same) family; 
``acc'' (respectively ``app.") stands for ``accurate" (resp.~``approximate"). 
We regard artificial interactions as approximate interactions between waves of different families. 
Since artificial fronts are not considered in the last two terms of ${\mathcal F}$, they are not included in the table. 

We fill up the cells of the table by relying on our interaction estimates on the strengths and the inner speed variation, 
and the definition of the number of fronts when the splitting strategy is applied. 
The following remarks are in order. \par
{1.} In the case of a monotone interaction between two fronts $(u_l,u_m)$ and $(u_m,u_r)$, 
the term $\bigl( \lamb_i^\text{min}(u_m,u_r) - \lamb_i^\text{min}(u_l,u_m) \bigr)_+$ in \eqref{WF.interactangle1} vanishes. 
This is due to the fact that fronts evolve at the minimum speed of the wave packet -- except in the case described in Remark~\ref{RemWF1}.
In the latter case, however, two fronts that meet would have met regardless of the slight change in speed.

{2.} When the $i$-th or $j$-th outgoing wave is split, we have (for instance for the $i$-wave $\alpha$):
\begin{itemize}
\item  In the case of a monotone interaction of waves of the same family, 
$$
\aligned 
\big[ \Ncal_i(t) \big]  
& = \lfloor \frac{\vartheta_{\eta_i}}{\varepsilon} \rfloor -1
\\
&  \leq -1 + \frac{\max(\vartheta_{\alpha},\vartheta_{\beta})}{\varepsilon} + \OO(\frac{1}{\varepsilon}) {\mathcal I}(t), 
\endaligned 
$$
\item In the case of a non-monotone interaction of waves of the same family, 
$$
\aligned 
\big[ \Ncal_i(t) \big]  
& = \lfloor \frac{\vartheta_{\eta_i}}{\varepsilon} \rfloor -1 
\\
& \leq -1 + \frac{\vartheta_{\alpha}}{\varepsilon} + \OO(\frac{1}{\varepsilon}) {\mathcal I}(t),
\endaligned 
$$
\item Otherwise, 
$$
\big[ \Ncal_i(t) \big] = \lfloor \frac{\vartheta_{\eta_i}}{\varepsilon} 
\rfloor  \leq  \frac{\vartheta_{\alpha}}{\varepsilon} + \OO(\frac{1}{\varepsilon}) {\mathcal I}(t).
$$
\end{itemize}
Moreover, we have 
\be
\big[ \Theta_{\varepsilon,i}(t) \big]
\leq 
\begin{cases} 
- \Big( \max(\vartheta_{\alpha},\vartheta_{\beta}) - \varepsilon\Big) / \varepsilon,    &  \text{ non-monotone case,} 
\\
- \Big( \vartheta_{\alpha} - \varepsilon\Big) / \varepsilon,         & \text{ otherwise.} 
\end{cases} 
\label{DiffGT1}
\ee 
The fact that the wave was split implies that $\vartheta_{\eta_i} \geq 2 \varepsilon$. From the estimate 
\be
\vartheta_{\eta_i} \leq 
\begin{cases} 
\max(\vartheta_{\alpha},\vartheta_{\beta}) +\GO(1){\mathcal I}(t),      &  \text{ non-monotone case},
\\
\vartheta_{\alpha} +\GO(1){\mathcal I}(t),       & \text{ otherwise,} 
\end{cases}
\nonumber
\ee 
we deduce, in both cases,
$$
\big[ \Theta_{\varepsilon,i}(t) \big] \leq -1 + \GO(\frac{1}{\varepsilon}) {\mathcal I}(t).
$$
Using again \eqref{DiffGT1}, this yields
$$
3 \big[ \Theta_{\varepsilon,i}(t) \big] \leq -1 + \GO(\frac{1}{\varepsilon}) {\mathcal I} -\begin{cases} 
\max(\vartheta_{\alpha},\vartheta_{\beta})/\varepsilon ,      &  \text{ non-monotone case}, \\
\vartheta_{\alpha}/\varepsilon ,& \text{ otherwise.} 
\end{cases}
$$

\noindent
\begin{table}[!ht]
\begin{tabular}{|c|c|c|c|c|}
\hline Interaction & Wave $k$ & Strategy  & $[{\mathcal N_k}]$ & $3 [{\Theta_{\varepsilon,k}}] \espacevertical$ \\
\hline
%
%
\multirow{3}{*}{\espaceverticalhautb $i \not = j$, acc.} & \multirow{2}{*}{\espaceverticalhaut $i$ or $j$} 
& splitting & $\leq \frac{\theta_\alpha}{\varepsilon} + \GO(1/\varepsilon) {\mathcal I}$ & 
$\begin{array}{l}
\espaceverticalhaut \leq -\frac{\theta_\alpha}{\varepsilon} -1 \\
\espacevertical \ \ \ \ + \GO(1/\varepsilon) {\mathcal I}
\end{array}$ \\ \cline{3-5}
 &  &no splitting& $0$ & $\espacevertical \leq \GO(1/\varepsilon) {\mathcal I}$ \\ \cline{2-5}
 & $k \not = i,j$ & splitting & $\espacevertical 1 + \GO(1/\varepsilon) {\mathcal I}$ & $0$ \\ \hline
%
%
%
\multirow{3}{*}{\espaceverticalhautb $i = j$, acc.} & \multirow{2}{*}{\espaceverticalhaut $i$} 
&
splitting 
& 
$\begin{array}{l}
\espaceverticalhaut \leq \frac{\max(\vartheta_{\alpha},\vartheta_{\beta}) }{\varepsilon} \\ 
\espacevertical \ \ -1+ \GO(1/\varepsilon){\mathcal I}
\end{array}$
& 
$\begin{array}{l}
\espaceverticalhaut \leq -\frac{\max(\vartheta_{\alpha},\vartheta_{\beta})}{\varepsilon} \\ 
\espacevertical \ \ \ \ -1 + \GO(1/\varepsilon) {\mathcal I} 
\end{array}$ \\ \cline{3-5}
 &  &no splitting& $-1$ & $\espacevertical \leq \GO(1/\varepsilon) {\mathcal I}$ \\ \cline{2-5}
 & $k \not = i$ & splitting & $\espacevertical 1 + \GO(1/\varepsilon) {\mathcal I}$ & $0$ \\ \hline
%
%
%
\multirow{3}{*}{\espaceverticalhautb $i \not = j$, app.} & \multirow{2}{*}{\espaceverticalhaut $i$ or $j$} 
& splitting & $\leq \frac{\theta_\alpha}{\varepsilon} + \GO(1/\varepsilon)$ & 
$\begin{array}{l}
\espaceverticalhaut \leq -\frac{\theta_\alpha}{\varepsilon} -1 \\
\espacevertical \ \ \ \  + \GO(1/\varepsilon) {\mathcal I}
\end{array}$ \\ \cline{3-5}
 & &no splitting& $0$ & $\espacevertical \leq \GO(1/\varepsilon) {\mathcal I}$ \\ \cline{2-5}
& $k \not = i,j$ & splitting & $0$ & $0 \espacevertical$ \\ \hline
%
%
%
\multirow{3}{*}{\espaceverticalhautb $i = j$, app.} & \multirow{2}{*}{\espaceverticalhaut $i$}
&
splitting
&
$\begin{array}{l}
\espaceverticalhaut \leq \frac{\max(\vartheta_{\alpha},\vartheta_{\beta}) }{\varepsilon} \\ 
\espacevertical \ \ -1+ \GO(1/\varepsilon){\mathcal I}
\end{array}$
&
$\begin{array}{l}
\espaceverticalhaut  \leq -\frac{\max(\vartheta_{\alpha},\vartheta_{\beta}) }{\varepsilon} \\ 
\espacevertical \ \ \ \ -1 + \GO(1/\varepsilon) {\mathcal I}
\end{array}$ \\ \cline{3-5}
 &  &no splitting& $-1$ & $\espacevertical \leq \GO(1/\varepsilon) {\mathcal I}$ \\ \cline{2-5}
& $k \not = i$ & splitting & $0$ & $0 \espacevertical$ \\ \hline
\end{tabular}
\caption{Estimates across each type of interaction}
\label{Table.1}
\end{table}

Now, the desired conclusion follows from the inequalities summarized in Table~\ref{Table.1}. 
At each interaction, we have $\big[V + C_0 \, Q \big](t) \leq -(C_0/2) \, {\mathcal I}(t)$, as follows from 
Proposition~\ref{WF-Glimm}. Hence, it is sufficient to take $C_1=\GO(1/\varepsilon)$, {\sl except for} 
the $N-2$ (or $N-1$) contributions of $1$ in ${\mathcal N}(t)$ due to new fronts of families $k \not = i,j$ 
in the accurate solver. But, such fronts appear only with the accurate solver; hence, when ${\mathcal I}(t) \geq \delta$, 
it suffices to take $C_1=\GO(1/\delta)$ to cover all cases.\quad 
\end{proof}


\subsubsection*{Number of fronts and interaction points.} 
\begin{itemize}
\item {\it Finite number of $j$-fronts ($1 \leq j \leq N$).} 
It follows immediately from Lemma~\ref{PropNGF} that the total number of $j$-fronts is bounded for all times. 
Moreover, observe that at all interactions generating new $j$-fronts 
an amount $\geq \min( 1, C_0 \, C_1 \delta/2)$ is used out of $\Fcal(t)$; this implies 
$$ 
\Ncal(t) \lesssim \frac{1}{\delta}.
$$ 
\item {\it Number of accurate interaction points.} From Proposition~\ref{WF-Glimm} 
it follows also that the number of accurate interactions is finite.

\item {\it Number of interaction points using the splitting strategy.} 
At such points, as seen in the table, ${\mathcal  F}$ decreases of an amount of $1$ at least;
 hence, this case arises finitely many times.

\item {\it Number of approximate $i$-interaction points using the no-splitting strategy.} 
This number is also finite since a $i$-front is lost at these points.

\item {\it Number of approximate $ij$ interaction points using the no-splitting strategy.} 
When tracing forward a given front by taking any front of the same family as its successor, we see that two fronts 
of differents families that have met will not meet again, since the speeds of different families are separated.
Since interaction points from which at least three $j$-fronts leave, are in finite number, 
this yields that approximate $ij$ points and artificial interaction points are finite in number.

\item {\it Number of artificial fronts.} 
There are only a finite number of approximate interactions and, therefore, 
the total number of artificial fronts is finite. 

\item{\it Number of artificial interaction points using the no-splitting strategy.} The same consideration as for the approximate $ij$ interaction points using the no-splitting strategy applies.
\end{itemize}


\subsection{Number of fronts of early generations.} 
We already have a bound on the total number of fronts, but 
it will be crucial in the next paragraph to have a sharp bound on their number, 
if we restrict to fronts of early generations, as we describe now.

\subsubsection*{Front generations.} We define the generation $g_\alpha$ of a front $\alpha$ as follows.

\begin{itemize}
\item[-] All fronts outgoing from the initial line have generation $1$.

\item[-] At an $ij$-interaction (an $i$-front $\alpha$ meeting a $j$-front $\beta$), the generation of outgoing $i$-fronts  
(resp. $j$, $k \not = i,j$) is fixed as $g_\alpha$ (resp. $g_\beta$, $\max(g_\alpha, g_\beta) +1$), 
whether or not these fronts are split.

\item[-] At an $i$-interaction (of two $i$-fronts $\alpha$ and $\beta$), 
the generation of the outgoing $i$-fronts (resp. $k \not = i$) is fixed as $\min(g_\alpha, g_\beta)$ 
(resp. $\max(g_\alpha, g_\beta )+1$), whether or not these fronts are split.

\item[-] At an artificial interaction point the generation of each family is preserved.
\end{itemize}


\subsubsection*{Main estimate.} We define 
\begin{gather*} 
\ONB^k(\tau):= \# \big\{ \text{ fronts at time } \tau, \text{ of generation } k \big\}, 
\\
{\mathcal P}^k := \# \big\{ (\alpha,\beta) \text{ fronts at time } 
\tau  /   \max(g_\alpha,g_\beta) = k, \ \alpha \text{ and } \beta \text{ approaching} \big\}, 
\\
\Theta_\varepsilon^{\leqslant k} := \sum_{\substack{\alpha \text{ front of} \\ 
\text{generation } \leq k}} \Big( \frac{\vartheta_\alpha - \varepsilon}{\varepsilon} \Big)_+,
\ \ \Theta_\varepsilon^{k} := \Theta_\varepsilon^{\leqslant k} -  \Theta_\varepsilon^{\leqslant k-1}.
\end{gather*}
Two fronts are said to be approaching if they are of the same family, or else if the front of the largest family is on the left. 
Hence, ${\mathcal P}^k$ represents the number of potential interactions between fronts, with largest  
generation number equal to $k$. 
Note that artificial fronts are taken into account in the above quantities (except for what concerns the inner speed variation).
Since there is no front of generation $0$, all these quantities vanish for $k=0$.

Estimates on the number of fronts of generation $\leq k$ are now derived. 

\begin{lemma}
\label{PropNbFrontsGen}
There exist a constant $C_2 = C_2(\varepsilon) >0$,
an increasing sequence $A_k = A_k(\varepsilon) >0$, and two decreasing sequences $B_k = B_k(\varepsilon) >0$ 
and $D_k = D_k(\varepsilon) >0$, such that for all $k \geq 0$ the functional
\be
\label{DefF}
{\mathcal F}_k := C_2 A_k (V+C_0 Q) + 3 {A_k} \Theta_\varepsilon^{\leqslant k} 
+ \sum_{i=1}^k B_i \ONB^i + \sum_{i=1}^{k-1} D_i {\mathcal P}^i 
\ee 
 is non-increasing in time. 
This is true regardless of the value of $\delta$ (small enough with respect to $\varepsilon$). 
Moreover, for each $k$, $A_k(\varepsilon)$ is bounded above by a polynomial expression in $1/\varepsilon$, while 
$B_k(\varepsilon)$ and $D_k(\varepsilon)$ are bounded below by polynomial expressions in $\varepsilon$.
\end{lemma}

\begin{corollary}
\label{CorNbFrontsGen}
For each $k$, the total number of fronts of generation less than or equal to $k$ is 
\be
\label{NbFtK}
\sum_{i=1}^k \ONB^i \leq J_k(\varepsilon),
\ee
where $J_k(\varepsilon)$ is polynomial in $1/\varepsilon$. 
\end{corollary}

\subsubsection*{Proof of Lemma \ref{PropNbFrontsGen} and Corollary \ref{CorNbFrontsGen}.} 
We start with two observations.
First, from the proof in Section~\ref{Subsec:TNF} that for some $C_2=C_2(\varepsilon)=\GO(1/\varepsilon)$, the functional
$$
{\mathcal G}_k:= C_2 (V+C_0 \, Q) + 3 \Theta_\varepsilon^{\leq k},
$$
is non-increasing in time. Hence, if the result in the lemma is established 
for a given $A_k$, then it remains true with a larger $A_k$. 
Moreover, as seen in 
Lemma~\ref{PropNGF}, ${\mathcal G}_k$ decreases by $1$, at least, when an incoming wave is split.

Next, when an interaction creates new fronts, the total number of new fronts is at most $\GO(1/\varepsilon)$, as follows from 
\eqref{IWSVForce} and the fact that the total strength of the waves is bounded.

Now, the proof is done by induction on $k$. The case where $k=1$ is essentially given by Lemma~\ref{PropNGF}. 
In fact, when the fronts of generation $g \geq 2$ have not to be considered, the estimate on $C_1$ can be replaced by $C_2=O(1/\varepsilon)$. 
(Note that there is no artificial front of generation $1$.) We can choose for instance $D_1=1$ and $B_0=1$. 

From now on, let us only consider the passage from $k$ to $k+1$, which is the heart of the proof. 
Suppose that for a certain $k \geq 1$, the expression in \eqref{DefF} is non-increasing. 
Whenever the property stated in Lemma~\ref{PropNbFrontsGen} is established at the rank $k$, Corollary \ref{CorNbFrontsGen} at the rank $k$
 follows immediately; hence, we can use \eqref{NbFtK} at the rank $k$.

The goal is to determine $A_{k+1}$ (large enough) and $B_{k+1}$, $D_k$ (small enough) so that the 
desired property is valid. 
We distinguish between several  cases when one of the ``new'' terms (that appear in ${\mathcal F}_{k+1}$ 
but did not appear in ${\mathcal F}_{k}$), that is, $\ONB^{k+1}$ and ${\mathcal P}^k$, can grow. 
This will provide us with some conditions on  
$A_{k+1}$, $B_{k+1}$ and $D_k$, that we will be able to fulfill by choosing these constants
sufficiently large or small.

Consider the interaction of two fronts with corresponding amount of interaction $I$. 
We introduce the same notation ${\mathcal I}$ as previously (including the case of artificial interaction), that is,
${\mathcal I}  :=Q$ (the potential at the interaction) 
unless the interaction is a non-monotone interaction  of waves of the same family, 
in which case ${\mathcal I}$ is the strength of the smallest wave.
There are two cases where ${\mathcal P}^k$ can grow, and two cases where $\ONB^{k+1}$ can grow.

\

\noindent
{\bf 1. Cases that increase ${\mathbf{\mathcal P}^k}$.}

The value ${{\mathcal P}^k}$ can grow in two ways: (1) 
an existing front of generation $ \leq k$ is split after an interaction with another wave, or (2) 
a new front of generation $ \leq k$ is created after an accurate interaction where the maximal generation of the incoming waves is exactly $k-1$. 
Note that we can measure the increase of ${\mathcal P}^k$ by
$$
[ {\mathcal P}^k ] \leq  (\text{number of new fronts }) \times {J}_k(\varepsilon).
$$

\begin{itemize}
\item {\it First case: splitting of a front of generation ${\mathit k}$.} 
In that case, the number of new fronts is of order $\GO(1/\varepsilon)$; we can compensate the increase of  
${{\mathcal P}^k}$ by taking $A_k \gtrsim J_{k}(\varepsilon)/\varepsilon$, since the decrease of ${\mathcal G}_k$ in that case is at least $1$.

\item {\it Second case: creation of a new front generation ${\mathit k}$.} We suppose that no incoming wave is split (if not, enhancing 
$A_k$ allows us to absorb this increase too.) Then, ${\mathcal P}^{k-1}$ decreases by $1$, at least, and it suffices 
to consider $D_{k} \lesssim \varepsilon D_{k-1} /J_{k}(\varepsilon)$.
\end{itemize}

\

\noindent
{\bf 2. Cases that increase $\mathbf{ \ONB^{k+1}}$.}

The value $\ONB^{k+1}$ can grow in two ways: (1) 
an existing front of generation $k+1$ is split after an interaction with another wave, or 
(2) a new front of generation $k+1$ is created after an accurate interaction where the maximal generation 
of the incoming waves is exactly $k$.

\begin{itemize}
\item {\it First case: splitting of a front of generation ${\mathit k+1}$.} In that case, the increase of the term $\ONB^{k+1}$ 
due to these new fronts is compensated by the decrease of $A_{k+1} {\mathcal G}_{k+1}$ (at least $A_{k+1}$), 
provided $A_{k+1} > B_{k+1}$ as seen in Table~\ref{Table.1}.

\item {\it Second case: creation of a new front generation ${\mathit k+1}$.} 
Again we consider only the case where there is no splitting of incoming waves. 
The increase of the term $\ONB^{k+1}$ due to these new fronts is of order $N- 2 +\GO(1/\varepsilon) {\mathcal I}$ 
(or $N- 1 +\GO(1/\varepsilon) {\mathcal I}$). In that case, ${\mathcal N}^{k}$ (if incoming waves are of same family) 
or  ${\mathcal P}^{k}$ (otherwise) decrease by $1$. Hence, 
taking $B_{k+1}$ small with respect to $B_k$ and $D_k$ allows to get the decrease. 
\end{itemize}

\noindent
This concludes the proof of Lemma~\ref{PropNbFrontsGen} and Corollary \ref{CorNbFrontsGen}. \quad \qed


\subsection{Conclusion}
The fact that the scheme converges to the entropy solution in the limit is a consequence of 
the following two properties:
\begin{itemize}
\item The $j$-fronts ($1 \leq j \leq N$) travel approximately at the correct speed given by the Rankine-Hugoniot relation.
\item The total strength of artificial fronts remains uniformly small.
\end{itemize}
Both properties are now discussed.

\subsubsection*{Accuracy of the speed of $j$-fronts.}  
From our construction it follows that, for any $j$-front $\alpha$,  
$$
\vartheta_\alpha \leq 2 \varepsilon.
$$ 
Indeed, we have: 
\begin{itemize}
\item When the front under consideration is generated by the accurate solver (or by the initial solver), then 
$\vartheta_\alpha \leq \varepsilon$. 
\item When the front is generated by an approximate solver, then either 
$\vartheta_\alpha \leq \varepsilon$ (if the solver happens to split the wave of the family of $\alpha$) 
or 
$\vartheta_\alpha \leq 2 \varepsilon$ (otherwise). 
\end{itemize}
Furthermore, the front travels at a speed which is the lowest speed in the wave packet up to a $\varepsilon$
error at most.
So, calling $\lam$ the speed of the front and $\sigma^\alpha(s)$ any speed in the wave packet, we 
also obtain 
$$ 
| \lam - \sigma^\alpha(s) | \leq 3 \varepsilon.
$$

\subsubsection*{Total strength of artificial fronts.}  
We follow here the argument known in the genuinely non-linear/linearly degenerate case 
(see, for instance, \cite{Bressan2}, Section~7.3.6). 

\ \\
{\it Strength of an artificial front.} From our previous discussion it follows that, for any artificial front $\alpha$, 
$$ 
| \eps_{N+1}(\alpha) | = \GO(1) \, \delta.
$$ 
Indeed, with the previous convention and thanks to the interaction estimates, the strength of any new artificial front 
is of order $\delta$. One can follow it during successive interactions.  
Calling $V_\alpha$ the total strength of fronts approaching $\alpha$ (that is, all 
$j$-fronts on its right), one gets (\cite{Bressan2}, p.~139) 
$$ 
|S_\alpha(t)| \leq \GO(1) \, \delta \exp( C' (V_\alpha + C_0 Q)).
$$ 
\ \\ 
{\it Total strength of artificial fronts.} For $k \geq 1$, we call $V_k$ (resp. $V^{art}_k$) the total strength of fronts
 (resp. of artificial fronts) of generation $\geq k$, and $Q_k$ defined as previously, where the sum is over
 all couple of fronts $\alpha$ and $\beta$ for which $\max(g_\alpha,g_\beta) \geq k$. \par
 
The estimates on the strengths of waves according to their generation remain valid;
hence, provided the total variation is small enough, 
we can deduce  
that, for some $\gam < 1$ and for all times $t$,
\be
\label{EstGenerations}
Q_k(t) \leq C_3 \gam^k, \qquad \qquad V_k(t) \leq C_4 \gam^k.
\ee

Now, given any integer $M$ we have 
\begin{eqnarray*}
V^{art}(t)= \sum_{k \leq M} V^{art}_k + \sum_{k > M} V^{art}_k.  
\end{eqnarray*}
The first term is estimated by
\begin{eqnarray*}
 \sum_{k \leq M} V^{art}_k \leq \GO(1) \, \delta J_M(\varepsilon),
\end{eqnarray*}
while, for the second one, 
\begin{eqnarray*}
 \sum_{k > M} V^{art}_k \leq C_4 \frac{\gam^{M+1}}{1- \gam}.
\end{eqnarray*}
Hence, given $\varepsilon>0$, we can choose $M$ so large that the second term above is less than $\varepsilon/2$, and 
next we choose $\delta$ small enough so that the first term is less than $\varepsilon/2$ too; 
hence
\be
\label{VartPetit}
V^{art}(t) \leq \varepsilon.  
\ee


\subsubsection*{Convergence.} 
Using the $L_t^\infty (BV_x)$ bound and the uniform bound 
on the wave speeds, it is a standard matter to derive a $\Lip_t(L^1_x)$ bound. 
Relying on Helly's theorem, these estimates allow us to extract a converging subsequence, say
$$ 
u^\eps \rightarrow u \in L^1_{loc}(\RN).
$$
We now check that $u$ is an entropy solution. We will consider here the case of a conservative system 
endowed with a convex entropy pair 
and prove that for any non-negative test-function $\varphi : \RR_+ \times \RR \to \RR$
and for every smooth, convex entropy / entropy-flux pair $(\eta,q) : \RN \to \RR \times \RN$.
\be
\label{EntrAprouver}
\liminf_{n \rightarrow +\infty} \int_{t=0}^\infty \int_{x \in \RR} 
\big( \eta(u^\eps) \, \varphi_t + q(u^\eps) \, \varphi_x \big) \, dtdx \geq 0.
\ee

Calling $I_n$ the above integral, for $\supp(\varphi) \subset [0,T] \times \RR$ we have 
$$ 
I_n = \int_0^T \sum_\alpha ( \dot{x}_\alpha(t) \, [\eta(u^\eps)]_\alpha(t) - [q(u^\eps)]_\alpha(t)) \varphi(t,x_\alpha) \, dt,
$$
where $x_\alpha$ is the trajectory of the front $\alpha$ and $[h]_\alpha$ is the jump of a function $h$ 
on the front $\alpha$. 
Thanks to the regularity of $\eta$ and (\ref{VartPetit}), the sum over artificial fronts is $\GO(\varepsilon)$. 
The statement \eqref{EntrAprouver} follows immediately from:   

\begin{lemma}
\label{LemErreurEntropie}
For any $j$-front $\alpha$ and all times $T_1$ and $T_2$ one has
\be
\label{ErreurEntropie}
\aligned 
& \int_{T_1}^{T_2} ( \dot{x}_\alpha(t) \, [\eta(u^\eps)]_\alpha(t) - [q(u^\eps)]_\alpha(t)) \varphi(t,x_\alpha) dt 
\\
& \gtrsim - \varepsilon |s_\alpha| |T_2 - T_1| \, \| \varphi \|_{C^1}.
\endaligned 
\ee
\end{lemma}

\medskip
\noindent
{\bf Proof.} 
Consider a front $\alpha$ connecting $u_-$ to $u_+=\psi_i(s_{\alpha},u_-)$, and denote by $\lamb_i(\cdot,s_\alpha,u_-)$ the 
corresponding wave speed function. Recall that the corresponding solution is given in \eqref{WF.PbRie}. 
Denote by $\tilde{\omega}$ the corresponding Riemann solution centered at the point $(T_1,x_\alpha(T_1))$. 
The solution satisfies  
\be
\label{ErreurEntropie2}
\int_{T_1}^{T_2} \int_{\RR} \left( \eta(\tilde{\omega}) \varphi_t + q(\tilde{\omega}) \varphi_x \right) \, dtdx
- \left[ \int_{\RR} \eta(\tilde{\omega}) \varphi(\cdot,x) \, dx \right]_{T_1}^{T_2}
\geq 0. 
\ee
We also introduce, for $t \in [T_1,T_2]$, 
$$
\omega_\alpha(t,x) 
= 
\begin{cases} 
u_-,     & x/t \leq \lam_\alpha,
\\
u_+,     &  x/t \geq \lam_\alpha,
\end{cases} 
$$
where $\lam_\alpha$ is the speed of the front $\alpha$ (that is, $\lam_\alpha=\lamb_i(0,s_\alpha,u_-)$ up to the additional 
change of speeds in Remark~\ref{RemWF1}).
Clearly the left-hand side of (\ref{ErreurEntropie}) is equal to the left-hand side of (\ref{ErreurEntropie2}) when 
we replace $\tilde{\omega}$ by $\omega_\alpha$. Hence, it is sufficient to prove that each term in the the left-hand side of
 (\ref{ErreurEntropie2})  yields an error of order $\varepsilon |s_\alpha| |T_2 - T_1| \,  \| \varphi \|_{C^1}$ 
when one replaces $\tilde{\omega}$ by $\omega_\alpha$. \par
For the first integral, the difference is supported in a triangle whose area is of order $\GO(1) \, \varepsilon (T_2 -T_1)$, 
and the difference of the integrand is of order $\GO(1) \, |s_\alpha|$. For the second integrals,  the difference is supported
 in an interval whose length is of order $\GO(1) \, \varepsilon (T_2 -T_1)$, and the difference of the integrand is again of 
 order $\GO(1) \, |s_\alpha|$, which concludes the proof of Lemma~\ref{LemErreurEntropie}. 
\quad\qed
\medskip

\section{Time-regularity of graph solutions}
\label{SY-0}

We now use the front tracking scheme to study the regularity of graph solutions introduced in \cite{LeFloch4}. 
To simplify the presentation it is 
convenient to assume that the flux $f(u)$ is defined for all $u \in \RN$. 
To begin with, we need a few definitions from \cite{LeFloch4}.
(Note that Lipschitz continuous representatives of all geometric maps under consideration
 are used throughout the present section.) 

\subsection{Geometric version of the front tracking scheme} 
A {\bf parametrized graph} is a map $(X,U) : \RR \to \RR \times \RN$
such that $X$ and $U$ are Lipschitz continuous, and 
$$
{\del_s X} \geq 0, \qquad 
\lim_{s \to \pm \infty} X(s) = \pm \infty.  
$$  
A maximal interval $[s_-,s_+]$ where $X$ remains constant will be refered to as 
a {\bf vertical segment} of the graph. We say that the parametrized graph $(X,U)$  
contains a single shock if the associated BV function $U\circ X^{-1}$
 is a step function with a single discontinuity point. 
A (Lipschitz continuous), time-dependent {\bf parametrized graph} is a map $(X,U) : \RR_+ \times \RR \to \RR \times \Bzero$
such that $X,U \in L^\infty(\RR_+, \Lip(\RR))$, $(X,U)(t)$ is a parametrized graph in the sense above, 
and for every continuously differentiable function $\theta$ with compact support and 
uniformly with respect to $t \geq 0$,  
\be 
\Big| \int_\RR U \, \del_t \theta(t,X) \, \del_s X  \, ds \Big| \lesssim \|\theta(t,\cdot)\|_{L^\infty(\RR)}. 
\label{SY.2} 
\ee  

A {\bf DLM family of paths} (Dal~Maso, LeFloch, and Murat \cite{DLM}) 
is a Lipschitz continuous map $\Phi : [0,1] \times \RN \times\RN$
such that for all $u_l,u_r,u_l',u_r' \in \RN$
\be
\begin{split} 
& \Phi(0; u_l,u_r) = u_l, \quad \Phi(1; u_l,u_r) = u_r,
\\
& \|\del_s \Phi(\cdot; u_l,u_r)\|_{L^\infty(0,1)} \lesssim |u_r - u_l|, 
\\
& \|\Phi (\cdot; u_l,u_r) - \Phi (\cdot; u_l',u_r')\|_{L^\infty(0,1)}  \lesssim |u_r - u_r'| + |u_l - u_l'|.  
\end{split}
\label{SY.3}
\ee   
A {\bf family of Riemann graphs} is a map $(\Lam, \Phi) : [0,1] \times \RN \times\RN \to \RR \times \RN$
such that :   
\begin{enumerate}
\item[$\bullet$] For any $u_l, u_r \in \RN$, $\bigl(\Lam(\cdot; u_l,u_r), \Phi(\cdot; u_l,u_r)\bigr)$ 
is a parametrized graph, and $\Phi$ is a DLM family of paths with moreover 
$$
\|\Lam (\cdot; u_l,u_r) - \Lam (\cdot; u_l',u_r')\|_{L^\infty(0,1)}  \lesssim |u_r - u_r'| + |u_l - u_l'|.   
$$  
\item[$\bullet$] The system of conservation laws \eqref{NC.system33} is satisfied: 
for every smooth function $\theta : \RR \to \RR$ 
and all $u_l,u_r \in \RN$ 
\be
\int_0^1  
\bigl( - \Lam(s; u_l,u_r) + Df \circ \Phi(s; u_l,u_r) \bigr) \, \del_s \Phi(s; u_l,u_r) \, \theta \circ \Lam(s; u_l,u_r) \, ds = 0. 
\label{SY.5}
\ee  
\item[$\bullet$] Along every vertical segment $[s_-,s_+]$,  
the Riemann graph $(\Lam,\Phi)(\cdot; u_-,u_+)$ associated with $u_\pm := \Phi(s_\pm; u_l,u_r)$
contains a single shock   
and  
\be 
\Phi(\cdot;u_l,u_r)_{| [s_-,s_+]} = \Phi(\cdot; u_-,u_+) \circ \beta, 
\label{RP.6}
\ee   
where $\beta : [s_-,s_+] \to [0,1]$ is the linear map satisfying $\beta(s_-) = 0$, 
$\beta(s_+) = 1$. 
\end{enumerate} 

\begin{defi} 
\label{SY-1}
A family of Riemann graphs $(\Lam,\Phi)$ being fixed, 
a parameterized graph $(X,U):\RR_+\times \RR \to \RR \times \RN$ is called 
a {\bf graph solution of \eqref{NC.system33} subordinate to } $(\Lam,\Phi)$ if: 
\begin{enumerate} 
\item[$\bullet$] For every test-function $\theta:\RR_+ \times \RR \to \RR$, 
\be 
\iint_{\RR_+\times\RR}
\Big( - U \, \del_t \theta(t,X) + f(U) \, \del_x \theta(t,X) \Big) \, \del_s X \, dsdt = 0.
\label{SY.7}
\ee
\item[$\bullet$] For almost every time $t$ and on every vertical segment $[s_-,s_+]$  
the graph solution coincides with a prescribed paths, more precisely 
\be 
U(t)_{\big| [s_-,s_+]} = \Phi(u_-, u_+) \circ \beta, 
\label{SY.8}
\ee 
where $u_\pm := U(t,s_\pm)$,  and $\beta : [s_-,s_+] \to [0,1]$ is the linear map satisfying $\beta(s_-) = 0$, 
$\beta(s_+) = 1$. 
\end{enumerate} 
\end{defi}

\begin{rema} 
\label{SY-2}
(i)  \,  The condition \eqref{SY.2} is equivalent to saying that the BV function $u = U \circ X^{-1}$, 
belongs to $Lip(\RR_+, L^1(\RR))$, that is 
$$
\|u(t) - u(t') \|_{L^1(\RR)} \lesssim |t - t'|.  
$$
\indent (ii) \, The condition \eqref{SY.7} is equivalent to say that $u = U \circ X^{-1}$ is a weak solution of \eqref{NC.system33} 
in the sense of distributions. 
\par 
(iii) \, In Definition~\ref{SY-1} the same parametrization (up to linear rescaling) is used for the graph solution 
and the prescribed paths. This definition will be sufficient for the purpose of the present paper. 
\end{rema}

To implement the front tracking scheme we will need to restrict the class of 
Riemann graphs, relying here on the interaction amount $Q(u_l, u_m, u_r)$ associated with three constant 
states. 

\begin{defi} 
\label{SY-3}
Consider the strictly hyperbolic systems of conservation laws \eqref{NC.system33}. 
A family of Riemann graphs $(\Lam, \Phi)$ is said to satisfy the interaction estimates if 
for all $u_l,u_m, u_r \in \RN$ such that $u_l$ is connected to $u_m$ by an $i$-wave fan 
and $u_m$ is connected to $u_r$ by an $i$-wave fan  
\be
\|(\Lam,\Phi) (\cdot; u_l,u_r) 
- (\Lam,\Phi) (\cdot; u_l,u_m) \vee (\Lam,\Phi) (\cdot; u_m,u_r)\|_{L^1(0,1)}  \lesssim 
Q(u_l, u_m, u_r), 
\label{SY.9}
\ee 
where $(\Lam,\Phi) (\cdot; u_l,u_m) \vee (\Lam,\Phi) (\cdot; u_m,u_r): [0,1] \to \RN$ 
is the arc-length parametrization of the concatenation of the two maps.  
\end{defi} 

By the result in the earlier section, there exists a family of Riemann graphs 
satisfying the interaction estimates.

\

The main result in this section is as follows. 

\begin{theorem} {\rm (Geometric version of the front tracking scheme.)} 
\label{SY-5} 
Let $(\Lam, \Phi)$ be a family of Riemann graphs satisfying the interaction estimates.   
Let $u^\eps=u^\eps(t,x)$ be a sequence of front tracking approximations for the Cauchy problem 
\eqref{NC.system33}-\eqref{WF.initial}. Then there exists a parametrized graph  
$(X^\eps, U^\eps)$ such that $X^\eps$ is Lipschitz continuous in both $(t,s)$, $U^\eps$ is continuous in the space variable $s$, 
and the following uniform estimates hold: 
\be 
\begin{split}
&  \|(\del_t X^\eps, \del_s X^\eps) \|_{L^\infty(\RR_+ \times \RR)} \leq  1, 
\qquad  
\|\del_s U^\eps \|_{L^\infty(\RR_+ \times \RR)} \lesssim TV(u_0), 
\\
& \|U^\eps(t) - U^\eps(t') \|_{L^1(\RR)} \lesssim |a^\eps(t) - a^\eps(t')|, 
\qquad t,t' \geq 0, 
\end{split}
\label{SY.10}
\ee 
where $a^\eps: \RR_+ \to \RR_+$ is a non-increasing, uniformly bounded function measuring the amount of 
cancellation and interaction in $u^\eps$. 
Moreover, $(X^\eps, U^\eps)$ converges uniformly for all but countably many times $t$ toward the graph solution 
$(X,U)$ of the Cauchy problem \eqref{NC.system33}-\eqref{WF.initial}. 

Furthermore, the map $X$ is Lipschitz continuous in both variables $(t,s)$ while $U$ is Lipschitz continuous 
in $s$ and satisfies 
\be
\|U(t) - U(t') \|_{L^1(\RR)} \lesssim |a(t) - a(t')|,  
\qquad t,t' \geq 0, 
\label{SY.11}
\ee 
where $a: \RR_+ \to \RR_+$ is the pointwise limit $a^\eps$.  
In particular, $U$ is continuous at all but countably many times $t$. 
\end{theorem}

\begin{rema} 
\label{SY-6} 
It is known that the BV solution $u$ of the problem \eqref{NC.system33}-\eqref{WF.initial}
belongs to the space $Lip(\RR_+, L^1(\RR))$, which is expressed by 
the condition \eqref{SY.5}. In contrast, the estimate 
\eqref{SY.11} provides a control of the component $U$ even in the vertical segments. 
\end{rema} 


\subsection{Proof of convergence} 

\subsubsection*{Step 1. } We begin by defining the parameterization. 
The standard approach consists of introducing the arc-length parametrization of the graph of $u^\eps$, which is based 
on the total variation map $x \mapsto TV_{-\infty}^x(u^\eps(t))$.  In order to ensure that the graph be continuous
in time and since the total variation of $u^\eps$ may change at times where two fronts meet, 
we need to modify the arc-length parametrization as follows. 
We will take into account the interaction and cancellation of waves, and take advantage of the decrease 
of the generalized Glimm functional. 
To the sequence $u^\eps$, let us associate the non-negative measure $\mu^\eps$  
consisting of Dirac masses concentrated at the points $(t_0,x_0)$ where two fronts $(u_l,u_m)$, $(u_m,u_r)$ meet. 
The point mass at $(t_0,x_0)$ is the scalar $\mu^\eps(t_0, x_0)$ given by 
\be
\mu^\eps(t_0, x_0) := - \big( Q^\eps(t_0+) - Q^\eps(t_0-) \big) \geq 0. 
\label{SY.12}
\ee 
In view of the results in previous sections the total mass of this measure is uniformly bounded
$$
|\mu^\eps|( \RR_+ \times \RR) \lesssim 1. 
$$

Let $\lam_\infty$ be a large positive constant, larger than all speeds $\lam_j$ as well as the speed $\lam_{N+1}$ 
of artificial waves introduced in the front tracking scheme. 
 Consider the map 
\be 
\sigma^\eps(t,x) := x + V^\eps (t,x) + \mu^\eps\bigl( \Omega_{t,x} \bigr), 
\label{SY.13}
\ee 
where $V^\eps$ is the total strength of waves at time $t$ on the left of $x$,  and we have introduced the triangular region 
$$
\Omega_{t,x} := \big\{ (s,y) \ \big/ \ 0 \leq s \leq t, \qquad y \leq x - \lam_\infty \, (t-s) \bigr\}. 
$$
Clearly, for every time $t$, the map $\sigma^\eps(t)$ is strictly increasing in $x$. Moreover, $\sigma^\eps$ is continuous in time 
except along polygonal lines which are transverse to the line $t=$constant. 

Let us now define $X^\eps: \RR_+ \times \RR \to \RR$ by 
\be 
X^\eps := (\sigma^\eps)^{-1}. 
\label{SY.14}
\ee 
By the properties of $\sigma^\eps$ just mentioned and the fact that $u^\eps$ is piecewise constant, it is not 
difficult to check that $X^\eps$ is Lipschitz continuous in both $t$ and $s$.

Finally, for every $t$ that is not an interaction time, following \cite{DLM}
we define $U^\eps(t)$ from $u^\eps(t)$ and $X^\eps(t)$ by completion based on the given family 
of paths $\Phi$.  
Precisely for each $t$, to the family of paths $\Phi$ and the BV function $u^\eps(t)$ one can associate 
a unique locally Lipschitz continuous, parametrized graph $(X^\eps,U^\eps)(t)$ called the $\Phi$-completion of $u^\eps(t)$, 
characterized by the conditions that  
\be 
u^\eps(t) = U^\eps \circ (X^\eps)^{-1}(t), 
\label{SY.15}
\ee 
and, on every maximal interval $[s_-,s_+]$ on which $X^\eps(t)$ is constant equal to some value $x$, 
\be 
U^\eps(t)_{\big| [s_-,s_+]} = \Phi(u^\eps(t,x-),u^\eps(t,x+)) \circ \beta, 
\label{SY.16}
\ee 
where $\beta: [s_-,s_+] \to [0,1]$ is the linear map satisfying $\beta(s_-) = 0$, 
$\beta(s_+) = 1$.

The regularity in space is clear by construction, since $u^\eps$ has uniformly bounded total variation in space
and we are using (a modification of) the arc-length parametrization. The function $\del_t X^\eps$ behaves essentially
like $\del_s X^\eps$ in a neighborhood of a single wave front. On the other hand, the function $X^\eps$ remains
constant at an interaction point. This yields easily the uniform bounds for 
$\del_s X^\eps, \del_t X^\eps, \del_s U^\eps$. 

To control $\del_t U^\eps$ we observe that 
$$
\aligned 
|U^\eps(t,s) - U^\eps(t',s)| 
& \lesssim |t - t'| \, \sup |\del_s U(t,s)| 
\\
& \lesssim |t - t'| \, \sup TV(u^\eps), 
\endaligned 
$$
provided there is no interaction point within the time interval $[t,t']$. Thus we only need 
to discuss the behavior of $\del_t U^\eps$ at interaction times.  

Consider a point $(t_0,x_0)$ of interaction of two wave fronts,
$(u_l,u_m)$ and $(u_m,u_r)$, and denote by $I_l:=[s_l,s_m]$, $I_r=[s_m,s_r]$ 
the $s$-intervals describing the incoming fronts, with 
$$
s_m - s_l = |\strength(u_l,u_m)|, \qquad s_r - s_m = |\strength(u_m,u_r)|. 
$$
In view of \eqref{SY.12} and \eqref{SY.13}, 
the interval $I=[s_l, s_r]$ is used to parameterize the outgoing wave fan. 
Let us begin with two waves of different families. 
Comparing the graphs before and after the interaction we can write 
$$
\begin{aligned} 
\| U^\eps(t_0+) - U^\eps(t_0-) \|_{L^\infty(\RR)}
& = 
\| U^\eps(t_0+) - U^\eps(t_0-)  \|_{L^\infty(I)} 
\\
& \lesssim  \min(s_m - s_l, s_r - s_m), 
\end{aligned}
$$ 
since whenever the left-hand side is Lipschitz continuous in the data $u_l$, $u_m$, $u_r$ 
and vanishes whenever the right-hand side vanishes. By integration in $s$ this
inequality leads us to the $L^1$ bound
$$
\begin{aligned} 
\| U^\eps(t_0+) - U^\eps(t_0-) \|_{L^1(\RR)}
& = \| U^\eps(t_0+) - U^\eps(t_0-) \|_{L^1(I)} 
\\
& \leq (s_r - s_l) \,   \| U^\eps(t_0+) - U^\eps(t_0-) \|_{L^\infty(I)}
\\
& \lesssim  (s_r - s_l) \, \min(s_m - s_l, s_r - s_m) 
\\
& \leq 2 \, (s_m - s_l) \, (s_r - s_m) \lesssim |u_m - u_l| \, |u_r - u_m|. 
\end{aligned} 
$$
Since the interaction potential between two waves of different families is quadratic we have obtained 
\be 
\| U^\eps(t_0+) - U^\eps(t_0-) \|_{L^1(\RR)}  \lesssim Q^\eps(t_0-) - Q^\eps(t_0+). 
\label{SY.17}
\ee

At a point where two fronts of the same family meet, we need a more precise bound 
and we rely on the property \eqref{SY.9} satisfied by the family of Riemann graphs.  
We have immediately 
\be
\aligned 
\| U^\eps(t_0+) - U^\eps(t_0-) ||_{L^1(\RR)} 
& = \| U^\eps(t_0+) - U^\eps(t_0-) ||_{L^1(I)} 
\\
& \lesssim  Q(t_0-) - Q(t_0+), 
\endaligned 
\label{SY.18}
\ee
which completes the derivation of uniform estimates.

\subsubsection*{Step 2. } We now justify the passage to the limit $\eps \to 0$.
On one hand, let us 
denote by $u$ the limit of the sequence of front tracking approximations $u^\eps$. 
On the other hand, consider the parametrized graphs $(X^\eps,U^\eps)$.
In view of the uniform Lipschitz bound on $X^\eps(t)$ (see \eqref{SY.10})  and after extracting a subsequence
if necessary we can assume that there exists a Lipschitz continuous function $X$ such that 
\be
\begin{split}
& X^\eps(t) \to X(t) \text{ uniformly on all compact subsets in } (t,s), 
\\
& \del_s X^\eps, \del_t  X^\eps \rightharpoonup \del_s X, \del_t X \text{ weak-}\!* \text{ in } L^\infty_{t,s}. 
\end{split}
\label{SY.19}
\ee 
On the other hand, since $U^\eps$ and $\del_s U^\eps$ are uniformly bounded 
we can find $U(t)$ defined for all rational numbers $t$ such that 
\be
U^\eps(t) \rightharpoonup U(t) \text{ weak-}\!* \text{ in } W^{1,\infty}_s 
\label{SY.20}
\ee
for all rational $t$ at least. Using the uniform bound on $\del_t U^\eps$ (see \eqref{SY.10}) 
it follows that $U$ can be actually defined for all times and that the convergence above holds
at all but countably many times $t$. Let us set $w := U \circ X$. 

For $t$ fixed the inverse map $(X^\eps)^{-1}(t)$ is monotone increasing and thus converges pointwise
some some limit $Y$. From the identity  $X^\eps \circ (X^\eps)^{-1}= id$ we deduce that 
$X \circ Y = id$, so that $Y$ is the 
(generalized) inverse of $X$.  In turn, passing to the limit in the relation $U^\eps \circ (X^\eps)^{-1} = u^\eps$
yields $w=U \circ X^{-1}= u$. In other words the BV function associated with the
limiting graph is exactly the function $u$. In particular, the graph satisfy the conservation law in the sense 
\eqref{SY.7}.

To conclude it remains to determine the vertical parts of $(X,U)$ and establish \eqref{SY.8}. 
We will check that $(X,U)$ coincides (up to re-parametrization) 
with the $\Phi$-completed graph of $u$ whose arc-length parameterization is denoted by 
$(Y,V_\Phi)$: 
\be 
(X,U) = (Y,V_\Phi) \circ \beta, 
\label{SY.21}
\ee
where $\beta : \RR \to \RR$ is an increasing, onto, and Lipschitz continuous map. 
In particular, \eqref{SY.21} implies that the vertical parts in both graph coincide
so that \eqref{SY.8} holds. 

The statement \eqref{SY.21} is a non-trivial property of the front tracking scheme which we will establish by relying on 
the well-known property of {\bf local uniform convergence :}   
for all but countably many times $t$ and for each $\eta>0$, any point $x \in \RR$ 
must be 
\begin{enumerate}
\item either a point of continuity of $u(t)$ for which there exist a neighborhood $N_\eta(x)$
and a real $\eps_\eta>0$ sufficiently small so that for all $y \in N_\eta(x)$ and $\eps < \eps_\eta$
\be
|u^\eps(t,y) - u(t,x)| + |u(t,y) - u(t,x)| < \eta, 
\label{SY.22}
\ee
\item or else a point of jump for $u(t)$ for which there exist a neighborhood $N_\eta(x)$,
a sufficiently small real $\eps_\eta>0$, and a sequence of points $x^\eps \in N_\eta(x)$, 
 so that for all $y \in N_\eta(x)$ and $\eps < \eps_\eta$
\be
\begin{aligned}
& |u^\eps(t,y) - u(t,x\pm)| + |u(t,y) - u(t,x\pm)| < \eta \, \text{ if } \, y \lessgtr x^\eps. 
\end{aligned}
\label{SY.23}
\ee  
\end{enumerate} 

Specifically, the argument presented now shows that the local uniform convergence of functions implies, 
for the associated $\Phi$-completed graphs, the convergence in the uniform sense of graphs.

Fix a time $t$ at which $u^\eps(t) \to u(t)$ uniformly locally. Given $\eta>0$ we can select finitely
many points $z_1, \ldots, z_m$ so that 
\be 
\sum_{\substack{x \in\RR} \\ {x \neq z_1, \ldots, z_m}} |u_+(t,x) - u_-(t,x)| < \eta.  
\label{SY.24}
\ee
To each $z_j$ we associate a neighborhood $N_\eta(z_j)$ and a point $z_j^\eps$ so that the property
\eqref{SY.23} holds at $z_j$ with $\eta$ replaced by $\eta/m$:
\be
\aligned
& |u^\eps(t,y) - u(t,z_j\pm)| + |u(t,y) - u(t,z_j\pm)| < \eta/m 
\\
& \text{ if } \, y \lessgtr z_j^\eps, \, y \in N_\eta(z_j). 
\endaligned
\label{SY.25}
\ee  

For each $R>0$ consider the compact set $K_R := [-R,R] \setminus \bigcup_j N_\eta(z_j)$. 
By \eqref{SY.22}, to each point $x \in K_R$ we can associate an open neighborhood $N_\eta(x)$ that should be sufficiently small
so it does not contain any of the $z_j$ and  
\be
|u^\eps(t,y) - u(t,x)| + |u(t,y) - u(t,x)| < \eta, \quad y \in N_\eta(x) 
\label{SY.26}
\ee
holds. By compactness we can extract finitely many points $z_{m+1}, \ldots z_p$ so that 
the whole family of $N_\eta(z_j)$ form a covering of $[-R,R]$: 
$$
[-R,R] \subset \bigcup_{j=1}^p N_\eta(z_j). 
$$  

Note that, in view of the Lipschitz continuity property of the map $\Phi$ (see \eqref{SY.3}), 
the graph distance associated with $u^\eps$ and $u$ is clearly controled in the $L^\infty$ norm. 
Such an estimate is used within each region of small oscillations, as follows. 
Set $u_j^{\pm,h} := u^\eps(t,z_j^\eps\pm)$ and $u_j^\pm := u(t,z_j\pm)$ for $j=1, \ldots, m$. 
By distinguishing between regions of large jump and regions of small oscillations 
we can estimate the graph distance between the approximate graphs $(X^\eps,U^\eps)$
and the $\Phi$-completion $(Y,V_\Phi)$, as follows (where the constants $C_1, C_2, \ldots$ depend
upon $\Lip(\Phi)$ only) 
$$
\dist  \bigl( (X^\eps,U^\eps), (Y,V_\Phi) \bigr)  \leq D_1 + D_2
$$
with, for the shock regions ($j=1, \ldots m$),
\be 
\begin{split} 
D_1 = & \, 2 \sup_{j=1, \ldots, m} \sup_{\substack{{x \in N_\eta(z_j), x \lessgtr  z_j } 
\\ 
{ y \in N_\eta(z_j^\eps), y  \lessgtr  z_j^\eps}}} | u^\eps(t,y) - u(t,x) |  
\\
& + 2 \, \Lip(\Phi) \, \sup_{j=1, \ldots, m} \sup_{x,y \in N_\eta(z_j), x,y  \lessgtr  z_j} | u(t,y) - u(t,x) | 
\\
&  + 2 \, \Lip(\Phi) \, \sup_{j=1, \ldots, m} \sup_{x,y \in N_\eta(z_j^\eps), x,y  \lessgtr  z_j^\eps} | u^\eps(t,y) - u^\eps(t,x) | 
\\
& + \sup_{j=1, \ldots, m} \sup_{s \in [0,1]} | \Phi(s; u_j^{-,h}, u_j^{+,h}) - \Phi(s; u_j^-, u_j^+) |
\end{split} 
\nonumber 
\ee 
and, for the regions with small oscillations ($j=m+1, \ldots p$), 
\be 
\begin{split} 
D_2 = & \, \sup_{j=m+1, \ldots, p} \sup_{x,y \in N_\eta(z_j)} | u^\eps(t,y) - u(t,x) |  
\\
& \, + \Lip(\Phi) \, \sup_{j=m+1, \ldots, p} \sup_{x,y \in N_\eta(z_j)} | u(t,y) - u(t,x) | + | u^\eps(t,y) - u^\eps(t,x) |. 
\end{split} 
\nonumber 
\ee 

Using again the Lipschitz continuity of $\Phi$ we can write 
$$
| \Phi(s; u_j^{-,h}, u_j^{+,h}) - \Phi(s; u_j^-, u_j^+) | \leq \Lip (\Phi) \, ( |u_j^{-,h} - u_j^-| + |u_j^{+,h} - u_j^+| ). 
$$ 
Therefore in view of \eqref{SY.25} and \eqref{SY.26} we find 
$$
\dist  \bigl( (X^\eps,U^\eps), (Y,V_\Phi) \bigr) \leq C \, \eta.
$$
where the constant $C$ only depends upon $\Lip (\Phi) $. This completes the proof of Theorem~\ref{SY-5}.


\begin{acknowledgement} 
 
This work was completed when the second author (P.G.L.) was visiting the Mittag-Leffler Institute 
(Stockholm) in 2005 during the Fall Semester Program on ``Nonlinear Wave Motion'' organized by
A. Constantin, C.M. Dafermos, H. Holden, K.H. Karlsen, and W. Strauss.

\end{acknowledgement}




\end{document}